\documentclass[12pt,a4paper]{amsart}
\usepackage{latexsym,amsfonts,amsmath,amssymb,amsthm,rotating,txfonts}
\usepackage{diagrams}
\usepackage{epic}
\usepackage{curves}

\newarrow{filto}{triangle}--->
\newarrow{Izom}C---{>>}
\newtheorem{thm}{Theorem}[section]
\newtheorem{cor}[thm]{Corollary}
\newtheorem{notation}[thm]{Notation}

\newtheorem{lemma}[thm]{Lemma}
\newtheorem{prop}[thm]{Proposition}
\theoremstyle{definition}
\newtheorem{defi}[thm]{Definition}
\newtheorem{exa}[thm]{Example}
\newtheorem{rem}[thm]{Remark}
\numberwithin{equation}{section}

\newcommand{\dif}{\mathrm{Diff}}
\newcommand{\diffdnk}{\dif_d(k)\times\dif_d(n)}
\newcommand{\res}{\operatornamewithlimits{Res}}
\newcommand{\R}{\mathbb{R}}
\newcommand{\C}{\mathbb{C}}

\renewcommand{\P}{\mathbb{P}}

\newcommand{\mapnk}{\mathcal{J}(n,k)}

\newcommand{\mapd}[2]{\mathcal{J}_{d}({#1},{#2})}
\newcommand{\mapdo}[2]{\mathcal{J}_{d}({#1},{#2})^0}

\newcommand{\hm}[2]{\mathrm{Hom}(\C^{#1},\C^{#2})}

\newcommand{\lin}{\mathrm{Lin}}

\newcommand{\mdeg}[1]{\mathrm{mdeg}[#1]}
\newcommand{\tensor}{\otimes}
\newcommand{\rd}[1]{\mathcal{J}_{d}(#1)}
\newcommand{\rdk}{\mapd kn}

\newcommand{\diff}[1]{\mathrm{Diff}_{d}({#1})}

\newcommand{\blp}{\Psi}
\newcommand{\bv}{\mathbf{v}}

\newcommand{\Thetab}{\overline{\Theta}}

\newcommand{\ddd}{,\dots,}
\newcommand{\epd}[1]{\mathrm{eP}[#1]}

\newcommand{\taugr}{{\tau_\mathrm{Gr}}}
\newcommand{\taum}{{\tau_M}}
\newcommand{\tautt}{{\tau_T}}
\newcommand{\evs}{{\ev_S}}

\newcommand{\snk}{\C[\boldsymbol{\lambda},\boldsymbol{\theta}]^{\sg
    n\times\sg k}}
\newcommand{\thetab}{\boldsymbol{\theta}}
\newcommand{\lambdab}{\boldsymbol{\lambda}}

\newcommand{\tp}[1]{\mathrm{Tp}_{#1}}
\newcommand{\gl}[1]{\mathrm{GL}_{#1}}
\newcommand{\sym}{\mathrm{Sym}}
\newcommand{\GL}{\mathrm{GL}}
\newcommand{\fl}{\mathrm{Fl}}

\renewcommand{\hom}{\mathrm{Hom}}
\renewcommand{\fl}{\mathrm{Fl}}

\newcommand{\codim}{\mathrm{codim}}
\newcommand{\fld}{\mathrm{Flag}_{d}(\C^n)}
\newcommand{\comb}{\mathrm{perm}}

\newcommand{\hnd}{\hom^\reg(\C_L^d,\C^n)}
\newcommand{\sol}{\mathrm{Sol}}
\newcommand{\grass}[2]{\mathrm{Gr}\left(#1,#2\right)}
\newcommand{\eqns}{\mathcal{E}}
\newcommand{\eqnstil}{\widetilde{\mathcal{E}}}
\newcommand{\euler}{\mathrm{Euler}}

\newcommand{\tc}{\hat T}
\newcommand{\bu}{\mathbf{u}}

\newcommand{\buref}{{\teps}_{\mathrm{ref}}}

\newcommand{\rf}{{\mathrm{ref}}}
\newcommand{\sg}[1]{\mathcal{S}_{#1}}
\newcommand{\reffl}{{{\mathbf{f}_\rf}}}
\newcommand{\bipi}{{\boldsymbol{\pi}}}
\newcommand{\birho}{{\boldsymbol{\rho}}}
\newcommand{\Bipi}{{\boldsymbol{\Pi}_d}}
\newcommand{\Bipio}{{\boldsymbol{\Pi}_\O}}
\newcommand{\ct}{\mathrm{RC}}
\newcommand{\ds}{{\mathrm{dst}}}
\newcommand{\dist}{{\mathrm{dst}}}
\renewcommand{\O}{\mathcal{O}}

\newcommand{\coeff}{\mathrm{coeff}}
\newcommand{\ires}{\res_{z_1=\infty}\res_{z_{2}=\infty}\dots\res_{z_d=\infty}}
\newcommand{\sires}{\res_{\mathbf{z}=\infty}}
\newcommand{\dbz}{\,d\mathbf{z}}
\newcommand{\bz}{\mathbf{z}}

\newcommand{\polye}{\C[u^\bullet]}

\newcommand{\mult}{\mathrm{mult}}
\newcommand{\Z}{\mathbb{Z}}
\newcommand{\defect}{\mathrm{defect}}
\newcommand{\n}{\mathfrak{n}}

\newcommand{\sign}{\mathrm{sign}}
\newcommand{\lrel}{I'_{\O}}
\newcommand{\dod}{\{1,\dots,d\}}

\newcommand{\uh}{\hat{u}}
\newcommand{\um}[1]{\mathrm{sum}(#1)}
\newcommand{\tfnk}{\mathcal{F}_d(n)}
\newcommand{\tfnkreg}{\mathcal{F}_d^{\mathrm{reg}}(n)}
\newcommand{\fnk}{\widetilde{\mathcal{F}}_d(n)}
\newcommand{\fnkreg}{\widetilde{\mathcal{F}}_d^{\mathrm{reg}}(n)}

\newcommand{\mapreg}{\mathcal{J}^\reg_d(1,n)}
\newcommand{\reg}{\mathrm{reg}}
\renewcommand{\epsilon}{\varepsilon}

\newcommand{\im}{\mathrm{im}}

\newcommand{\hofi}{\hom^\triangle(\C^d_R,\symdot)}

\newcommand{\symdot}{\sym_d^{\bullet}\C^n}
\newcommand{\newd}{\mathrm{Ym}^\bullet\C_L^d}
\newcommand{\homnewd}{\hom^\triangle(\C_R^d,\newd)}

\newcommand{\taufl}{\tau_{\mathrm{Fl}}}

\newcommand{\phie}{\phi_{\eqnstil}}
\newcommand{\feq}{\mathrm{Ind}(\eqnstil)}
\newcommand{\cu}{u_{\pi}^{l}}
\newcommand{\qfl}{Q_{\fl}}

\newcommand{\emu}{\mathrm{emult}}
\newcommand{\lie}{\mathrm{Lie}}
\newcommand{\cpt}{\mathrm{cpt}}
\newcommand{\thom}{\mathrm{Thom}(W)}
\newcommand{\mP}{\mathbb{P}}
\newcommand{\qd}[1]{\QQ_{#1}}
\newcommand{\QQ}{\widehat{Q}}
\newcommand{\OO}{\widehat{\O}}
\newcommand{\erf}{\hat{\epsilon}_\rf}
\newcommand{\prz}{\mathrm{pr}_0}
\newcommand{\KK}{\widehat{\mathcal{T}}}
\newcommand{\NN}{\mathcal{N}}

\newcommand{\tilf}{\check{f}}
\newcommand{\gnk}{\mathrm{Diff}(\C^k)\times\mathrm{Diff}(\C^n)}
\newcommand{\tpo}{\mathrm{Tp}_O}
\newcommand{\ePd}{equivariant Poincar\'e dual }
\newcommand{\DD}{\mathrm{TD}}
\newcommand{\porb}{T\!\cdot\! p}
\newcommand{\lead}{\mathrm{lead}}
\newcommand{\relz}[2]{\mathrm{Rel}(#1,#2)}
\newcommand{\trp}[2]{\langle #1\leftrightarrow #2\rangle}

\newcommand{\J}{\mathcal{J}}
\newcommand{\ev}{\mathrm{ev}}
\newcommand{\evm}{\mathrm{ev}_{\!M}}
\renewcommand{\Upsilon}{\Theta'}
\newcommand{\jddiff}{\mapreg/\diff1}
\newcommand{\phip}{\psi}
\newcommand{\phif}{\phi_{\widetilde{\mathcal{F}}}}
\newcommand{\phit}{\phi_{\widetilde{\mathcal{F}}}}
\newcommand{\teps}{{\tilde\epsilon}}
\newcommand{\solfnk}{\sol_{\widetilde{\mathcal{F}}}}
\newcommand{\eref}{\gamma_{\rf}}
\newcommand{\homreg}{\hom^{\mathrm{reg}}(\C^d_L,\C^n)}
\newcommand{\phigr}{\phi_{\mathrm{Gr}}}
\newcommand{\phic}{\alpha}
\newcommand{\pfl}{\pi_{\fl}}
\newcommand{\Gammac}{\gamma}
\newcommand{\ind}[1]{\mathrm{Ind}(#1)}
\newcommand{\kappat}{\tilde\kappa}
\newcommand{\epsref}{\epsilon_{\rf}}

\newcommand{\Nh}{\widehat{\NN}}
\newcommand{\phat}{\widehat{\mathrm{pr}}}

\newcommand{\io}{I_{\O}}
\newcommand{\ehat}{\hat{\epsilon}_\rf}
\newcommand{\TT}{\mathrm{T}}
\newcommand{\Ts}{\mathrm{Ts}}
\newcommand{\abi}{\alpha_\bipi}
\newcommand{\tepsi}{{\tilde\epsilon}_\bipi}
\newcommand{\Lie}{\mathrm{Lie}}
\newcommand{\Ss}[1]{S_{\!#1}}
\newcommand{\pre}{\mathrm{pr}_{\eqns}}

\newcommand{\tcl}{\overline{\Theta}_d}

\newcommand{\Cw}{\mathrm{CW}}
\newcommand{\qmr}{Q_d(n)}
\newcommand{\grnk}{\grass{-dk}{\mapd nk}}
\newcommand{\Lambdal}{L}

\begin{document}
\begin{flushright}
{\em In memoriam Raoul Bott}
\end{flushright}
\vspace{.5in}
\title{Thom polynomials of Morin singularities}
\author{Gergely B\'erczi}
\address{Mathematical Institute, Oxford University, 24-29 St Giles, OX1 3LB,
Oxford, UK}
\email{email: berczi@maths.ox.ac.uk}
\author{Andr\'as Szenes}
\address{Mathematics Institute, Budapest University of Technology,
   Egry J u 1, 1111 Budapest, Hungary}
\address{Universit\'e de Gen\`eve, Section de Math\'ematiques, 2-4 rue du
Li\`evre, 1211 Gen\`eve, Suisse}
\email{andras.szenes@unige.ch}
\thanks{The support of OTKA and FNS is gratefully acknolwedged}
\maketitle

\setcounter{section}{-1}
\section{Introduction}
\label{sec:intro}

We begin with a quick summary of the notions of global singularity
theory and the theory of Thom polynomials. For a more detailed review
we refer the reader to \cite{arnold,kazarian}.

Consider a holomorphic map $f:N\to K$ between two complex manifolds,
of dimensions $n\leq k$. We say that $p\in N$ is a {\em singular}
point of $f$ if the rank of the differential $df_p:\TT_pN\to
\TT_{f(p)}K$ less than $n$. 

Topology often forces $f$ to be singular at some points of $N$, and we
will be interested in studying such situations.  Before we proceed, we
introduce a finer classification of singular points. Choose local
coordinates near $p\in N$ and $f(p)\in K$, and consider the resulting
map-germ $\tilf_p:(\C^n,0)\to(\C^k,0)$, which may be thought of as a
sequence of $k$ power series in $n$ variables without constant
terms. The group of infinitesimal local coordinate changes $\gnk$ acts
on the space $\mapnk$ of all such map-germs. We will call
$\gnk$-orbits or, more generally, $\gnk$-invariant subsets
$O\subset\mapnk$ {\em singularities}.  For a singularity $O$ and
holomorphic $f:N\to K$, we can define the set
\[ Z_O[f] = \{p\in N;\; \tilf_p\in O\}, \] which is independent of any
coordinate choices.  Then, under some additional technical assumptions
(compact $N$, appropriately chosen closed $O$, and sufficiently
generic $f$), $Z_O[f]$ is an analytic subvariety of $N$.  The
computation of the Poincar\'e dual class $\alpha_O[f]\in H^*(N,\Z)$ of
this set is one of the fundamental problems of global singularity
theory. This is indeed useful: for example, if we can prove that
$\alpha_O[f]$ does not vanish, then we can guarantee that the
singularity $O$ occurs at some point of the map $f$.

This problem was first studied by Ren\'e Thom (cf. \cite{thom, haef})
in the category of smooth varieties and smooth maps; in this case
cohomology with $\Z/2\Z$-coefficients is used. Thom discovered
that to every singularity $O$ one can associate a bivariant
characteristic class $\tau_O$, which, when evaluated on the pair
$(TN,f^*TK)$ produces the Poincar\'e dual class $\alpha_O[f]$.  One of
the consequences of this result is that the class $\alpha_O[f]$
depends only on the homotopy class of $f$.

A similar result, which we will call {\em Thom's principle}, has been
used in the holomorphic category (cf. \cite{kazarian,fr} and
\S\ref{sec:epdthom} of the present paper). To formulate it in more
concrete terms, denote by $\snk$ the space of those polynomials in the
variables $(\lambda_1\ddd \lambda_n,\theta_1\ddd\theta_k)$ which are
invariant under the permutations of the $\lambda$s and the
permutations of the $\theta$s.  According to the structure theorem of
symmetric polynomials, $\snk$ itself is a polynomial ring in the
elementary symmetric polynomials:
\[   \snk = \C[c_1(\lambdab)\ddd c_n(\lambdab),
c_1(\thetab)\ddd c_k(\thetab)].
\]
Using the Chern-Weil map, a polynomial $Q\in \snk$, and a pair
of bundles $(E,F)$ over $N$ of ranks $n$ and $k$, respectively,
produces a characteristic class $Q(E,F)\in H^*(N,\C)$. Then the
complex variant of Thom's principle reads: \\ {\em For appropriate
  $\gnk$-invariant $O$ of codimension $m$ in $\mapnk$, there exists a
  homogeneous polynomial $\tpo\in\snk$ of degree $m$, such that for an
  arbitrary, sufficiently generic map $f:N\to K$, the cycle
  $Z_O[f]\subset N$ is Poincar\'e dual to the characteristic class
  $\tpo(TN,f^*TK)$.}

A  precise version of this statement is described in
\S\ref{sec:epdthom}.  The polynomial $\tpo$ is called the {\em Thom
  polynomial} of $\O$, and  the computation of these polynomials is a
central problem of singularity theory.

The structure of the $\gnk$-action on $\mapnk$ is rather complicated;
even the parametrization of the orbits is difficult.  There is,
however, a simple invariant on the space of orbits: to each map-germ
$\tilf:(\C^n,0)\to(\C^k,0)$, we can associate the finite-dimensional
nilpotent algebra $A_{\tilf}$ defined as the quotient of the algebra
of power series $\C[[x_1\ddd,x_n]]$ by the ideal generated by the
pull-back subalgebra $\tilf^*(\C[[y_1\ddd,y_k]])$. This algebra
$A_{\tilf}$ is trivial if the map-germ $\tilf$ is nonsingular, and it
does not change along a $\gnk$-orbit (cf. \S2 more details).

Combining Thom's principle with this observation,  to
each finite-dimensional nilpotent algebra $A$ and pair of integers
$(n,k)$, one can associate a doubly symmetric polynomial $\tp A^{n\to
  k}\in \snk$; in the sense described above, this will serve as a
universal Poincar\'e dual of those points in the
source spaces of holomorphic maps whose local nilpotent algebra is $A$.

The computation of Thom polynomials associated to nilpotent algebras
is a difficult problem. A few structural statements are known, however
(cf. \S\ref{sec:structural} for more details).

First, as discovered by Damon and Ronga (\cite{damon,sigma11}) in
the 70's, the polynomial $\tp A^{n\to k}$ lies in the subring of
$\snk$ generated by the relative Chern classes defined by the
generating series
\[ 1+c_1q+c_2q^2+\dots =
\frac{\prod_{j=1}^k(1+\theta_jq)}{\prod_{i=1}^n(1+\lambda_iq)}.
\]
Next, the Thom polynomial, expressed in terms of these relative Chern
classes, only depends on the codimension $j=k-n$. More precisely,
there is a unique polynomial $\DD^j_A(c_1,c_2,\dots)$ such that
\[ \tp A^{n\to k}(\boldsymbol{\lambda},\thetab)
=\DD^{k-n}_A(c_1(\boldsymbol{\lambda},\thetab),
c_2(\boldsymbol{\lambda},\thetab),\dots).
\]

Finally, in a recent paper, Feh\'er and Rim\'anyi observed \cite{fr}
that performing the substitution $c_i\mapsto c_{i-1}$ in $\DD^{j}_A$
produces $\DD^{j-1}_A$.  This implies that to each nilpotent algebra $A$
one can associate a power series in infinitely many variables, which
encodes all of the Thom polynomials associated to $A$. This
observation served as the starting point for the present work.

In this paper, we will concentrate on the so-called Morin singularities
\cite{morin}, which correspond to the situation when the algebra $A$
is generated by a single element. The list of these algebras is
simple: $A_d=t\C[t]/t^{d+1}$, $d=1,2,\dots$.

  The goal of our paper is to compute the Thom polynomial $\tp
  {A_d}^{n\to k}$ for arbitrary $d,n$ and $k$. For simplicity of notation,
  we will denote this polynomial by $\tp d^{n\to k}$, or sometimes simply
  by $\tp d$, omitting the dependence on the parameters $n$ and
  $k$.

  The problem of calculating $\tp d^{n\to k}$ goes back to Thom
  \cite{thom}.  The case $d=1$ is the classical formula of Porteous:
  $\tp 1=c_{k-n+1}$. The Thom polynomial in the $d=2$ case was
  computed by Ronga in \cite{sigma11}. More recently, in \cite{bfr},
  the authors proposed a formula for $\tp3$; P. Pragacz has given a
  sketch of a proof for this conjecture \cite{prag}. Finally, using
  his method of restriction equations, Rim\'anyi \cite{rimanyi} was
  able to treat the $n=k$ case, and computed $\tp d^{n\to n}$ for $d\le8$
  (cf. \cite{gaffney} for the case $d=4$).

  Our approach combines the test-curve model of Porteous \cite{por}
  with localization techniques in equivariant cohomology
  \cite{bv,rossmann,voj}. We obtain a formula which reduces the
  computation of $\tp d^{n\to k}$ to a certain problem of commutative
  algebra which only depends on $d$. This problem is trivial for
  $d=1,2,3$, hence we instantly recover all results known for
  arbitrary $n\leq k$. An important feature of our formula is that it
  manifestly satisfies all three properties listed above. In
  particular, we obtain a tentative geometric interpretation for the
  Thom series introduced by Feh\'er and Rim\'anyi.

  The paper is structured as follows: we describe the basic setup and
  notions of singularity theory in \S\ref{sec:basicsing}, essentially
  repeating the above construction using more formal notation. Next,
  in \S\ref{sec:epdthom} we recall the notion of equivariant
  Poincar\'e dual, which provides us with a convenient language for
  describing Thom polynomials. We also present the localization
  formulas of Berline-Vergne \cite{bv} and Rossmann \cite{rossmann},
  which are crucial to our computations. In \S\ref{sec:loctech} we
  develop a calculus, localizing equivariant Poincar\'e duals by
  combining the localization principles with Vergne's integral formula
  for equivariant Poincar\'e duals.  With these preparations, we
  proceed to describe the test curve model for Morin singularities in
  \S\ref{sec:model}. The key part of our work is \S\ref{sec:compact},
  where we reinterpret this model using a double fibration in a way
  which allows us to compactify our model space and apply the
  localization formulas. The following section, \S\ref{sec:apploc} is
  a rather straightforward application of the localization techniques
  of \S\ref{sec:epdthom} to the double fibration constructed in
  \S\ref{sec:compact}. The resulting formula \eqref{fixedtwo}, in
  principle, reduces the computation of our Thom polynomials to a
  finite problem, but this formula is difficult to use for concrete
  calculations. Remarkably, however, the formula undergoes through
  several simplifications, which we explain in
  \S\ref{sec:vanishing}. At the end of \S\ref{sec:apploc}, we
  summarize our constructions and results in a diagram, which will
  hopefully orient the reader.

  The simplifications bring us to our main result: Theorem \ref{dathm}
  and formula \eqref{thethom}. While this formula is rather simple, it
  still contains an unknown quantity: a certain homogeneous polynomial
  $\QQ_d$ in $d$ variables, which does not depend on $n$ and $k$.  The
  list of these polynomials begins as follows:
\[  \QQ_1=\QQ_2=\QQ_3=1,\;\QQ_4(z_1,z_2,z_3,z_4)=2z_1+z_2-z_4,\dots
\]
In principle, $\QQ_d$ may be calculated for each concrete $d$ using a
computer algebra program, but, at the moment, we do not have an
efficient algorithm for performing such calculations for large $d$. We
discuss certain partial results in the final section of our paper;
these, in particular, allow us to compute $\QQ_5$ by hand, and $\QQ_6$
using the computer algebra program Macaulay.  We will elaborate on
this method in a forthcoming publication.

We end the paper with an application of our theorem to positivity of
Thom series. Rim\'anyi conjectured in \cite{rimanyi} that the Thom
polynomials $\tp d$ expressed in terms of relative Chern classes have
positive coefficients. Our formalism suggests a stronger positivity
conjecture, which we formulate in \S\ref{sec:positive}, and check for
the first few values of $d$. A list of notations is provided in
\S\ref{sec:listnot} to help the reader navigate the paper.

In closing, we note that Morin singularities are special cases of the
so-called Thom-Boardman singularities
\cite{thom,boardman,mather}. These are parametrized by finite
nonincreasing sequences of integers, and Morin singularities
correspond to sequences starting with 1.  Our method extends to a
wider class of Thom-Boardman singularities; we hope to report on
new results in this direction in a later publication.

{\bf Acknowledgments}. We would like to express our gratitude to
Rich\'ard Rim\'anyi for introducing us to the subject, and explaining
this problem to us. We are greatly indebted to Mich\`ele Vergne, whose
ideas profoundly influenced this paper. In particular, most of
\S\ref{sec:loctech} is based on her suggestions. Finally, useful
discussions with L\'aszl\'o Feh\'er, Maxim Kazarian, Andr\'as
N\'emethi, Felice Ronga and Andr\'as Sz\H ucs are gratefully
acknowledged.

\section{Basic notions of singularity theory}
\label{sec:basicsing}

\subsection{The setup}\label{sec:setup}
We start with a brief introduction to singularity theory. We suggest
\cite{mather2},\cite{arnold},\cite{thom} as references for the
subject.

Let $(e_1,\dots,e_n)$ be the basis of $\C^n$, and denote the
corresponding coordinates by $(x_1,\dots,x_n)$. Introduce the
notation $\mathcal{J}(n)=\{h\in\C[[x_1\ddd x_n]];\; h(0)=0\}$ for
the algebra of power series without a constant term, and let $\rd n$
be the space of $d$-jets of holomorphic functions on $\C^n$ near the
origin, i.e. the quotient of $\mathcal{J}(n)$ by the ideal of those
power series whose lowest order term is of degree at least $d+1$.
As a linear space, $\rd n$ may be identified with polynomials on
$\C^n$ of degree at most $d$ without a constant term.

In this paper, we will call an algebra {\em nilpotent} if it is
finite-dimensional, and there exists a positive integer $N$ such
that the product of any $N$ elements of the algebra vanishes. The
algebra $\rd n$, in particular, is nilpotent, since $\rd n^{d+1}=0$.

Our basic object is $\mapd nk$, the space of $d$-jets of holomorphic
maps $(\C^n,0)\to(\C^k,0)$. This is a finite-dimensional complex
vector space, which one can identify $\rd n\tensor\C^k$; hence
$\dim\mapd nk =k\binom{n+d}{d}-k$.  We will call the elements of
$\mapd kn$ {\em map-jets of order} $d$, or simply map-jets. In this
paper we will always assume $n\leq k$.

One can compose map-jets via substitution and elimination of terms
of degree greater than $d$; this leads to the composition maps
\begin{equation}
  \label{comp}
\mapd nk \times\mapd mn\to\mapd mk,\;\;  (\Psi_2,\Psi_1)\mapsto
\Psi_2\circ\Psi_1.
\end{equation}
When $d=1$, $\mathcal{J}_1(m,n)$ may be identified with $n$-by-$m$
matrices, and \eqref{comp} reduces to multiplication of matrices.  By
taking the linear parts of jets, we obtain a map
\[
\lin:\mapd nk\to\hom(\C^n,\C^k),
\] which is compatible with the compositions \eqref{comp} and matrix
multiplication.

Consider now the  set
\[
\diff n = \{\Delta\in\mapd nn;\;\lin(\Delta) \text{ invertible}\}.
\]
The composition map \eqref{comp} endows this set with the structure
of an algebraic group, which has a faithful representation on $\rd
n$. Using the compositions \eqref{comp} again, we obtain the
so-called {\em left-right} action of the group $\diff k\times\diff
n$ on $\mapd kn$:
\[ [(\Delta_L,\Delta_R),\Psi]
 \mapsto \Delta_L\circ\Psi\circ \Delta_R^{-1}
 \] Note that the action of $\diff n$ is linear, while the action of
 $\diff k$ is not. {\em \sl Singularity theory}, in the sense that we
 are considering here, studies the left-right-invariant algebraic
 subsets of $\mapd nk$.

 A natural way to form such subsets is as follows. Observe that to
 each element $\Psi=(P_1,\dots,P_k)\in\mapd nk$, where $P_i\in\rd n$
 for $i=1\ddd k$, we can associate the quotient algebra $A_\Psi=\rd
 n/I\langle P_1,\dots,P_k\rangle$: the algebra $\rd n$ modulo the
 ideal generated by the elements of the sequence. Since $\rd
 n^{d+1}=0$, we also have $A_\Psi^{d+1}$=0. We will call $A_\Psi$
 the {\em nilpotent algebra}\footnote{Instead of this algebra, it is
   customary to use the so-called {\em local algebra} of $\Psi$, which
   is simply the augmentation of $A_\Psi$ by the constants.} of the
 map-jet $\Psi$.  For $\Psi=0$ this nilpotent algebra is $\rd n$,
 while for a generic $\Psi$ (in fact, as soon as
 $\mathrm{rank}[\lin(\Psi)]=n$) we have $A_\Psi=0$.

 Now let $A$ be a nilpotent algebra, as defined above.  Consider
 the subset
\begin{equation}
  \label{defconea}
\Theta_A^{n\to k} = \{(P_1,\dots,P_k)\in\mapd nk;\;\rd n/I\langle
P_1,\dots,P_n\rangle \cong A\}
  \end{equation}
  of the map-jets of order $d$.  Again, the dependence on the
  parameters $d,n$ and $k$ will be usually omitted. 

  It is easy to show that $\Theta_A$ is $\diff k\times\diff
  n$-invariant. A key observation is that although two map-jets with
  the same nilpotent algebra may be in different $\diff k\times\diff
  n$-orbits, there is a group acting on $\mapd nk$ whose orbits are
  exactly the sets $\Theta_A^{n\to k}$ for various nilpotent algebras
  $A$. This group is defined as the semidirect product
\begin{equation}
  \label{defk}
\mathcal{K}_d(n,k) = \GL_k(\C\oplus\rd n)\rtimes \diff n,
\end{equation}
using the natural action of $\diff n$ on $\rd n$; the algebra
$\C\oplus\rd n$ is the augmentation of $\rd n$ by constants. The
vector space $\rd n$ is naturally a module over $\C\oplus\rd n$, and
hence $\mathcal{K}_d(n,k)$ acts on $\mapd nk$ via
\begin{equation}
  \label{via}
   [(M,\Delta),\Psi]\mapsto \left(M\cdot\Psi\right)\circ\Delta^{-1},
\end{equation}
where ``$\cdot$'' stands for matrix multiplication.

\begin{prop}[\cite{mather},\cite{mather2},\cite{arnold}]
\label{korbit} Two map-jets in $\mapd
  kn$ have the same nilpotent algebra if and only if they are in the
  same $\mathcal{K}_d$-orbit.
\end{prop}
\begin{rem}
  Two jets in the same $\mathcal{K}_d$-orbit are called {\em contact
    equivalent}, or $\mathcal{K}$-equivalent (cf. \cite{arnold}). The
  term $V$-equivalence is also used (e.g. \cite{martinet}). The
  varieties $\Theta_A$ are called \emph{contact singularity classes}
  or simply \emph{contact singularities}.
\end{rem}

Using the fact that  $\mathcal{K}_d$ is connected, 
it is not difficult to derive the following properties of $\Theta_A$.
  \begin{prop}[\cite{arnold}]
    \label{propgena} Let $A$ be a nilpotent algebra such that
    $A^{d+1}=0$ and $n\geq\dim(A/A^2)$.  Then for $k$ sufficiently
    large, $\Theta_A^{n\to k}$ is a nonempty, $\diff k\times\diff
    n$-invariant, irreducible quasiprojective algebraic variety of
    codimension $(k-n+1)\dim(A)$ in $\mapd nk$.
\end{prop}
Note that the codimension of $\Theta_A$ depends only on the difference
$k-n$ and does not depend on $d$.
    
In the present paper, we will study certain  rough topological
invariants of  contact singularities; these invariants depend only
on the closure of the singularity locus in $\mapd nk$. As it turns
out, in an asymptotic sense, the closures of contact orbits are also
closures of left-right orbits, hence, from our point of view, these
two types of singularity classes are closely related.

While we will not need this statement, we describe it in some
details for reference.  Roughly, we claim that for fixed $A$ and
$r$, and sufficiently large $n$, there is a dense left-right orbit
in $\Theta_A^{n\to n+r}$.

Let $r$ be a nonnegative integer. An {\em unfolding} of a map-jet
$\Psi\in\mapd nk$ is a map-jet $\widehat\Psi\in\mapd{k+r}{n+r}$ of
the form
\[(x_1,\ldots ,x_n,y_1,\ldots, y_r)\mapsto
(F(x_1,\ldots, x_n,y_1,\ldots, y_r),y_1,\ldots, y_r)
\]
where $F\in \mapd {n+r}k$ satisfies
\[F(x_1\ldots, x_n,0,\ldots,
0)=\Psi(x_1,\ldots, x_n).\] The {\em trivial unfolding} is the
map-jet
\[(x_1,\ldots, x_n, y_1,\ldots, y_r)\to (\Psi(x_1,\ldots
,x_n),y_1, \ldots ,y_r).\]

\begin{defi}[\cite{arnold},\cite{mather2}]
  A map-jet $\Psi \in \mapd nk$ is {\em stable} if all unfoldings of
  $\Psi$ are left-right equivalent to the trivial unfolding.
\end{defi}

Informally, a germ of a holomorphic map $f:N\rightarrow K$ of
complex manifolds at a point $x\in N$ is stable if for any
small deformation $\tilde{f}$ of $f$, there is a point in the
vicinity of $x$ at which the germ of $\tilde{f}$ is left-right
equivalent to the germ of $f$ at $x$.

Now we can formulate the relationship between contact and left-right
orbits precisely.

\begin{prop}[\cite{arnold},\cite{mather2}]
\label{propstable}
  \begin{enumerate}
  \item If $\,\widehat\Psi$ is an unfolding of $\Psi$, then
    $A_{\widehat\Psi}\cong A_\Psi$.
  \item Every map germ has a stable unfolding.
  \item If a map germ is stable, then its left-right orbit is dense in
    its contact orbit.
    \end{enumerate}
\end{prop}

\subsection{Morin singularities}
\label{sec:adsing}

In this paper, we will focus on nilpotent algebras $A$ generated by a
single element. Such algebras form a one-parameter family:
\[ A_d = t\C[t]/t^{d+1},\quad d=1,2,\dots
\]
The corresponding singularity classes are called the {\em
  $A_d$-singularities} or {\em Morin singularities}
\cite{arnold},\cite{morin}. We introduce the simplified notation
\begin{equation}
  \label{thetasimp}
\Theta_d^{n\to k}  \text{ instead of } \Theta_{A_d}^{n\to k}
\end{equation}
for these varieties, and we will omit the parameters $n$ and $k$
when this causes no confusion.

Let us specialize the results quoted in  the previous paragraph to
the case of the  $A_d$ algebras. We have
\begin{itemize}
\item  $(A_d)^{d+1}=0$, hence we can work in $\mapd nk$.
\item The variety $\Theta_d^{n\to k}$ is nonempty for any $n\leq k$. For
  $n=k=1$, we simply have $\Theta_d[1,1]=\{0\}$, the constant zero germ
  in $\J_d(1,1)$.  This germ is not stable.
\item There are stable map-jets in $\mapd nk$ with nilpotent algebra
  $A_d$, whenever $n\ge d$. An example in $\mathcal{J}_N(d,d)$ for
  $N\geq d$  with minimal source dimension $n=d$ is
\begin{equation}\label{stableexample}
(x_1 \ldots, x_d)\mapsto (x_d^{d+1}+x_1x_d^{d-1}+x_2x_d^{d-2}+\ldots
+x_{d-1}x_d,x_1,\ldots,
  x_{d-1}).
\end{equation}
\end{itemize}

Finally, we recall that the $A_d$-singularities fit into the wider
family of  so-called {\em Thom-Boardman} singularity classes.
(\cite{boardman},\cite{arnold}). A Thom-Boardman class is specified
by a nonincreasing sequence of positive integers $i_1\ge \ldots\ge
i_d$; the class corresponding to the special values $i_1=\ldots
=i_d=1$ contains exactly those maps with nilpotent algebra
isomorphic to $A_d$.

As the description of $\Theta_d$ as a Thom-Boardman class is rather
different from \eqref{defconea}, we provide it  for reference.
Observe that
\begin{itemize}
\item eliminating the terms of degree $d$ results in an algebra
  homomorphism $\pi_{d\to d-1}I:\rd n\to \mathcal{J}_{d-1}(n)$, and
\item partial differentiation $f\mapsto \partial f/\partial x_j$ is a
  well-defined map $\rd n\to \mathcal{J}_{d-1}(n)$ for $j=1\ddd n$.
\end{itemize}
Now, given a proper ideal $I$ in the algebra $\rd n$, denote by
$\delta I$ the ideal in $\mathcal{J}_{d-1}(n)$ generated by $\pi_{d\to
  d-1}I$ together with the determinants of the $n$-by-$n$ matrices of
the form
\[     \det\left(\frac{\partial Q_i}{\partial x_j}\right)_{i,j=1}^n\in\mathcal{J}_{d-1}(n),
\]
 with arbitrary $Q_1,\dots,Q_{n}\in I$.

\begin{prop}\label{charlocalg} Denoting by $ I\langle
  P_1,\dots,P_k\rangle$ the ideal in $\rd n$ generated by the elements
  $P_1\ddd P_k$, we have
  \begin{multline} \Theta_d^{n\to k}=\{(P_1,\dots,P_k)\in\mapd nk;\;
  \codim(\delta^{d-1}
  I\langle P_1,\dots,P_k\rangle \subset \mathcal{J}_{1}(n))=1\}.
  \end{multline}
\end{prop}

\section{Equivariant Poincar\'e duals and Thom polynomials}
\label{sec:epdthom}

The goal of this paper is to compute certain topological invariants
of the subvarieties $\Theta_d^{n\to k}$ introduced in the previous
section. In this section, we define and describe these invariants in
detail.

Let $T$ be a complexified torus: $T\cong(\C^*)^r$.  The {\em \ePd}
is an invariant $\Sigma\mapsto\epd\Sigma$ associated to algebraic or
analytic $T$-invariant subvarieties of $T$-modules; this invariant
takes values in homogeneous polynomials on the Lie algebra
$\mathrm{Lie}(T)$ of $T$.  The central objects of the present work,
Thom polynomials, are special cases of equivariant Poincar\'e duals
(cf.  \cite{rimanyi},\cite{kazarian}). We review the definitions and
properties of equivariant Poincar\'e duals in some detail here in
order to prepare ourselves for the localization formulas of the next
section.

The \ePd has appeared in the literature in several guises: as Joseph
polynomial, equivariant multiplicity, multidegree, etc.  One of the
first definitions was given by Joseph \cite{joseph}, who introduced
it as the polynomial governing the asymptotic behavior of the
character of the algebra of functions on the subvariety. Rossmann in
\cite{rossmann} defined this invariant for analytic subvarieties via
an integral-limit representation, and then used it to write down a
very general localization formula for equivariant integrals. This
formula will play an important role in our computations.

We begin with an explicit formula in \ref{sec:epdmult}, then turn to
an axiomatic definition of the invariant in \S\ref{subsec:axiomatic},
Following the algebraic treatment of \cite{milsturm}. This will
provide us with some useful computational tools. After considering an
example in \S\ref{sec:basic}, and recording a few technical statements
in \S\ref{sec:sometech}, we turn to the analytic picture.  We first
give an overview of Rossmann's localization formula, then we describe
Vergne's integral representation, which places the \ePd in the proper
context of equivariant cohomology.  Finally, in \S\ref{sec:thompoin}
we define Thom polynomials as \ePd s, and we justify this definition;
this allows us to formulate our problem precisely. In the final
paragraph, we collect what is known about the general structure of
Thom polynomials of contact singularities.

\subsection{ Equivariant Poincar\'e duals, Multidegrees}
\label{sec:epdmult} Denote the weight lattice of $T=(\C^*)^r$ by
$\Lambdal$; this is the lattice in $\Lie(T)^*=\C^r$ generated by the
standard weights (the coordinate vectors)
$\lambda_1,\dots,\lambda_r$. Let $W$ be an $N$-dimensional complex
vector space endowed with an action of $T$. This action is
diagonalizable, hence one can choose coordinates $y_1,\dots,y_N$ on
$W$ in such a way that the action in the dual basis is diagonal;
denote the respective weights by $\eta_1\ddd\eta_N$.

Note that we will {\em not} restrict ourselves to the so-called {\em
  convergent} case (cf. \cite{rossmann,milsturm}), i.e. we will not
assume that the weights $\eta_1\ddd\eta_N$ all lie in an open
half-space of $\Lambdal\otimes_{\Z}\R\subset\Lie(T)^*$; hence the
$\Lambdal$-graded pieces of the ring $S=\C[y_1,\ldots, y_N]$ of
polynomial functions on $W$ might be infinite-dimensional.

Let $\Sigma$ be a closed $T$-invariant algebraic subvariety of $W$,
and denote by $I(\Sigma)\subset S$  the ideal of polynomials
vanishing on $\Sigma$.  This ideal is {\em reduced}, i.e. has the
property that $f^n\in I(\Sigma)\Rightarrow f\in I(\Sigma)$. Our plan
is to define an extended invariant: $I\mapsto\mdeg{I,S}$, called the
{\em multidegree} of $I$, where $I$ is an arbitrary $T$-invariant
ideal in $S=\C[y_1\ddd y_N]$. Then we can simply define the
equivariant Poincar\'e of a variety as the the multidegree of the
corresponding ideal (cf. Definition \ref{defepd} below). 
Now we sketch an explicit and an axiomatic definition of the
multidegree.

For the construction, let $D$ be the codimension of the variety
defined by the ideal $I\subset S$, and consider a finite, $T$-graded
resolution of $S/I$ by free $S$-modules:
\[
\oplus_{i=1}^{j[M]}  Sw_i[M]\to\dots\to  \oplus_{i=1}^{j[m]}
 Sw_i[m]\to\dots\to \oplus_{i=1}^{j[1]} S w_i[1] \to S \to S/I\to 0;
\]
where $w_i[m]$ is a free generator of degree $\eta_i[m]\in\Lambdal$ for
$i=1,\dots j[m]$, $m=1\ddd M$. Then
\begin{equation}
  \label{expldefmdeg}
  \mdeg{I,S} = \frac1{D!}\sum_{m=1}^M\sum_{i=1}^{j[m]} (-1)^{D-m} \eta_i[m]^D.
\end{equation}
\begin{defi}\label{defepd} Let  $\Sigma \subset W$ be $T$-invariant
 closed subvariety as in
  \S\ref{sec:epdmult}. Then we define the {\em $T$-equivariant Poincar\'e
    dual of $\Sigma$ in $W$} by
\[       \epd{\Sigma,W}_T=\mdeg{I(\Sigma),\C[y_1\ddd y_N]}.
\]
\end{defi}
We will usually omit the lower index $T$ when this does not cause
confusion.  Note that the multidegree, and hence the equivariant
Poincar\'e dual, is manifestly a homogeneous polynomial of degree
$D$.

While \eqref{expldefmdeg} is explicit, its meaning is not
transparent, and we note that, usually, it is rather difficult to
write down free resolutions of ideals. Hence we turn to an axiomatic
description, which is more intuitive, and provides us with a more
algorithmic understanding of the invariant as well.

\subsection{Axiomatic definition}
\label{subsec:axiomatic}

We follow the treatment of \cite{milsturm} to give the axiomatic
definition: we describe 3 characterizing properties of the
multidegree, and then we prove that these properties indeed
determine the polynomial.

The monomials $\mathbf{y}^{\mathbf{a}}=\prod_{i=1}^Ny_i^{a_i}\in
S=\C[y_1,\ldots, y_N]$ are parametrized by the integer vectors
$\mathbf{a}=(a_1\ddd a_N)\in\Z_+^N$. A {\em monomial order} $<$ on
$S$ is a total order of the monomials in $S$ such that for any three
monomials $m_1,m_2,n$ satisfying $m_1>m_2$, we have $nm_1>nm_2>m_2$
(see \cite[\S 15.2]{eisenbud} ).

An ordering of the coordinates $y_1,\ldots, y_N$ induces the
so-called {\em  lexicographic} monomial order of the monomials, that
is, $\mathbf{y}^{\mathbf{a}}>\mathbf{y}^{\mathbf{b}}$ if and only if
$a_i>b_i$ for the first index $i$ with $a_i \neq b_i$. We will use
this lexicographic monomial order throughout this paper.

Now let $I\subset S$ be a $T$-invariant ideal. Define the {\em initial
  ideal} $\mathrm{in}_<(I)\subset S$ to be the ideal generated by
the monomials $\left\{\mathrm{in}_<(p):p\in I \right\}$, where
$\mathrm{in}_<(p)$ is the largest monomial of $p$ w.r.t $<$. There is
a flat deformation of $I$ into $\mathrm{in}_<(I)$ (\cite{eisenbud},
Theorem 15.17.), and the first axiom says that
$\mdeg{I}$ does not change under this deformation:\\
{\bf 1. Deformation invariance}:
$\mdeg{I,S}=\mdeg{\mathrm{in}_<(I),S}$.

To describe the second axiom, we define the multiplicity of a
maximal-dimensional component of a non-reduced variety. Let $I\subset
S$ be an ideal, and denote $\Sigma(I)$ the variety of common zeros of
the polynomials in $I$:
\[   \Sigma(I) = \{p\in W;\; f(p)=0\;\forall f\in I\}.
\]
Denote by $\Sigma_1,\Sigma_2,\ldots, \Sigma_m$ the
maximal-dimensional irreducible components of $\Sigma(I)$. Then each
$\Sigma_i$ corresponds to a prime ideal $\mathfrak{p}_i\subset S$,
and one can define a positive integer $\mult(\mathfrak{p}_i,I)$, the
{\em
  multiplicity of} $\Sigma_i$ with respect to $I$, as the length of
the largest finite-length $S_{\!\mathfrak{p}_i}$-submodule in
$(S/I)_{\mathfrak{p}_i}$, where $S_{\!\mathfrak{p}_i}$ (resp.
$(S/I)_{\mathfrak{p}_i}$) is the localization of $S$ (resp. $S/I$)
at $\mathfrak{p}_i$ (see section II.3.3 in \cite{eisenbudharris}).
Then
we have \\
{\bf 2. Additivity:}
\begin{equation}
  \label{addide}
\mdeg{I,S}=\sum_{i=1}^m \mult(\mathfrak{p}_i,I) \cdot
\mdeg{\mathfrak{p}_i,S}.
\end{equation}

The last axiom describes the multidegree for the case of coordinate
subspaces:\\
 {\bf 3. Normalization:}
for every subset $\mathbf{i}\subset\{1\ddd N\}$ we have
\begin{equation}
\mathrm{mdeg}\left[\left\langle
y_{i},\,i\in\mathbf{i}\right\rangle,S\right]=
\prod_{i\in\mathbf{i}}\eta_i,
\end{equation}
where $\langle\cdot\rangle$ stands for the ideal generated by the
polynomials listed in the angle brackets.

A special case of the normalization axiom is the case
$\Sigma=\{0\}$. We will often use the notation $\euler^T(W)$ for
$\epd{\{0\},W}$, since, indeed, this is the equivariant Euler class
of $W$ thought of as a $T$-vector bundle over a point. We have thus
  \begin{equation}
    \label{eec}
\epd{\{0\},W}_T =  \euler^T(W)  =\prod_{i=1}^N\eta_i.
  \end{equation}

\begin{rem} \label{passtoeuler}

  Using this notation, the normalization axiom may be recast in a
  geometric form as follows: given a surjective equivariant linear map
  $\gamma:W\to E$ from $W$ to another $T$-module $E$, we have
\begin{equation}
  \label{compint}
  \epd{\gamma^{-1}(0),W}=\euler^T(E).
\end{equation}
\end{rem}

 Consider the following three examples:
\begin{enumerate}
\item Set $N=4$, and consider the ideal $I=\langle y_1^2,y_2^3,y_3\rangle$ in
  $S=\C[y_1,y_2,y_3,y_4]$. This is the line $\{y_1=y_2=y_3=0\}$ with
  multiplicity $6$, so its multidegree is
\begin{equation*}
\mdeg{I,S}=6\eta_1\eta_2\eta_3.
\end{equation*}
\item The ideal $I=\langle y_1^2y_2^3y_3\rangle$ in
  $S=\C[y_1,y_2,y_3]$ corresponds to the union of the hyperplanes
  $y_1=0,y_2=0,y_3=0$ with multiplicities $2,3,1$, respectively. By
  the normalization and additivity properties
\begin{equation*}
\mdeg{I,S}=2\eta_1+3\eta_2+\eta_3
\end{equation*}
\item The ideal $I=\langle y_1y_2,y_2y_3,y_1y_3 \rangle = \langle y_1,y_2
\rangle \cap \langle y_2,y_3 \rangle \cap \langle y_1,y_3 \rangle $
in $S=\C[y_1,y_2,y_3]$ has three components with multiplicity 1,
corresponding to the given decomposition, so
\[\mdeg{I,S}=\eta_1\eta_2+\eta_2\eta_3+\eta_1\eta_3\]
\end{enumerate}

Following \cite{milsturm} \S8.5, now we sketch an algorithm for
computing $\mdeg{I,S}$, proving that the axioms determine this
invariant.

An ideal $M\subset S$ generated by a set of monomials in
$y_1,\ldots, y_N$ is called a \emph{monomial ideal}. Since
$\mathrm{in}_<(I)$ is such an ideal, by the deformation invariance
it is enough to compute $\mdeg{M}$ for monomial ideals $M$. If the
codimension of $\Sigma(M)$ in $W$ is $s$, then the maximal
dimensional components of $\Sigma(M)$ are codimension-$s$ coordinate
subspaces of $W$. Such subspaces are indexed by subsets
$\mathbf{i}\in\{1\ddd N\}$ of cardinality $s$; the corresponding
associated primes $\mathfrak{p}[\mathbf{i}]=\langle y_i:i\in
\mathbf{i} \rangle$.

It is not difficult to check that
\begin{equation}
\label{primemon} \mult(\mathfrak{p}[\mathbf{i}],M)=
\left|\left\{\mathbf{a}\in\Z_+^{[{\mathbf{i}}]};\;
\mathbf{y}^{\mathbf{a}+\mathbf{b}}\notin M\text{ for all }
\mathbf{b}\in\Z_+^{[\hat{\mathbf{i}}]}\right\}\right|,
\end{equation}
where $\Z_+^{[\mathbf{i}]}=\{\mathbf{a}\in \Z_+^N;a_i=0 \text{ for }
i\notin \mathbf{i}\}$, $\hat{\mathbf{i}}=\{1\ddd
N\}\setminus\mathbf{i}$, and $|\cdot|$, as usual, stands for the
number of elements of a finite set.

Then by the normalization and additivity axiom we have
\begin{equation}
  \label{epdmon}
\mdeg{M,S} =
\sum_{|\mathbf{i}|=s}\mult(\mathfrak{p}[\mathbf{i}],M)
\prod_{i\in\mathbf{i}}\eta_i.
\end{equation}

\subsection{An example}
\label{sec:basic} A simple way to construct $T$-invariant
subvarieties of $W$ is to take the orbit closures of points in $W$.

Consider the following example: let $W=\C^4$ endowed with a
$T=(\C^*)^3$-action, whose weights $\eta_1, \eta_2,\eta_3$ and
$\eta_4$ span $\lie(T)^*$, and satisfy
$\eta_1+\eta_3=\eta_2+\eta_4$. In other words, the four weights,
$\eta_i$, $i=1,\dots,4$, form the vertices of a parallelogram in
$\lie(T)^*$ lying in a hyperplane which does not pass through the
origin. Choose $p=(1,1,1,1)\in W$; then the closure of the $T$-orbit
of $p$ is given by a single equation:
\begin{equation}\label{toriceq4}
\overline{\porb}=\{(y_1,y_2,y_3,y_4)\in\C^4;\; y_1y_3=y_2y_4\}.
\end{equation}
We will compute the equivariant Poincar\'e dual of this subvariety
in a number of ways.

\textsc{Method 1: Deformation of the variety.} We use the axioms
listed in \S\ref{subsec:axiomatic}.

\[I(\Sigma)=\langle y_1y_3-y_2y_4 \rangle \subset
S=\C[y_1,y_2,y_3,y_4]\] has initial ideal
\[\mathrm{in}_<(I(\Sigma))=\langle
y_1y_3 \rangle\] with respect to the lexicographic monomial order
corresponding to the order $y_1>y_2>y_3>y_4$ on the variables (see
\cite{eisenbud} \S 15.2). Note that $\mathrm{in}_<(I(\Sigma))$
defines the union of two hyperplanes: $\{y_1=0\}$ and $\{y_3=0\}$
with multiplicity $1$. Then, using the additivity and the
normalization axioms, we arrive at the result that the equivariant
Poincar\'e dual is $\epd{\Sigma}=\eta_1+\eta_3=\eta_2+\eta_4$, hence
\begin{equation}\label{simpleform}
\epd{\overline{\porb}}=\eta_1+\eta_3.
\end{equation}

\subsection{Some technical statements}
\label{sec:sometech} In the previous paragraphs we sketched the
construction and properties of the equivariant Poincar\'e dual. Here
we will discuss a few simple consequences of these properties.

We retain the notation of \S\ref{sec:epdmult}:  $W$ is a 
$T$-module endowed with coordinates $y_1\ddd y_N$, which are of
weight $\eta_1\ddd\eta_N$, respectively. The following technical
lemma will be crucial in our computations.
\begin{lemma}\label{divisible}
  Let $I\subset \C[y_1\ddd y_N]$ be a $T$-invariant ideal, and assume
  that for some $j$, $1\leq j\leq N$, there is an element $R\in I$
  which expresses the variable $y_j$ as a polynomial of the remaining
  variables:
  \begin{equation}
    \label{yl}
    R:\;  y_j=f(y_1\ddd y_{j-1},y_{j+1}\ddd y_N).
  \end{equation}
  Then $\mdeg{I,\C[y_1\ddd y_N]}$ is divisible by $\eta_j$. More
  precisely,
 \begin{equation}
   \label{ydis}
   \mdeg{I,\C[y_1\ddd y_N]} = \eta_j\cdot \mdeg{I_j,\C[y_1\ddd
     y_{j-1},y_{j+1}\ddd     y_N]}
 \end{equation}
 where $I_j$ the ideal in $\C[y_1\ddd y_{j-1},y_{j+1}\ddd y_N]$
 obtained from $I$ by performing the substitution \eqref{yl}.
\end{lemma}
\begin{proof}
  Let $<$ be the lexicographic monomial order on S induced by the
  ordering of the coordinates in the following way: $y_j>y_1> \ldots
  >y_{j-1}>y_{j+1}> \ldots>y_N$. Then $y_j$ is the initial monomial
  of $R$, therefore it is a generator of $\mathrm{in}_<(I)$. As we saw
  before, the prime monomial ideals $\mathfrak{p}[\mathbf{i}]$ of
  $\mathrm{in}_<(I)$ are indexed by subsets $\mathbf{i}\subset
  \{1,\ldots ,N \}$, and
  \begin{equation}\label{conclusion}
  y_j \in \mathrm{in}_<(I) \Rightarrow j\in \mathbf{i},
  \end{equation}
  by \eqref{primemon}. As a result, each nonvanishing term of
  the sum in \eqref{epdmon} will contain the factor $\eta_j$.  The
  second statement follows from the fact that
  \[\mathrm{in}_<(I_j)=\mathrm{in}_<(I)\cap\C[y_1\ddd y_{j-1},y_{j+1}\ddd
  y_N].\]
\end{proof}

\begin{rem}\label{geominterpr}
  The geometric version of Lemma \ref{divisible}, corresponding to the
  case when $I$ is reduced, reads as follows. Let $\Sigma\subset W$ be
  a closed $T$-invariant subvariety, and assume that the conditions of
  Lemma \ref{divisible} hold for $I(\Sigma)$.  Let $\pi_j:W\to W_j$
  denote the projection onto the hyperplane
  $W_j=\left\{y_j=0\right\}$. Then $\pi_j(\Sigma)$ is a closed
  subvariety in $W_j$ and
\[\epd{\Sigma,W}=\eta_j\cdot \epd{\pi_j(\Sigma),W_j}\]
Again, note that in this case the polynomial $\epd{\Sigma,W}$ is
divisible by $\eta_j$.
\end{rem}

\subsection{Integration and equivariant multiplicities}
\label{sec:intem}

In \cite{rossmann}, Rossmann made the important observation that
the notion of \ePd may be extended to the case of analytic
$T$-invariant varieties defined in a neighborhood of the origin in
$T$-representations, and further, to nonlinear actions, as we
explain below.

Let $Z$ be a complex manifold with a holomorphic $T$-action, and let
$M\subset Z$ be a $T$-invariant analytic subvariety with an isolated
fixed point $p\in M^T$. Then one can find local analytic coordinates
near $p$, in which the action is linear and diagonal. Using these
coordinates, one can identify a neighborhood of the origin in $\TT_pZ$
with a neighborhood of $p$ in $Z$. We denote by $\tc_pM$ the part of
$\TT_pZ$ which corresponds to $M$ under this identification;
informally, we will call $\tc_pM$ the $T$-invariant {\em tangent cone}
of $M$ at $p$. This tangent cone is not quite canonical: it depends on
the choice of coordinates; the equivariant Poincar\'e dual of
$\Sigma=\tc_pM$ in $W=\TT_pZ$, however, does not. Rossmann named this
\ePd the {\em equivariant multiplicity of $M$ in $Z$ at $p$}:
\begin{equation}\label{emult}
   \emu_p[M,Z] \overset{\mathrm{def}}= \epd{\tc_pM,\TT_pZ}.
\end{equation}

\begin{rem}
In the algebraic framework one might need to pass to the {\em
tangent
  scheme} of $M$ at $p$ (cf. \cite{fulton}). This is canonically
defined, but we will not use this notion.
\end{rem}

An important application of the equivariant multiplicity is
Rossmann's localization formula \cite{rossmann}.  The reader will
find the necessary background material about equivariant
differential forms and equivariant integration in \cite{gs,bgv}. For
technical reasons, we need to pass to the compact versions of our
reductive groups. We will use the notation $G_\circ$ for the compact
form of the complex reductive group $G$; for example $T_\circ$ will
be a product of copies of the circle group $U(1)$. The introduction
of these groups into our framework means an implicit choice of an
Hermitian metric.

Let $\mu:\lie(T_\circ)\to\Omega^\bullet(Z)$ be a holomorphic
equivariant map with values in smooth differential forms on $Z$.
Then {\em Rossmann's
  localization formula} states that
\begin{equation}
  \label{rossform}
  \int_M\mu=\sum_{p\in M^T}\frac{\emu_p[M,Z]}{\euler^T(\TT_pZ)}\cdot\mu^{[0]}(p),
\end{equation}
where $\mu^{[0]}(p)$ is the differential-form-degree-zero component
of $\mu$ evaluated at $p$.  Recall that $\euler^T(\TT_pZ)$ stands
for the product of the weights of the $T$-action on $\TT_pZ$.

This formula generalizes the equivariant integration formula of
Berline and Vergne \cite{bv}, which applies when $M$ is smooth. In
this case the tangent cone of $M$ at $p$ is a well-defined linear
subspace $\TT_pM\subset\TT_pZ$, and $\emu_p[M]$ is the \ePd of this
subspace.  Then the fraction in \eqref{rossform} simplifies: the
ambient space $Z$ is eliminated from the picture, and one arrives at
(cf. \cite{bv})
\begin{equation}
  \label{bvform}
  \int_M\mu = \sum_{p\in M^T}\frac{\mu^{[0]}(p)}{\euler^T(\TT_pM)}.
\end{equation}
Rossmann proves \eqref{rossform} by first expressing the equivariant
multiplicity in terms of an integral-limit, and then applying an
adaptation of Stokes theorem, following the method of Bott
\cite{bott}.

As showed by Vergne \cite{voj}, such a local integration formula for
equivariant Poincar\'e duals may be given in the framework of
equivariant cohomology. To describe this formula, we return to our
setup of a $T$-invariant subvariety $\Sigma$ in a complex vector
space $W$ of dimension $N$. The starting point is the Thom
isomorphism in equivariant cohomology:
\begin{equation}
  \label{thomg}
   H^*_{T_\circ,\cpt}(W) = H^*_{T_\circ}(W)\cdot\mathrm{Thom}_{T_\circ}(W),
\end{equation}
which presents compactly supported equivariant cohomology as a
module over usual equivariant cohomology. The class
$\mathrm{Thom}_{T_\circ}(W)\in H^{2N}_{T_\circ,\cpt}(W)$ may be
represented by an explicit equivariant differential form with
compact support (cf. \cite{mq,dv}). Then {\em
  Vergne's integration formula} (cf. \cite{voj}) reads as follows:
\begin{equation}
  \label{vergneepd}
  \epd\Sigma = \int_\Sigma\mathrm{Thom}_{T_\circ}(W).
\end{equation}

Compared to Rossmann's formula \eqref{rossform}, this result turns
things upside down, and describes $\epd\Sigma$ as an integral in
equivariant cohomology. As we explain in the next section, this
allows us to localize the equivariant Poincar\'e dual near fixed
points of torus actions.

We complete this review by noting that a consequence of
\eqref{vergneepd} is the following formula. For an equivariantly
closed differential form $\mu$ with compact support, we have
\[  \int_\Sigma\mu = \int_W \epd\Sigma\cdot\mu.
\]
This formula serves as the motivation for the term {\em equivariant
  Poincar\'e dual.}

\subsection{Thom polynomials and equivariant Poincar\'e duals}
\label{sec:thompoin} Let us apply our new-found invariant to the setup
of global singularity theory described in
\S\ref{sec:basicsing}. Recall that, for integers $d$ and $n\leq k$, we
have an irreducible variety $\Theta_d\subset\mapd nk$, which is
invariant under the natural action of the group $\diffdnk$.

Now observe that the quotient map $\lin:\diff n\to\dif_1(n)=\gl n$
has a canonical section, consisting of linear substitutions. In
other words, we have a canonical group embedding
\[ \gl n \hookrightarrow \diff n,\] and we can restrict the action of
the diffeomorphism groups $\diff k\times\diff n$ on $\mapd nk$ to
the canonical subgroup $\gl k\times \gl n$. Denoting the subgroups
of diagonal matrices of $\gl k$ and $\gl n$ by $T_k$ and $T_n$,
their basic weights by $\thetab=(\theta_1\ddd\theta_k)$ and
$\lambdab=(\lambda_1\ddd\lambda_n)$, respectively, we can introduce
the central object of our paper.
\begin{defi}  Let $A$ be a nilpotent algebra. The {\em Thom
    polynomial} of the $A$-singularity from $n$-to-$k$ dimension is
  \label{deftp}
  \begin{equation}
    \label{thompoldefi}
 \tp A^{n\to k}(\lambdab,\thetab)\overset{\mathrm{def}}=
\epd{\overline{\Theta}_A,\mapd
     nk}_{T_k\times T_n }.
  \end{equation}
\end{defi}
According to Proposition \ref{propgena}, this is a homogeneous
polynomial of degree $(k-n+1)\dim A$ in the variables
$\theta_1\ddd\theta_k,\lambda_1\ddd\lambda_n$.  Note that in case the
torus action extends to the action of the general linear group, the
symmetric group $\sg n$, thought of as the Weyl group, naturally acts
on the weights of $T$ by permuting the $\lambda$s. Thus we can
conclude the following.
\begin{lemma}\label{prop_symmetric}
  Let $T=(\C^*)^n$ be the subgroup of diagonal matrices of the complex
  group $\gl n$, and denote by $\lambda_1\ddd\lambda_n$ its basic
  weights. If $\Sigma$ is a $\gl n$-invariant subvariety of the $\gl
  n$-module $W$, then the equivariant Poincar\'e dual $\epd{\Sigma,W}_T$ is
  a {\em symmetric} polynomial in $\lambda_1\ddd\lambda_n$.
\end{lemma}
Clearly, this Lemma applies to our situation, hence we have
\begin{cor}\label{corsymm}
  The Thom polynomials $ \tp A^{n\to k}(\lambdab,\thetab)$ are {\em
    symmetric} in $\theta_1\ddd\theta_k$ and in
    $\lambda_1\ddd\lambda_n$.
\end{cor}

Starting with the next section we will focus on the computation of the
polynomial $ \tp A^{n\to k}(\lambdab,\thetab)$ for the case
$A=\C[t]/t^{d+1}$.  In the remainder of this paragraph, however, we
would like to argue that this polynomial is a reasonable candidate for
the universal class satisfying Thom's principle quoted in
\S\ref{sec:intro}.  This is standard for the experts
(cf.  \cite{rimanyi,kazarian,fr,prag}), but good references are hard
to come by. In any case, we would like to stress that this material is
not necessary for understanding the rest of the paper. The reader
comfortable with Definition \ref{deftp} may safely skip to
\S\ref{sec:structural}.

When comparing Thom's principle from \S\ref{sec:intro} to Definition
\ref{deftp}, we come up against several difficulties.  First: how to
relate equivariant Poincar\'e duals such as in \eqref{thompoldefi} to
the usual Poincar\'e class of corresponding cycles on $N$?  Next, how
can the replacement of the symmetry group $\diff k\times\diff k$ by
$\gl k\times \gl n$ in \eqref{thompoldefi} be justified?  And finally, in
the holomorphic category one cannot always deform a function into a
transversal position. What is the meaning of this polynomial in
this case?  We address the first question in Proposition \ref{sect}, and
the second in Proposition \ref{jets}. For more details we direct the
reader to the references listed above.

Now fix the notation $G=\gl n$ and $G_\circ=U_n$ for its maximal
compact subgroup. Let $F$ be a principal $G_\circ$-bundle over a
compact oriented manifold $M$. Then, using the Chern-Weil map, any
symmetric polynomial $P\in\C[\lambda_1,\dots,\lambda_n]^{\sg n}$
defines a characteristic class $P(F)\in H^*(M,\C)$.  Now let $\Sigma$
be $G$-invariant subvariety of the $G$-module $W$, and, denote by
$W_F$ the associated vector bundle $F\times_{G_\circ}W$ over $M$, and
by $\Sigma_F$ the subset of $W_F$ corresponding to $\Sigma$.
\begin{diagram}[LaTeXeqno,labelstyle=\textstyle]
\label{charclass} F\times_{G_\circ}W  & =  &  W_F & \lInto &
\Sigma_F & =
&  F\times_{G_\circ}\Sigma\\
  & & \uTo^s \dTo &  \ldTo  &  &  & \\
  &  & M &    &  &  &
\end{diagram}
Then by Poincar\'e duality on the manifold $W_F$, there is a
cohomology class \\
$\alpha_\Sigma\in H^{2\codim(\Sigma)}(W_F)$ such
that
\[   \int_{W_F} \alpha_\Sigma\cdot\beta = \int_{\Sigma_F} \beta
\]
for any compactly supported cohomology class on $W_F$. Thus the answer
to our first question maybe written as follows:
\begin{equation}\label{alphasigma}
\alpha_\Sigma=\epd{\Sigma,W}(F) \text{ in } H^*(W_F),
\end{equation}
i.e. the Chern-Weil image of the equivariant Poincar\'e dual is the
ordinary Poincar\'e dual of the induced variety.

We will prove this statement in a geometric form which is more
convenient for our purposes. In this setup $\epd{\Sigma,W}(F)$ will
appear as the Poincar\'e dual of $s^{-1}(\Sigma_F)$ in $M$ for an
appropriate section $s:M\to\Sigma_F$. To make this more precise, we
make the following
\begin{defi}
  \label{deftransversal}
  Consider the diagram \eqref{charclass}, and assume for simplicity
  that $\Sigma$ is equidimensional. We say that a {\em smooth} section
  $s:M\to W_F$ is {\em transversal} to $\Sigma_F$ at some point $p\in M$ if
  $s(p)$ is a smooth point of $\Sigma_F$ and the intersection
  $ds(\mathrm{T}_pM)\cap \mathrm{T}_{s(p)}\Sigma_F$ of vector spaces
  in $\mathrm{T}_{s(p)}W_F$ has the smallest possible dimension. We
  say that $s:M\to W_F$ is {\em generically transversal} to $\Sigma_F$
  if the we have
  \[ \overline{\{p\in M; s\text{ is transversal to }\Sigma_F\text{ at
    }p\}} = s^{-1}(\Sigma_F).
\]
\end{defi}

Armed with this technical notion, we reformulate
\eqref{alphasigma} as follows.
\begin{prop}
  \label{sect} For a smooth section $s:M\to W_F$ generically
  transversal to $\Sigma_F$, the cycle $s^{-1}(\Sigma_F)\subset M$ is
  Poincar\'e dual to the characteristic class $\epd\Sigma(F)$ of $F$
  corresponding to the symmetric polynomial
  $\epd {\Sigma,W}$.
\end{prop}
\begin{proof}
  Considering \eqref{vergneepd} as the definition of the equivariant
  Poincar\'e dual,  this statement becomes almost tautological. Indeed,
  recall Cartan's correspondence, which associates to an equivariantly
  closed differential form $\mu$ on a $G$-manifold $X$ an ordinary
  closed differential form $\Cw(\mu)$ on the manifold
  $X_F=F\times_GX$. There is a simple construction of this
  correspondence, which uses the Weyl algebra model for equivariant
  cohomology; the only necessary input is a connection on $F$
  \cite{bgv}. In particular, when $X=\mathrm{pt}$, then $\Cw$ reduces
  to the usual Chern-Weil correspondence. As  $\Cw$ clearly
  commutes with integration and restriction, considering forms with
  compact support, we obtain the following commutative diagram:
\begin{diagram}[LaTeXeqno,labelstyle=\textstyle]
\label{cartanmap}
H_{G,\mathrm{cpt}}^*(W)  & \rTo^{\Cw}  &    H_{\mathrm{cpt}}^*(W)   \\
\dTo^{\int_\Sigma} &   & \dTo_{\pi_*^\Sigma} \\
 H_G^*(\mathrm{pt}) & \rTo^{\Cw} &    H^*(M)
\end{diagram}
The symbol $\int_\Sigma$ here stands for integrating on
$\Sigma\subset W$, while $\pi_*^\Sigma$ is the push-forward along
the fibers of the bundle $\Sigma_F\to M$.

Now, starting with $\mathrm{Thom}_{G_\circ}(W)\in
H_{G,\mathrm{cpt}}^*(W)$ defined by \eqref{thomg} in the upper left
corner of the diagram, we arrive exactly at our statement. Indeed,
according to \eqref{vergneepd}, we have
\[   \Cw\left(\int_\Sigma\mathrm{Thom}_{G_\circ}(W)\right)=
\Cw(\epd\Sigma)=\epd\Sigma(F).
\]
On the other hand, the Cartan correspondence takes
$\mathrm{Thom}_{G_\circ}(W)$ to the Thom class of the bundle $W_F\to
M$, which is also the Poincar\'e dual of $M$ thought of as the zero
section in $W_F$. Now, using the properties of the Poincar\'e dual
(cf. \cite{botttu}), it is a simple exercise to check that the
push-forward is Poincar\'e dual to $s^{-1}(\Sigma_F)\subset M$ for a
section $s:M\to W_F$, generically transversal to $\Sigma_F$.
\end{proof}

Now let us look at the situation of singularity loci of holomorphic
maps described in the introduction; this appears to be similar to
the setup we have just considered.

Indeed, for complex manifolds $N$ and $K$ of dimensions $n$ and $k$,
respectively, and a positive integer $d$, consider the principal
$\diff k\times \diff n$-bundle $\diff K\times \diff N$ over the
product space $N\times K$ consisting of local coordinate changes up
to order $d$. Denote by $\mathcal{J}_d(N,K)$ the bundle over
$N\times K$ associated to the representation $\mapd nk$ of the group
$\diff k\times \diff n$. Note that even though the space $\mapd nk$
has a linear structure the action of the group $\diff k\times\diff
n$ on it is not linear, and hence this bundle is not a vector
bundle. Then any holomorphic map $f:N\to K$ induces a section
$s_f:N\to (1\times f)^*\mathcal{J}_d(N,K)$ of the bundle pulled back
from graph.

Now, for a nilpotent algebra $A$ satisfying $A^{d+1}=0$, consider
the subvariety
\begin{equation}
  \label{jetsub}
 \mathcal{J}_d(\Theta_A^{N\to K})\subset \mathcal{J}_d(N,K),
\end{equation}
associated to the subvariety $\Theta_A^{n\to k}\subset\mapd nk$.

Now we can state the main technical statement of this paragraph:
\begin{prop}
  \label{jets}
  Let $N,K,A$ and $d$ be as above. Let $f:N\to K$ be a {\em smooth} map
  and $s:N\to (1\times f)^*\mathcal{J}_d(N,K)$ be an arbitrary smooth
  section, generically transversal to $(1\times
  f)^*\mathcal{J}_d(\Theta_A^{N\to K})$.  Next, denote by
  $Q_A(\lambda_1\ddd\lambda_n,\theta_1\ddd\theta_k)$ the polynomial $\tp
  A^{n\to k}$ defined in \eqref{thompoldefi}. Then the
  cohomology class $Q_A(TN,f^*TK)\in H^*(N)$ is Poincar\'e dual
  to the subvariety $s_f^{-1}( (1\times f)^*\mathcal{J}_d(\Theta_A^{N\to K}))$.
\end{prop}
\begin{proof}
  One can repeat the above construction replacing the group $\diff
  k\times\diff n$ by its subgroup $\gl k\times\gl n$; then the
  subvariety \eqref{jetsub} is replaced by a subvariety
  $\tilde{J}_d(\Theta_A^{N\to K})$ of the tensor bundle
  $\hom(\oplus_{m=1}^d\sym^mTN,TK)$. For this pair, the statement of
  Proposition is an immediate consequence of Proposition \ref{sect}.

  Now the Proposition immediately follows from the structure group of
  the bundle $\mathcal{J}_d(N,K)$ considered in the smooth category,
  reduces to $\gl k\times \gl n$. This can be seen using that $\diff
  k\times\diff n$ is homotopy equivalent to $\gl k\times \gl n$ or,
  alternatively, by directly presenting the reduction using, for
  example, Hermitian metrics on $TN$ and $TK$ (cf. \cite[\S2.2]{kazarian}).
\end{proof}

\subsection{Thom polynomials of contact singularities}
\label{sec:structural}

One of the natural questions to ask is how the Thom polynomials for
fixed $A$ and different pairs $(n,k)$ are related.  We collect the
known facts \cite{arnold,damon,fr} in Proposition \ref{collect} below.
For simplicity, we will formulate the statements for the algebra
$A_d=t\C[t]/t^{d+1}$ we study, although essentially the same
properties are satisfied by the Thom polynomials of any other contact
singularity (see \cite{fr} for details).

Denote the ring of bisymmetric polynomials in the $\lambda$s and
$\theta$s by $\snk$, and recall from \S2.\ref{sec:basicsing} that for
$1\le d$ and $1\le n\le k$, $\Theta_d=\Theta_d^{n\to k}$ is a nonempty
subvariety of $\mapd nk$ of codimension $d(k-n+1)$. Consider the
infinite sequence of homogeneous polynomials $c_i\in\snk$, $\deg
c_i=i$, defined by the generating series
\begin{equation}
\label{deftc}
  \ct(q)=1+c_1q+c_2q^2+\dots =
\frac{\prod_{m=1}^k(1+\theta_mq)}{\prod_{l=1}^n(1+\lambda_lq)};
\end{equation}
we will call $c_i$ the $i$th {\em relative Chern class}.
\begin{prop}[\cite{fr}]
\label{collect} Let $1\le d$ and $1\le n\le k$.
 Then for each nonnegative integer $j$, there is a polynomial
    $\DD^j_d(b_0,b_1,b_2,\dots)$ in the indeterminates $b_0,b_1,b_2,\dots$
    with the following properties
    \begin{enumerate}
    \item $\DD^j_d$ is homogeneous of degree $d$, and
    \item if we set $\deg(b_i)=i$, then $\DD^j_d$ is homogeneous of
      degree $d(k-n+1)$;
    \item for all $1\leq n\leq k$, we have
      \begin{equation}
        \label{tptd}
 \tp d^{n\to k}(\boldsymbol{\lambda},\thetab)
=\DD^{k-n}_d(1,c_1(\boldsymbol{\lambda},\thetab),
c_2(\boldsymbol{\lambda},\thetab),\dots),
              \end{equation}
where the polynomials $c_i(\lambdab,\thetab)$, $i=1,\dots$, are
defined by (\ref{deftc});
            \item the polynomial $\DD^{j-1}_d$ may be
              obtained from $\DD^{j}_d$ via the following substitution:
\[ \DD^{j-1}_d(b_0,b_1,b_2,\dots)=\DD^{j}_d(0,b_0,b_1,b_2,\dots),
\]
  \end{enumerate}

\end{prop}
The notation TD stands for Thom-Damon polynomial.  The 3rd property
\eqref{tptd} is an older result of Damon and Ronga
(\cite{damon,sigma11}), while the 4th is a theorem of Feh\'er and
Rim\'anyi \cite{fr}.

There is a somewhat confusing aspect of \eqref{tptd}, which we would
like to clarify now.  For fixed $j$ and sufficiently large $n$ and
$k$, the polynomials $c_i(\lambdab,\thetab)\in\snk$, $i=1,\ddd d(j+1)$
are algebraically independent.  This means that for fixed codimension
$j$ and large enough $n$, the Thom polynomial $\tp d^{n\to
  n+j}(\boldsymbol{\lambda},\thetab)$ determines $\DD^j_d$. However,
for small values of $n$, the natural map
\[     \C[c_1,c_2,\dots]\to\snk
\]
is not surjective in degree $d(k-n+1)$, and in this case there are
several expressions of the Thom polynomial in terms of relative Chern
classes. Only one of these expressions remains valid for all $n$.

\begin{exa}
For $d=4,n=1,k=1,$
\[\mathrm{RC}(q)=\frac{1+\theta q}{1+\lambda
  q}=1+(\theta-\lambda)q-\lambda(\theta -\lambda)q^2+\ldots ,\] thus
we have
\[
\quad  c_0(\theta,\lambda)=1,\quad c_1(\theta,
\lambda)=\theta-\lambda, \quad c_2(\theta,
\lambda)=-\lambda(\theta-\lambda),\]
\[ \quad c_3(\theta,
\lambda)=\lambda^2(\theta-\lambda), \quad c_4(\theta,
\lambda)=-\lambda^3(\theta-\lambda)\dots
\]

We have (cf. \cite[Theorem 2.2]{gaffney}, also \S\ref{negyothat})
\[\DD_4^0=c_1^4+6c_1^2c_2+2c_2^2+9c_1c_3+6c_4c_0,\]
and for $n>1$, this is the only possible expression for the Thom
polynomial in terms of the relative chern classes.
However, since for $n=k=1$, 
\[c_1(\theta,\lambda)c_3(\theta,\lambda)=c_2(\theta,\lambda)^2,\]
we can conclude that
\[\tp 4^{1\to 1}(\theta,\lambda)=c_1^4+6c_1^2c_2+\alpha c_2^2+ (11-\alpha)
c_1c_3+6c_4c_0\] holds for any $\alpha \in \mathbb{R}$.
\end{exa}

Next, following \cite{fr}, observe that property (4) allows us to
define a universal object, the Thom series $\Ts(a_i,\,i\in\Z)$,
which is an infinite formal series in infinitely many variables with
the following properties:
\begin{itemize}
\item $\Ts(a_i,\,i\in\Z)$ is homogeneous of degree $d$;
\item setting $\deg(a_i)=i$ for $i\in\Z$, the series
  $\Ts_d(a_i,\,i\in\Z)$ is homogeneous of degree 0;
\item the Thom-Damon polynomial maybe expressed via the following
  substitution:
\[  \DD^{j}_d(b_0,b_1,b_2,\dots)=
\Ts_d\left(
  \begin{cases}
    a_i=b_{i+k-n+1}\, \text{ if }i\geq -(k-n+1),\\
    a_i=0\,\text{ otherwise}.
  \end{cases}
\right)
\]
\end{itemize}
For example, in this language Porteous's formula is simply $\Ts_1=a_0$,
while Ronga's formula takes the form $\Ts_2=a_0^2+\sum_{i=0}^\infty
2^{i-1}a_ia_{-i}$. This suggestive way of expressing Thom polynomials,
found by Feh\'er and Rim\'anyi, served as a starting point for our
work. We obtained a rather satisfactory answer, which manifestly has
the structure described above; the final result \eqref{thethom} even
gives some insight into the geometric meaning of the coefficients of
the Thom series.

\section{Localizing Poincar\'e duals}
\label{sec:loctech}

In this section we develop the idea introduced at the end of
\S\ref{sec:intem}: the localization of equivariant Poincar\'e duals
based on Vergne's integration formula. Roughly, we show that if the
$T$-invariant subvariety $\Sigma \subset W$ is equivariantly fibered
over a parameter space $M$, then the \ePd $\epd{\Sigma,W}$ may be
read of from local data near fixed points of the $T$ action on $M$.
The final form of the statement is Proposition \ref{locepdsing}. We
will start, however, with the more regular case of a smooth
parameter space.
\subsection{Localization in the smooth case} Let
$\Sigma$ be a $T$-invariant closed subvariety of the $T$-module $W$.
Consider the following diagram:
\begin{diagram}[LaTeXeqno,labelstyle=\textstyle]
\label{diagsmooth}
  & & & & W & \lInto &  \Sigma\\
  & & & \ruTo^{\evm} & \uTo^{\evs} &  & \\
 S_{\!M^T} & \rInto & S_{\!M} & \rTo & S &  & \\
 \dTo^{\tautt}    & & \dTo^{\taum} & &\dTo^{\taugr}   & &\\
 M^T  & \rInto^{\iota_T} & M & \rTo^\phi & \grass mW & &
\end{diagram}

Here
\begin{itemize}
\item $\grass mW$ is the Grassmannian of $m$-planes in $W$, $S$ is the
  tautological bundle over $\grass mW$, and $\taugr:S\to\grass mW$ is
  the tautological projection; observe that the tautological
  evaluation map $\evs:S\to W$ is proper.
\item $M$ is a smooth compact complex manifold, endowed with a
  $T$-action; as usual, the notation $M^T$ stands for the set
  $\{y\in M;\; Ty=y\}$ of fixed points of the $T$-action; assume
that $M^T$ is a finite set of points. The embedding
$M^T\hookrightarrow M$ is denoted by $\iota_T$.
\item  $\phi:M\to\grass mW$ be a $T$-equivariant map, and introduce
  the pull-back bundles $S_{\!M}=\phi^*S$ and
  $S_{\!M^T}=\iota_T^*S_{\!M}$; we denoted by $\evm$ the induced
  evaluation map $S_{\!M}\to W$.
\item For clarity, we indexed our spaces and maps, but these indices
  will be omitted whenever this does not cause confusion. For example
  if $p\in M$, then we will denote by $\Ss p$ the fiber of the bundle
  $\Ss M$ over the point $p$.
\end{itemize}

Literally, to say that $\Sigma$ is fibered over $M$ would mean that
the map $\evm:S_{\!M}\to W$ establishes a diffeomorphism of
$S_{\!M}$ with $\Sigma$. Since this essentially never happens, we
weaken this condition as follows.

Recall (see e.g. \cite{botttu}) that to a smooth proper map $f:X\to
Y$ between connected oriented manifolds of equal dimensions one can
associate an integer $\deg(f)$ called the {\em degree}.  This
constant may be defined via the equality
\begin{equation}
  \label{degsmooth}
\int_X   f^*\mu=\deg(f) \int_Y\mu,
\end{equation}
which holds for any compactly supported form $\mu$ on $Y$.

An alternative definition of $\deg(f)$ is the signed sum of the
preimages of a regular value; the sign associated to a preimage
depends on whether the map is orientation-preserving or reversing at
the point. Since a holomorphic map is orientation-preserving
everywhere, we have the following simple statement.
\begin{lemma}
  Let $f$ be a proper holomorphic map between complex manifolds. Then
  $f$ is of degree 1 if and only if there is dense open $U\subset X$
  such that $f$ restricted to $U$ is a biholomorphism onto a dense
  open subset of $Y$.
\end{lemma}

 The definition of a degree-1 map may be extended to the following
 situation.
\begin{defi}\label{degree1}
  Let $f:X \to Y$ be a smooth, proper map between complex manifolds,
  and $U\subset X$ and $V\subset Y$ not necessarily smooth closed
  analytic subvarieties.  We say that $f$ establishes a degree-1 map
  between $U$ and $V$ if there are Zariski open subsets $U^o\subset U$
  and $V^o\subset V$, not containing singular points, such that
  $f|_{U^o}:U^o \to V^o$ is biholomorphic. Here Zariski open means
  that the complement is a closed analytic subvariety.
\end{defi}

Another convenient way to describe our notion is
\begin{prop}\label{biholomorphic}
Let $f:X \to Y$ be a proper map of complex manifolds, $U\subset X$
possibly singular closed analytic subvariety. Suppose that there is
$U^o \subset U$ Zariski open subset, not containing singular points,
such that $f|_{U^o}$ is injective. Then $f$ establishes a degree-1
map between $U$ and $f(U)$.
\end{prop}

\begin{proof}
  Since $f$ is proper, $f(U)$ is a closed analytic subvariety of $Y$,
  (see \cite{griffiths}, page 34). Injectivity implies that
  $\dim(U^o)=\dim(V^0)$, and hence there is a possibly smaller Zariski
  open $U'\subset U^o$ such that $f(U')$ is in the smooth part of
  $f(U)$. Since an injective holomorphic map between manifolds is
  biholomorphic, can conclude that $f$ restricted to $U'$ is a
  biholomorphism, and this completes the proof.
\end{proof}

Now, we would like to extend the property \eqref{degsmooth} for the
singular degree-1 case.  A key fact is that integration of
differential forms with compact support may be extended to not
necessarily smooth analytic subvarieties of complex manifolds.

Let $\mu$ be a differential form with compact support on a complex
manifold $X$, and let $U\subset X$ be a closed analytic subvariety,
whose set of smooth points we denote by $U^s$,
$\iota:U^s\hookrightarrow U$. Then one defines
\begin{equation}
  \label{analintdef}
\int_{U}\mu \overset{\mathrm{def}}= \int_{U^s}\iota^*\mu.
\end{equation}

\begin{prop}
 The integral on the right hand side of \eqref{analintdef} is
 absolutely convergent, and vanishes if $\mu$ is exact.
\end{prop}
The reason for this is that  that in a local chart, with respect to
the euclidean metric, the submanifold  $U^s$
has finite volume in bounded regions (cf. \cite[\S2, p. 32]{griffiths}).  

The following two corollaries will be important for us.

\begin{cor}
  \label{cordegint} If  the map $f:X\to Y$ establishes a degree-1 map
  between $U$ and $V$ as in Definition \ref{degree1}, then
\[\int_{U}f^*\mu=\int_{V}\mu \quad\text{for every compactly supported
smooth form  }\mu \text{ on }Y.\]
\end{cor}

\begin{cor}\label{familyversion}
  Let $M$ be a complex manifold, $V$ be a complex vector bundle over
  $M$, and let $S\hookrightarrow V$ be a locally trivial subbundle
  with fibers which are possibly singular analytic subvarieties of the
  corresponding linear fibers of $V$.  Denote by $\pi:S \to M$ the
  projection. Then for any smooth compactly
  supported form $\mu$ on $V$, the push-forward of the restriction:
  $\pi_*\mu$ is a smooth form on $M$, moreover,
\[
\int_S\mu=\int_M\pi_*\mu.
\]
\end{cor}

Now we are ready to formulate our first localization formula.
\begin{prop}
  \label{locsmooth} Assume that in diagram \eqref{diagsmooth} the
  fixed point set $M^T$ is finite, and $\evm$ establishes a degree-1
  map from $S_{\!M}$ to $\Sigma$. Then we have
 \begin{equation}
   \label{loclin0}
\epd{\Sigma,W}=\sum_{p\in
M^T}\frac{\epd{\evm(S_{\!p}),W}}{\euler^T(\TT_pM)}.
 \end{equation}
\end{prop}
\begin{rem}
\begin{enumerate}
\item The most natural situation is when $M$ is a smooth submanifold
  of $\grass mW$. The more general setup we are considering in
  Proposition \ref{locsmooth} works, however, even when the image
  $\phi(M)$ is singular.
\item Since the space $\evm(S_{\!p})$ is a linear $T$-invariant
  subspace of $W$ for $p\in M^T$, the polynomial $\epd{\evm(
    S_{\!p})}$ is determined by the normalization axiom: it simply
  equals the product of those weights of $W$ which are not weights of
  $\evm(S_{\!p})$ (with multiplicities taken into account).
\item The equivariant Euler class in the denominator is also a product
  of weights (cf. \eqref{eec}), hence each term in the sum is a
  rational function. After the summation, however, the denominators
  cancel, and one ends up with a polynomial result.
\end{enumerate}
\end{rem}
\begin{proof}
  Vergne's integral formula, \eqref{vergneepd} combined with our
  assumption that $\evm:S_{\!M}\to\Sigma$ is degree-1,
  implies that
\[  \epd\Sigma = \int_{S_{\!M}}\ev^*_M\thom.
\]
Integrating first along the fibers, we obtain that
\[  \epd\Sigma = \int_M\tau_*\evm^*\thom,
\]
where the integrand $\tau_*\evm^*\thom$ is a smooth equivariant form
on $M$. Now we apply the Berline-Vergne equivariant integration
formula \eqref{bvform} to this form, and obtain that
\begin{equation}
  \label{temps}
  \epd\Sigma = \sum_{p\in
    M^T}\frac{\left(\tau_*\evm^*\thom\right)^{[0]}(p)}
{\euler^T(\TT_pM)},
\end{equation}
where, as usual, we denote by $\mu^{[0]}$ the
differential-form-degree-zero part of the equivariant form $\mu$.
Since $\evm$ is a linear injective map on each fiber, the numerator
of \eqref{temps} is simply the integral $\int_{\evm(S_{\!p})}\thom$.
Now, using Vergne's formula \eqref{vergneepd} once again, we arrive
at \eqref{loclin0}.
\end{proof}

In the remainder of this section we present examples of using this
formula, and also give a few variants of this result.

We first note that using remark \ref{passtoeuler}, formula
\eqref{loclin0} may be rewritten as follows.  Let $E$ be an
equivariant vector bundle over $M$, and let $\gamma_p:W\to E_p$ for
$p\in M$ be an equivariant family of surjective linear maps. Assume,
that this establishes a degree-1 map between the subbundle
\[   \{(p,w)\in M\times W;\; \gamma_p(w)=0\}
\]
and $\Sigma$. Then according to Remark \ref{passtoeuler}, we have
$\epd{\evm(S_{\!p}),W}=\euler^T(E_p)$, which leads to the following
variant of (\ref{loclin0}):
  \begin{equation}
    \label{euler}
     \epd{\Sigma}=\sum_{p\in
  M^T}\frac{\euler^T(E_p)}{\euler^T(\TT_pM)}
  \end{equation}

As a quick application, we give yet another way of computing the
\ePd for the example introduced in \S\ref{sec:basic}.\\
\textsc{Method 2: localization on the projectivized
  cone}. Consider the smooth, $T$-invariant projective variety
$\mP\Sigma\subset\mP^3$ cut out by the homogeneous equation
$x_1x_3=x_2x_4$. In the notation of \eqref{diagsmooth}, we have
$M=\mP\Sigma$, $m=1$ and $W=\C^4$. Then the fixed point set
$\mP\Sigma^T$ consists of the  four fixed points on $\mP^3$
corresponding to the four coordinate axes.

Pick one of these fixed points, say,  $p=(1:0:0:0)$, which
corresponds to the coordinate line $S_{\!p}=\{x_2=x_3=x_4=0\}$.
Using the normalization axiom, we have then
$\epd{S{\!_p}}=\eta_2\eta_3\eta_4$.

Turning to the denominator in \eqref{loclin0}, it is not hard to see
that
\[ \euler^T(\TT_p\mP\Sigma)=(\eta_2-\eta_1)(\eta_4-\eta_1).\] Indeed,
this is the standard yoga of toric geometry: consider the
parallelogram formed by the weights $\eta_1,\eta_2,\eta_3$ and
$\eta_4$; the fixed points of the torus action correspond to the
vertices of this parallelogram, and the weights at a particular
fixed point are the edge-vectors emanating from the associated
vertex.

The contributions at the other fixed points may be computed
likewise, and the result is the following complicated formula for
the equivariant Poincar\'e dual:
\begin{eqnarray}
\epd{\Sigma}=\frac{\eta_2\eta_3\eta_4}{(\eta_2-\eta_1)(\eta_4-\eta_1)}+
\frac{\eta_1\eta_3\eta_4}{(\eta_1-\eta_2)(\eta_3-\eta_2)}+ \nonumber \\
\frac{\eta_1\eta_2\eta_4}{(\eta_2-\eta_3)(\eta_4-\eta_3)}+
\frac{\eta_1\eta_2\eta_3}{(\eta_1-\eta_4)(\eta_3-\eta_4)}.
\end{eqnarray}
This rational function is {\em not} a polynomial, however, assuming
$\eta_1+\eta_3=\eta_2+\eta_4$ holds, it can be easily shown to
reduce to the simple form \eqref{simpleform}.

We note that this procedure may be applied, inductively, to more
general toric varieties, and, again, the data may be read off the
corresponding polytope. However, if the polytope is not simple, then
the prescription is more involved.

\subsection{An interlude: the case of $d=1$} \label{sec:interlude}

In this paragraph, we consider the case $d=1$ of the
$A_d$-singularities introduced in \S\ref{sec:adsing}, and recover
the classical result of Porteous.

We have $\mathcal{J}_1(n,k)=\hom(\C^n,\C^k)$, and
$\Theta_1\subset\mathcal{J}_1(n,k)$ consists of those linear maps
$\C^n\to\C^k$ whose kernel is 1-dimensional. These maps may be
identified with $k$-by-$n$ matrices, and the weight of the action on
the entry $e_{ji}$ is equal to $\theta_j-\lambda_i$. Then the
closure $\overline{\Theta}_1$ consist of those $k$-by-$n$ matrices
which have a nontrivial kernel:
\begin{equation}
  \label{model1}
  \overline{\Theta}_1 =
\{A\in\hm kn;\; \exists v\in\C^n,\,v\neq0:\,Av=0\}.
\end{equation}

This description immediately suggests us an equivariant birational
fibration of $\overline{\Theta}_1$ over $\mP^{n-1}$, fitting the
conditions of Proposition \ref{locsmooth}: the fiber over a point
$[v]\in\P^{n-1}$ is the linear subspace
$\{A;\;Av=0\}\subset\overline{\Theta}_1$; where $[v]$ stands for the
point in $\mP^{n-1}$ corresponding to the nonzero vector $v\in\C^n$.

Again, we simply need to collect our fixed-point data, and then
apply \eqref{loclin0}. There are $n$ fixed points, $p_1\ddd p_n$ in
$\mP^{n-1}$, corresponding to the coordinate axes. The weights of
$\TT_{p_i}\mP^{n-1}$ are $\{\lambda_s-\lambda_i;\;s\neq i\}$. The
fiber at $p_i$ is the set of matrices $A$ with all entries in the
$i$th column vanishing. Again, using the normalization axiom, this
shows that the \ePd of the fiber at $p_i$ is
$\prod_{j=1}^k(\theta_j-\lambda_i)$, so our localization formula
looks as follows:
\begin{equation}\label{sigma1}
\epd{\overline{\Theta}_1}=\sum_{i=1}^n\frac{\prod_{j=1}^k(\theta_j-\lambda_i)}{\prod_{s\neq
i}(\lambda_s-\lambda_i)}
\end{equation}
This is a finite sum for fixed $n$, but as $n$ increases, the number
of terms also increases. There is a way, however, to further
``localize'' this expression, and obtain a formula, which only
depends on the local behavior of a certain function at a single
point.

Indeed,  consider the following rational differential form on
$\mP^1$:
\[
-\frac{\prod_{j=1}^{k}(\theta_j-z)}
{\prod_{i=1}^{n}(\lambda_i-z)}\,dz.\] Observe that the residues of
this form at finite poles: $\{z=\lambda_i;\;i=1\ddd n\}$ exactly
recover the terms of the sum \eqref{sigma1}. Then, applying the
Residue theorem, we obtain
\[   \epd{\overline{\Theta}_1}=\res_{z=\infty}
\frac{\prod_{j=1}^{k}(\theta_j-z)}
{\prod_{i=1}^{n}(\lambda_i-z)}\,dz.\]

Finally, after the change of variables $z\to-1/q$, we end up with
\[ \epd{\overline{\Theta}_1}=\res_{q=0}\frac{\prod_{j=1}^{k}(1+q\theta_j)}
{\prod_{i=1}^{n}(1+q\lambda_i)}\,\frac{dq}{q^{k-n+2}},
\]
which, according to \eqref{deftc}, is exactly the relative Chern
class $c_{k-n+1}$. Thus we recovered the well-known Giambelli-Thom
-Porteous formula (\cite{forporteous}; \cite{griffiths} Chapter
I.5).

As a final remark, note that our basic example introduced in
\S\ref{sec:basic} is a special case of $\overline{\Theta}_1$,
corresponding to the values $n=k=2$. Hence this computation provides
us with a {\sc 3rd method } of arriving at \eqref{simpleform}. This
computation uses localization, similarly to the {\sc 2nd method},
but the two constructions are different.

\begin{equation}
\epd{\Sigma}=\frac{\eta_1 \eta_2}{\eta_3-\eta_2}+\frac{\eta_3
\eta_4}{\eta_2-\eta_3}
\end{equation}
Using $\eta_1+\eta_3=\eta_2+\eta_4$, we arrive to the formula
\eqref{simpleform}.

\subsection{Variations of the localization formula}

We will need to amend and generalize Proposition \ref{locsmooth} in
two ways in order to be able deal with $\Theta_d$ for $d>1$: we will
drop the assumption on that the fibers are linear, and we will also
allow $M$ to be singular.

\subsubsection{Nonlinear fibers} \label{sec:nonlinfib}
Next, observe that, during the proof of Lemma \ref{locsmooth}, we
never used the assumption that the fibers are linear spaces. In
fact, using Corollary \ref{familyversion}, the same formula and the
same argument holds if the fibers of $S$ are  possibly singular
analytic subvarieties.

\begin{prop}
\label{locsingfibr} Let $\Sigma$ be a closed subvariety of the
complex vector space $W$. Assume that $M$ is a smooth compact
complex manifold, $V$ is a complex vector bundle over $M$, and let
$S\hookrightarrow V$ be a locally trivial subbundle
  with fibers which are possibly singular analytic subvarieties of the
  corresponding linear fibers of $V$. Suppose that we have a proper
map: $\ev_V:V \to W$, which establishes a degree-1 map from $S$ to
$\Sigma$. Then
\begin{equation}
   \label{loclin}
\epd{\Sigma,W}=\sum_{p\in
M^T}\frac{\epd{\ev_V(S_{\!p}),W}}{\euler^T(\TT_pM)}.
 \end{equation}
\end{prop}
\noindent We will use this variant of the localization in
\ref{subsec:locind}, for the localization on a flag variety.

\subsubsection{Fibrations over a singular base} Finally, we remove
the assumption that $M$ is smooth. For brevity, below, without
explicitly stating this, we will assume that every space and map is
in the $T$-equivariant category. We will apply the following
proposition for the localization on $\O$ in \ref{subsec:locfib}.

\begin{prop}
  \label{locepdsing}
  \begin{enumerate}
  \item  Let $\Sigma$ be a closed subvariety of the
  complex vector space $W$.  Assume that $Z$ is a compact, smooth
  complex manifold, and $M\subset Z$ is a possibly singular, closed
  subvariety with a finite set of fixed points $M^T$.
  Consider the following analog of diagram \ref{diagsmooth}:
  \begin{diagram}[LaTeXeqno,labelstyle=\textstyle]
\label{diagsing}
  & & & & W & \lInto &  \Sigma\\
  & & & \ruTo^{\ev_Z} & \uTo^{\evs} &  & \\
 S_{\!M} & \rInto & S_{\!Z} & \rTo & S &  & \\
 \dTo^{\tau_M}    & & \dTo^{\tau_Z} & &\dTo^{\taugr}   & &\\
 M  & \rInto^{\iota_M} & Z & \rTo^\phi & \grass mW & &
\end{diagram}
  Assume that $\ev_Z$ establishes a degree-1 map
  between $\tau_Z^{-1}(M)$ and $\Sigma$.

Then
\begin{equation}
   \label{locsing}
\epd{\Sigma}=\sum_{p\in M^T}
\frac{\epd{\ev_Z(S_p)}\,\emu_p[M,Z]}{\euler^T(\TT_pZ)}.
\end{equation}
\item Assume that there is a T-equivariant vector bundle $E$ over $M$,
and an equivariant family of surjective linear maps $\gamma_p:W\to
E_p$ for $p\in M$, such that the set
\[   \{(p,w)\in M\times W;\; \gamma_p(w)=0\}
\]
is a subbundle of the trivial bundle $M\times W$, and it maps to
$\Sigma$ in a birational fashion. Then
\[\epd{\Sigma}=\sum_{p\in M^T}
\frac{\euler^T(E_p)\,\emu_p[M,Z]}{\euler^T(\TT_pZ)}.\]
\end{enumerate}
\end{prop}

\begin{proof}
The second part is the combination of the first part and
\eqref{euler}. The proof of the first part is analogous to that of
Proposition \ref{locsmooth}; when passing to \eqref{temps}, however,
one needs to use Rossmann's integration formula \eqref{rossform}.
 \end{proof}

\section{The test curve model}
\label{sec:model}

In \S\ref{sec:basicsing}, we described the variety $\Theta_d$ in two
different ways: as an example of a contact singularity class defined
in (\ref{defconea}), and as the Boardman class corresponding to the
sequence $(1,1,\dots,1)$ (cf. Prop. \ref{charlocalg}).  In this
section, we recall another, birationally equivalent description of
$\Theta_d$ -- the so-called ``test curve model'' -- which goes back
to the works of Porteous, Ronga, and Gaffney
\cite{por,ronga,gaffney}. Roughly, the idea of the construction is
to generalize \eqref{model1} to $d>1$ by requiring that the map-jet
$\Psi\in\mapd nk$ carry a $d$-jet of a curve in $\C^n$ to zero.  As
we have not found a complete proof of the appropriate statement
(Theorem \ref{model}) in the literature, we give one below.

Recall the notation $\lin:\mapd nk\to \hom(\C^n,\C^k)$ for the
linear part of map-jets. A $d$-jet of a curve in $\C^n$ is simply an
element of $\mapd 1n$.  We will call such a curve $\Gammac$ {\em
regular} if $\lin(\Gammac)\neq0$; introduce the notation $\mapreg$
for the set of these curves:
\begin{equation} \label{mapreg} \mapreg\overset{\mathrm{def}}= \left\{
    \Gammac\in\mapd 1n;\, \lin(\Gammac)\neq0\right\}.
\end{equation}

Now define the set
\begin{equation}\label{modelset}
\Upsilon_d = \left\{\Psi\in\mapd nk;\;
    \exists\Gammac\in\mapreg\text{ such that } \Psi\circ\Gammac=0\right\}.
\end{equation}
In words: $\Upsilon_d$ is the set of those $d$-jets of maps, which
take at least one regular curve to zero. By definition, $\Upsilon_d$
is the image of the closed subvariety of the quasi-projective $\mapd
nk \times \mapreg$ defined by the algebraic equations $\Psi \circ
\Gammac=0$, under the projection to the first factor. By a theorem of
Chevalley (see \cite{hartshorne}, Ex. 3.19, page 94), the set
$\Upsilon_d$ is constructible. We will not use the set $\Upsilon_d$
itself in this paper, rather its Zariski closure: the variety
$\overline{\Upsilon}_d\subset\mapd nk$.
\begin{thm}\label{model} The Zariski closures of $\Theta_d$
  and $\Theta_d'$ in $\mapd nk$ coincide.
\end{thm}
\begin{proof}
  Recall from Proposition \ref{korbit} that $\Theta_d$ is an orbit of
  the complex algebraic group $\mathcal{K}_d$ defined in
  \eqref{defk}. To prove the theorem, it is then sufficient to show
  that
\begin{itemize}
\item $\Upsilon_d$ is $\mathcal{K}_d$-invariant,
\item $\Upsilon_d \cap \Theta_d$ is nonempty,
\item $\codim(\overline{\Upsilon}_d)=\codim(\Theta_d)$ in $\mapd nk$, and that
\item the subvariety $\overline{\Upsilon}_d\subset\mapd nk$ is irreducible.
\end{itemize}

Indeed, to see that these 4 statements are sufficient, we observe that
according Propositions \ref{korbit} and \ref{propgena}, $\Theta_d$ is
a single, irreducible $\mathcal{K}_d$-orbit.  This fact, with the
first two properties above induces that $\Theta_d \subset \Upsilon_d$,
so $\overline{\Theta}_d \subset \overline{\Upsilon}_d$. Since
$\overline{\Upsilon}$ is irreducible of the same dimension as
$\overline{\Theta}$, $\overline{\Theta}_d = \overline{\Upsilon}_d$
must hold.

To show the $\mathcal{K}_d$-invariance of $\Upsilon$, observe that
if $\Gammac\in\mapd 1n$ is regular and $\Delta\in\diff n$, then
$\Delta\circ\Gammac$ is also regular. Indeed, in this case
\[\lin(\Delta \circ \Gammac)=\lin (\Delta) \cdot \lin (\Gammac) \neq 0.\]
Now, if $\Psi\in\mapd nk$ such that $\Psi\circ\Gammac=0$ for some
regular $\Gammac$, and $(M,\Delta)\in\mathcal{K}_d$, then recalling
the action  \eqref{via}, we have
\[
[(M,\Delta)\cdot\Psi] \circ(\Delta \circ \Gammac)= (M\cdot\Psi)
\circ \Delta^{-1}\circ (\Delta \circ \Gammac)=(M\cdot \Psi) \circ
\Gammac =(M\circ\Gammac)\cdot(\Psi\circ\Gammac)=0.
\]
This shows that $\Delta \circ \Gammac$ is an appropriate test curve
for the transformed map-jet $(M,\Delta)\cdot\Psi$.

To find an element in the intersection of $\Theta_d$ and
$\Upsilon_d$, consider the map-jet
\[\Psi_{0}(x_1,\ldots, x_n)=(0,x_2,\ldots, x_n,0,\ldots, 0).
\]
It obviously belongs to $\Theta_d$; on the other hand, for the test
curve $\Gammac(t)=(t,0,\ldots, 0)$, we have $\lin (\Gammac) \neq 0$
and $\Psi_{0}\circ \Gammac=0$ in $\mapd nk$, hence $\Psi_{0}\in
\Upsilon_d$.

Regarding the codimensions, we have $\codim(\Theta_d)=d(k-n+1)$
according to Proposition \ref{propgena}. The proof of the
irreducibility of $\Upsilon_d$ and the computation of its
codimension (cf. Proposition \ref{lincors}) will follow from the more
detailed study of its structure, to which we devote the rest of this
section.

\end{proof}

Our first project is to write down the equation $\Psi\circ\Gammac=0$
in coordinates. This is a rather mechanical exercise, and we will
spend some time setting up the notation.

A curve $\Gammac\in\mapd 1n$ is parametrized by $d$ vectors
$v_1,\dots,v_d$ in $\C^n$:
\begin{equation}
  \label{gammaform}
\Gammac(t) = t v_1 + t^2 v_2 + \dots + t^d v_d,
\end{equation}
In this explicit form, the condition of regularity,
$\lin(\Gammac)\neq0$, simply means that $v_1\neq 0$.

Next, we switch to a new parametrization of our space $\rdk$.
Separating the similar homogeneous components of the $k$
polynomials, $P_1\ddd P_k$, and thinking of a homogeneous degree-$l$
polynomial as an element of $\hom(\sym^l \C^n,\C)$, we may represent
$\blp\in\rdk$ as a linear map
\begin{equation}
  \label{psidef}
\blp=(\blp^1,\ldots ,\blp^d):\oplus_{l=1}^d\sym^l \C^n \to \C^k.
\end{equation}

The standard basis of the vector space $\oplus_{l=1}^d\sym^l \C^n$
may be parametrized by nondecreasing sequences of positive integers,
or, alternatively -- and this is the language we will prefer -- by
{\em
  partitions}. Namely, to the partition $[i_1,\ddd i_l]$ of the
integer $i_1+\dots+i_l$ with $1\leq i_m\leq n$, we associate the
basis element $e_{i_1}\cdots e_{i_l}\in\sym^l \C^n$.

In what follows, certain integer characteristics of partitions will
be used.
\begin{notation}\label{partition}
 For a partition $\tau=[i_1\ddd i_l]$ of the integer
$i_1+\ldots +i_l$, introduce
\begin{itemize}
\item
 the
\textsl{length}: $|\tau|=l$,
\item  the \textsl{sum}: $\um\tau=i_1+\ldots +i_l$,
\item the \textsl{maximum}: $\max(\tau)=\max(i_1,\ldots,i_l)$,
\item and the \textsl{number of permutations}: $\comb(\tau)$, which is
  the number of different sequences consisting of the numbers
  $i_1,\dots, i_l$; e.g. $\comb([1,1,1,3])=4$.
\end{itemize}
\end{notation}

Denoting the set of all nonempty partitions by $\Pi$, we can
parametrize the basis elements of $\oplus_{l=1}^d\sym^l \C^n$ by the
finite set
\begin{equation}
  \label{finiteq}
  \{\tau\in\Pi;\;|\tau|\leq d,\,\max(\tau)\leq n\}.
\end{equation}
We will also use the notation $\Pi[m]$ for the set of all partitions
of the positive integer $m$:
\begin{equation}\label{pim}
    \Pi[m] =
\{\tau\in\Pi;\;\um\tau=m\}.
\end{equation}

Next, for a map-jet $\Psi\in\mapd nk$, a sequence
$\bv=(v_1,v_2,\dots)$ of vectors in $\C^n$, and a partition
$\tau=[i_1,\dots,i_l]$ satisfying $l\leq d$, $\max(\tau)\leq n$,
introduce the shorthand
\begin{equation}
  \label{notationv}
\bv_\tau = \prod_{j=1}^l v_{i_j}\in\sym^l\C^n\;\text{ and }\;
\Psi(\bv_\tau)=\Psi^l(v_{i_1},\dots,v_{i_l})\in\C^k.
\end{equation}

Armed with this new notation, we can write down the equation
$\Psi\circ\Gammac=0$ more explicitly, as follows.
\begin{lemma}\label{explgp} Let  $\Gammac\in\mapd 1n$ be given
  in the form \eqref{gammaform}. Then, using the notation
  \eqref{notationv}, the equation $\Psi\circ\Gammac=0$ is equivalent to
  the following system of $d$ linear equations with values in $\C^k$
  on the components $\Psi^l$ of $\,\Psi\in\mapd nk$, $l=1\ddd d$:
\begin{equation}
  \label{modeleq}
\sum_{\tau\in\Pi[m]} \comb(\tau) \,\Psi(\bv_\tau)=0,\quad
m=1,2,\dots, d.
\end{equation}
\end{lemma}

Let us see what the system of equations \eqref{modeleq} looks like
for small $d$. To make the formulas easier to follow, we will use
the $l$th capital letter of the alphabet for the symmetric
multi-linear map $\Psi^l$ introduced in \eqref{psidef}: we will
write $A$ for the linear part $\Psi^1$ of $\Psi$, $B$ for its second
order part, etc. With this convention (see also \eqref{gammaform}),
the system of equations for $d=4$ reads as follows:
\begin{eqnarray} \label{eqn4}
& A(v_1)=0, \\ \nonumber & A(v_2)+B(v_1,v_1)=0, \\ \nonumber
& A(v_3)+2B(v_1,v_2)+C(v_1,v_1,v_1)=0, \\
&
A(v_4)+2B(v_1,v_3)+B(v_2,v_2)+3C(v_1,v_1,v_2)+D(v_1,v_1,v_1,v_1)=0.
\nonumber
\end{eqnarray}

For a curve $\Gammac\in\mapreg$, introduce the notation
$\epsilon(\gamma)$ for the system of equations
\eqref{modeleq}\footnote{We will give a more formal meaning to
  $\epsilon$ in the next section.}, and
\begin{equation}\label{solepsilon}
\sol_{\epsilon(\Gammac)} \text{ for the space of solutions of this
system.}
\end{equation}
Then, according to \eqref{modelset},
\begin{equation}
  \label{upsdef}
  \Upsilon_d=\bigcup\left\{\sol_{\epsilon(\Gammac)};
\;\Gammac\in\mapreg\right\}.
\end{equation}

In the following Proposition, we collect some simple facts about the
system \eqref{modeleq}.
\begin{prop}
\label{modelprops}
\begin{enumerate}
\item \label{normalize} Let $0\neq v\in\C^n$, and assume that
  $\Gammac\in\mapreg$ is such that  $\lin(\Gammac)$ is parallel to $v$. Pick a
 hyperplane $H$ in $\C^n$ which is complementary to $v$. Then there is
 a unique $\delta\in\diff1$ such that
 \begin{equation}
   \label{gammadelta}
\Gammac\circ\delta=tv+t^2v_2+\dots+t^dv_d\quad \text{with}\;
v_2,v_3\ddd v_d\in H.
    \end{equation}
\item For \,$\Gammac\in\mapreg$,  the set of
  solutions $\sol_{\epsilon(\Gammac)}\subset\mapd nk$ is a linear
  subspace of codimension $dk$.
\item Introduce the set
\[\mapdo nk=\left\{ \Psi \in \mapd nk|\dim\ker(\lin(\Psi))=1\right\}.\]
Then for any $\Gammac \in \mapreg$, $\sol_{\epsilon(\Gammac)}\cap \mapdo
nk$ is a dense subset of $\sol_{\epsilon(\Gammac)}$.
\item If $\Psi \in \mapdo nk$, then $\Psi$ may belong
  to at most one of the spaces $\sol_{\epsilon(\Gammac)}$. More precisely,
\[
\text{if }\Gammac,\Gammac'\in\mapreg,\;\dim(\ker\lin(\Psi))=1,
\text{ and }
  \Psi\circ\Gammac=\Psi\circ\Gammac'=0,
\]
then there exists
  $\delta\in\diff1$ such that $\Gammac'=\Gammac\circ\delta$.
\item Given $\Gammac,\Gammac'\in\mapreg$, we have
  $\sol_{\epsilon(\Gammac)}=\sol_{\epsilon(\Gammac')}$ if and only if there is a
  $\delta\in\diff1$ such that $\Gammac'=\Gammac\circ\delta$.
\end{enumerate}
  \end{prop}
\begin{proof}

For (1), write explicitly $\Gammac(s) = sw_1+\dots+s^dw_d$ and
$\delta = \lambda_1t+\dots+\lambda_dt^d$. After performing the
substitution $s\mapsto \delta,$ we obtain a curve
$\Gammac\circ\delta=tv+t^2v_2+\dots+t^dv_d$, where $v_l=\lambda_l
w_1+$terms with $\lambda$s which have lower indices than $l$; this
clearly implies the statement.

The second statement follows from the presence of the term
$\Psi^l(v_1,\dots v_1)$ in the $l$th equation of \eqref{modeleq},
which is clearly linearly independent of the rest of the terms in
the first $l$ equations.

For statement (3), let $\Gammac=(v_1,\ldots, v_d)$, $v_1\neq
0$ be as in \eqref{gammaform}, and consider the linear map
\begin{equation}
  \label{linsur}
\lin:\sol_{\epsilon(\gamma)}\to \mathrm{Hom}(\C^n,\C^k)
\end{equation}
 associating
to each solution  $\Psi =(A,B,\ldots )\in
\sol_{\epsilon(\Gammac)}$ of the system  \eqref{eqn4} its component
$A$. Using the same argument as in the proof of statement (2), we can
see that for each fixed $A$ with $v_1\in \ker (A)$, the system
\eqref{eqn4}, becomes a system of $d-1$ linear equations with values
in $\C^k$, whose solution is a $(d-1)k$-codimensional linear subspace
in the space of the rest of the components  $(B,C,\ldots )$. In
particular, this shows that \eqref{linsur} is surjective, and this
implies statement (3).

To prove statement (4), we  assume that $\Gammac$ and $\Gammac'$ are
normalized according to \eqref{gammadelta} with respect to some
$v\in\ker(\lin(\Psi))$; then we show that $\Gammac=\Gammac'$ using
induction. Assume, for example, that the two curves coincide up to
the third order, i.e. $v_1=v_1',\,v_2=v_2',\,v_3=v_3'$. Then we see
from \eqref{eqn4} that $A(v_4)=A(v_4')$. We have $A=\lin(\Psi)$ and
$\ker(A)=\C v_1$, hence $v_4,v_4'\in H$ and $A(v_4)=A(v_4')$ imply
$v_4=v_4'$. This completes the inductive step.

The last statement is an immediate consequence of statement (2),(3)
and (4).
\end{proof}

The construction of this section are summarized in the following
diagram:
\begin{diagram}[LaTeXeqno,labelstyle=\textstyle]
\label{diagups}
\Theta_{d}& \rInto &\Upsilon_d & \rInto & \mapd nk  & \lTo^{\ev_S} & S\\
&&&&&\ldTo_{\taugr} & \\
 &  & \mapreg & \rTo^{\tilde{\phi}} &
\grnk  & \lTo^{\phigr}  & \qmr
\end{diagram}

Explanations:
\begin{itemize}
\item Each space in the diagram carries an action of the group
  $GL(k)\times GL(n)$, and the maps are equivariant with respect to
  this action.
\item As usual, we denote by $S$ the tautological bundle over the
  Grassmannian, and by $\ev_S$ the tautological evaluation map
  (cf. diagram \ref{diagsmooth}). To streamline our notation, we
  denote by $\grnk$ the variety of linear subspaces of {\em codimension}
  $dk$ in $\mapd nk$; hence, the rank of the bundle $S$ equals to
  $\dim(\mapd nk)-dk$.
\item $\tilde{\phi}:\mapreg \to \grass{-dk}{\mapd nk};\
  \Gammac\mapsto\sol_{\epsilon(\Gammac)}$ was introduced in
  \eqref{solepsilon}.
\item The space $\qmr$ denotes the topological quotient $\jddiff$ and
  the map $\phigr$ is induced by $\tilde{\phi}$ (see Propositions
  \ref{lincors} (1) and \ref{propmodel} below).
\end{itemize}

Now we have
\begin{prop}\label{lincors}
  \begin{enumerate}
  \item The map $\tilde{\phi}$ is $\diff1$-invariant, and the
  induced map $\phigr$ on the orbits is injective.
  \item The map $\ev_S$ restricted to $[\tau_{\mathrm{Gr}}]^{-1}\overline{\im(\tilde{\phi})}$ is of degree 1 onto
  $\overline{\Upsilon_d}$.
  \item \label{irreducibility}
$\overline{\Upsilon_d}$ is an irreducible subvariety of $\mapd nk$.
  \item \label{codimension}
$\codim(\overline{\Upsilon}_d)=d(k-n+1)$.
  \end{enumerate}
\end{prop}

\begin{proof} The first statement immediately follows from
Proposition \ref{modelprops} (2) and (5), while the second is a
consequence of Proposition \ref{modelprops} (3),(4) and Proposition
\ref{biholomorphic}: $\ev_S$ is injective on
$\mapdo nk\cap [\tau_{\mathrm{Gr}}]^{-1}\im(\tilde{\phi})$.

To prove the third statement, we rewrite \eqref{upsdef} in terms of
diagram \eqref{diagups}:
   \begin{equation}
    \label{upsdescr}
   \Upsilon_d = \ev_S\left([\taugr]^{-1}\im(\tilde{\phi})\right)
  \end{equation}
  As the map $\ev_S$ is proper, we have
  \begin{equation}\label{upsdescr2}
   \overline{\Upsilon_d} =
   \ev_S\left([\taugr]^{-1}\overline{\im(\tilde{\phi})}\right)
  \end{equation}
  Now, the Zariski closure of the image of an irreducible variety
  under a morphism is irreducible, and so is a vector bundle over an
  irreducible variety.  Applying this to the morphisms $\tilde{\phi}$
  and $\ev_S$, and to the restriction of the vector  bundle $S$ to
  $\im(\tilde{\phi})$, we obtain (3).

Finally, note that the fibers of this vector bundle are
codimension-$dk$ vector spaces in $\mapd nk$ (Proposition
\ref{modelprops}, (2)), while the base
$\overline{\im(\tilde{\phi})}$ has dimension $d(n-1)$ by the first
statement: the dimension of $\mapreg$ is $dn$, and the dimension of
$\diff1$ is $d$. Hence the codimension of $\Upsilon_d$ equals
$dk-d(n-1)=d(k-n+1)$.
\end{proof}

In what follows, the second statement of Proposition \ref{lincors}
will be crucial, as it provides us with a fibered model of the
singularity locus $\overline{\Theta}_d$. In view of Proposition
\ref{lincors}, it is natural to try to endow the quotient $\qmr$ with
a complex structure such that $\phigr$ is a morphism; then, in our
model \eqref{upsdescr}, we could replace $\im(\phi)$ by the image of
the {\em injective} morphism $\phigr$.

This is indeed possible, as we show below. First, however, we recall
some basic facts related to quotienting of complex manifolds
(e.g. \cite[\S9]{lee}). 
\begin{prop}
  \label{quotienting}
  A free action of a complex Lie group $G$ on a complex manifold $M$
  is proper if and only if the topological quotient $M/G$ may be
  endowed with the structure of a complex manifold such that the
  canonical map $\pi:M\to M/G$ is holomorphic.  In this case, this
  complex structure is unique, and any $G$-invariant holomorphic map
  $f: M\to K$ factors through $M/G$, i.e. the unique map
  $\tilde{f}:M/G\to K$ for which $f=\tilde{f}\circ\pi$ is holomorphic.
\end{prop}

\begin{prop}
  \label{propmodel}
  There is a smooth algebraic bundle with affine fibers $\qmr\to
  \P^{n-1}$ and a holomorphic map $\rho:\mapreg \to \qmr$ which is
surjective, $\diff1$-invariant  and separates the $\diff1$ orbits.
\end{prop}
\begin{proof}
  It will be convenient to identify $\mapd 1n$ with $\hom(C^d,\C^n)$,
  i.e with the set of $n$-by-$d$ matrices.
  Then $\mapreg$ is the set of matrices with nonvanishing first column, while the action of
  $\diff1$ is represented by multiplication by $d$-by-$d$
  matrices (cf. Lemma \ref{hgroup}). For a curve $\gamma\in\mapd 1n$ we will denote the
  $(i,m)$th entry of the corresponding matrix by $\gamma[i,m]$; this is
  the same as the $i$th coordinate of the the vector $v_m$ in the
  parametrization \eqref{gammaform}.

  Now we can formalize the first part of Proposition \ref{modelprops}
  as follows.  Let
\[\mapreg_i=\left\{\gamma\in \mapreg:
\gamma[i,1]\neq 0 \right\}\]
and
\[  U_i = \{\gamma\in\mapreg;\; \gamma[i,1]\neq0
\text{ and }\gamma[i,m]=0\text{ for }m>1\}
\]

According to Proposition \ref{modelprops} (1), for each
$\gamma\in\mapreg_i$ there exists a unique $h^i(\gamma)\in \diff1$
such that $\gamma \cdot h^i(\gamma)\in U_i$.  Moreover, it is clear
from matrix form in Lemma \ref{hgroup} that the entries of
$h^i(\gamma)$ are polynomials in the entries of $\gamma$ and
$\gamma[i,1]^{-1}$. This defines a map
\[\rho_i:\mapreg_i \rightarrow U_i,\quad \rho_i:\gamma\mapsto\gamma \cdot
h^i(\gamma),
\]
which establishes a one-to-one correspondence between the
$\diff1$-orbits of $\mapreg_i$ and $U_i$.

This allows us to construct an algebraic manifold $\qmr$ with
coordinate patches $U_i$, $i=1\ddd n$, and transition functions
\[\phi_{i,j}:U_i\cap\{\gamma;\;\gamma[j,1]\neq0\} \to U_j,\quad
\gamma \mapsto \gamma \cdot h^j(\gamma). \]
Since these transition functions are compatible with those of the
projective space $\P^{n-1}$, we can conclude that $\qmr$ has the
structure of an algebraic bundle over $\P^{n-1}$ with affine fib-res.

Now, by the construction, the maps $\rho_i$, $i=1\ddd n$ assemble into
an algebraic map
\[\rho:\mapreg \to Q_d,\]
which establishes  a one-to-one correspondence between the
$\diff1$-orbits of $\mapreg$ and $\qmr$, which is what we needed to
show.
\end{proof}

\begin{cor}\label{corqmr}
  The map $\tilde{\phi}$ on diagram \eqref{diagups} induces an
  injective holomorphic map
\[  \phigr: \qmr\to \grnk
\]
such that $\im(\phigr)=\im(\tilde{\phi})$.
\end{cor}

Note that in view of  Proposition \ref{lincors} (2) and Corollary
\ref{corqmr},  diagram \eqref{diagups} seems to fit the scheme of diagram
\eqref{diagsing}, with $\qmr$ playing the role of $M$.

Recall, however, that the localization formulas of \S\ref{sec:loctech}
apply to compact manifolds.  While the injective map $\phigr$ suggests
a reasonable compactification of $\qmr$: the closure of $\im(\phigr)$
in $\grnk $, the corresponding localization
computations would be very difficult. The choice of the
compactification is very important from the point of view of the
efficiency of resulting formulas, and we will be very careful in
constructing one. This is the subject of the next section.

Another approach would be finding a general quotienting procedure
resulting in a compact space representing the quotient of $\mapd 1n$
with respect to the action of the nonreductive group $\diff1$. The
problem of finding such an analog of the Geometric Invariant Theory of
Mumford \cite{mumford} is addressed in the recent work by Brent Doran
and Frances Kirwan \cite{kirwan}; the comparison of our constructions
with their results should provide us with new insights. Thus we hope
that our work represents a step in the direction of creating an
effective theory of localization on nonreductive quotients.

\section{The compactification}
\label{sec:compact}

As we observed at the end of the previous section, the morphism
$\tilde{\phi}$ in diagram \eqref{diagups} may be used to compactify
$\qmr=\jddiff$, and, in principle, allows us to apply the localization
techniques of \S\ref{sec:loctech}. The resulting formulas turn out
to be intractable, however, and the purpose of this section is to
replace the Grassmannian by a ``smaller'' space, which provides us
with a better compactification and, hopefully, with more efficient
formulas.

The constructions of this section form the backbone of the paper; we
will employ two ideas. The first is straightforward: we note that the
system of equations \eqref{modeleq} has a special form respecting a
certain filtration, and thus not every $dk$-codimensional linear
subspace of the Grassmannian may appear as the solution space of a
system of our equations.  These special systems give us a smaller
space to consider (cf. \S\ref{sec:embed}).

The second idea, detailed in \S\ref{subsec:fibr}, is a bit more
involved. The main features of this construction are removing a
certain part of the space of regular curves, thus breaking the
$\diff1$-symmetry, and then fibering the remainder over the space of
full flags of $d$-dimensional subspaces of $\C^n$. This leads to a
double fibration, whose study we are able to reduce to that of a
single fiber.

\subsection{Embedding into the space of equations}
\label{sec:embed} We start by rewriting the linear system
$\Psi\circ\Gammac=0$ associated to $\Gammac\in\mapd 1n$ in a dual
form (cf. Lemma \ref{explgp}). The system is based on the standard
composition map \eqref{comp}:
\[ \mapd nk\times\mapd 1n \longrightarrow\mapd 1k,
\]
which, in view of $\mapd nk=\mapd n1\tensor\C^k$, is derived from
the map
\begin{equation*}
 \mapd n1\times\mapd 1n \longrightarrow\mapd11
\end{equation*}
via tensoring with $\C^k$. Observing that composition is linear in
its first argument, and passing to linear duals, we may  rewrite
this correspondence in the form
\begin{equation}
  \label{mapson}
 \phip:\mapd 1n\longrightarrow\hom(\mapd11{}^*,\mapd n1{}^*).
\end{equation}
To present this map explicitly, we recall (cf. \eqref{gammaform})
that a $d$-jet of a curve $\Gammac\in\mapd 1n$ is given by a
sequence of $d$ vectors in $\C^n$, and thus, as a vector space, we
can
\begin{equation}
  \label{ident}
\text{identify }   \mapd 1n\text{ with }\hom(\C^d,\C^n).
\end{equation}
Also, according to \eqref{psidef}, the dual of $\mapd n1$ is the
vector space $\symdot=\oplus_{l=1}^d\sym^l \C^n$, hence a system of
$d$ linear equations on $\mapd n1$ may be thought of as a linear map
$\epsilon\in\hom(\C^d,\symdot)$; the solution set of this system is
the linear subspace orthogonal to the image of $\epsilon$:
$\im(\epsilon)^\perp\subset\mapd n1$ (cf. Definition \ref{epsdef}
below).

Using these identifications, we can recast the map $\phip$ in
\eqref{mapson} as
\begin{equation}
  \label{homs}
\phip:\hom(\C_L^d,\C^n) \longrightarrow\hom(\C_R^d,\symdot),
\end{equation}
which may be written out explicitly as follows (cf. \eqref{eqn4}):
\[   \phip:(v_1\ddd v_d)\longmapsto
\left(v_1,v_2+v_1^2,v_3+2v_1v_2+v_1^3,\dots,\sum_{\um\tau=m}\comb(\tau)\,
\bv_\tau,\dots\right).
\]
Note that in \eqref{homs} -- anticipating what is to come - we
marked the two copies of $\C^d$ with different indices: $L$ for left
and $R$ for right (cf. {\bf Convention} after Lemma
\ref{thecorrespondence} below).

The constructions of this section will be based on the observation
that the spaces of map germs $\mapd n1$ and $\mapd11$ -- and hence
their duals -- have natural filtrations, and these filtrations are
preserved by the map $\phip$.

The filtration on the dual of $\mapd n1$ (cf. \eqref{psidef}) is
\begin{equation}\label{filtsymdot}
\symdot = \oplus_{l=1}^d\sym^l\C^n\supset
\oplus_{l=1}^{d-1}\sym^l\C^n\supset\dots
\supset\C^n\oplus\sym^2\C^n\supset\C^n;
\end{equation}
setting $n=1$, this reduces to $\C^d$ with the standard filtration:
\begin{equation}
  \label{filtcd}
\C^d\supset\oplus_{l=1}^{d-1}\C e_l\supset\dots\supset\C e_1
\oplus\C e_2\supset\C e_{1}.
\end{equation}

Now introduce the notation $\hom^\triangle(\cdot,\cdot)$ for the
linear space of morphisms of filtered vector spaces. Then we have
\begin{equation}
  \label{filtered}
 \hofi =
\{\epsilon\in \hom(\C^d_R,\symdot);\;
\epsilon(e_l)\in\oplus_{m=1}^l\sym^m\C^n, l=1\ddd d\}.
\end{equation}
We will also need two open subsets of $\hofi$: the set of {\em
  nondegenerate systems}
\begin{equation}\label{fdn}
 \tfnk=\{\epsilon\in\hofi;\;\ker(\epsilon)=0\},
\end{equation}
and the set of {\em regular nondegenerate systems}
\begin{equation}
 \label{fdnreg}
 \tfnkreg =
\{\epsilon\in \hofi;\;
\epsilon(e_l)\notin\oplus_{m=1}^{l-1}\sym^m\C^n, l=1\ddd d\}.
\end{equation}

The following property of the map $\phip$ is manifest (cf.
Proposition \ref{modelprops}(2)):
\begin{lemma}
\label{thecorrespondence}
  The correspondence $\phip$ given in \eqref{homs} takes values in
  $\hofi$.
\end{lemma}

\noindent {\bf Convention:} The group of linear automorphisms of
$\C^d$ will be denoted, as usual by $\GL_d$, its subgroup of diagonal
matrices by $T_d$, and its subgroup of upper-triangular matrices by
$B_d$. In what follows, the two (left and right) copies of $\C^d$
appearing in (\ref{homs}) will play rather different roles. To avoid
any confusion, we will use the following notation for the corresponding
groups:
\[   T_L\subset B_L\subset \gl L\quad\text{and}\quad
 T_R\subset B_R\subset \gl R.
\]
\smallskip

The space $\hofi$ carries a left action of $\gl n$, and also a right
action of the Borel subgroup $B_R$ of $\gl R$ preserving the
filtration \eqref{filtcd}. Indeed, we have
\begin{equation}
  \label{bddef}
  B_R = \{b\in\hom^\triangle(\C_R^d,\C_R^d);\; b \text{ invertible}\}.
\end{equation}

\begin{lemma}
  \label{fdnquot} The subspaces $\tfnk$ and $\tfnkreg$ of $\hofi$ are
  invariant under both $\gl n$ and $B_R$. The quotient
  $\fnk=\tfnk/B_R$ is a compact, smooth manifold endowed with a $\gl
  n$-action, while
  $\fnkreg=\tfnkreg/B_R\subset\fnk$ is a $\gl n$-invariant open subset.
\end{lemma}

\begin{proof}
  To check the invariance with respect to the group actions is
  straightforward. The quotient $\tfnk/B_R$ may be described as the
  total space of a tower of $d$ fibrations as follows. The base of the
  tower is $\P(\C^n)$, and a fiber of the first fibration over a line
  $l_1\in \P(\C^n)$ is $\P((\C^n \oplus \sym^2\C^n)/l_1)$. Next, the
  fiber of the second fibration over a point $(l_1,l_2)\in(\P(\C^n),
  P(\C^n \oplus \sym^2\C^n/l_1)$) is $\P(\C^n \oplus \sym^2\C^n\oplus
  \sym^3\C^n)/(l_1+l_2)$, etc. This tower, which we denote by $\fnk$,
  is clearly a smooth, compact manifold.  More formally, this
  construction defines a surjective holomorphic map $\tfnk\to \fnk$,
  which is a bijection on the orbits, and hence (cf. Proposition
  \ref{quotienting}) $\fnk$ is the quotient $\tfnk/B_R$. Finally, since
  $\tfnkreg$ is open in $\tfnk$, then so is $\fnkreg$ in $\fnk$.
 \end{proof}

\begin{rem}
\label{remfibr} The space $\fnk$ may also be thought of as a
Schubert variety in the flag variety of  the  partial flag manifold
of full flags of $d$-dimensional subspaces of $\symdot$:
\begin{equation}\label{flag}
\mathrm{Flag}(\symdot)=\left\{ 0=F_0\subset F_1\subset\dots\subset
F_d\subset \symdot,\; \dim F_l=l\right\}.
\end{equation}
Lakshmibai and Sandhya in \cite{laksan} (see also \cite{gasrei},
Theorem 1.1) give combinatorial-type conditions under which a
Schubert variety is smooth, and $\fnk$ satisfies these conditions.
\end{rem}

Before proceeding, we introduce some notation associated with the
quotient in Lemma \ref{fdnquot}.
\begin{defi}
\label{epsdef}
  For $\epsilon\in\hofi$, thought of as a system of
equations, introduce the notation
\begin{itemize}
\item $\sol_\epsilon$ for the solution set
  $\im(\epsilon)^\perp\tensor\C^k\subset\mapd nk$,
  (cf. \eqref{solepsilon}) and
\item $\teps$ for the point in $\fnk$ corresponding to $\epsilon$.

\item Clearly, $\sol_\epsilon=\sol_{\epsilon b}$ for
  $\epsilon\in\hofi$ and $b\in B_R$, hence to each element
  $\teps\in\fnk$ we can associate a solution space $\sol_\teps$.
\end{itemize}
\end{defi}

The family of subspaces $\sol_\teps$ forms a holomorphic bundle over
$\fnk$ as the following statement shows.

  \begin{lemma}
\label{vident} Consider the bundle $V$ over $\fnk$ associated to the
standard
  representation of $B_R$:   $V = \tfnk\times_{B_R}\C_R^d$.
Then the canonical pairing
\begin{equation}\label{canpairing0}
  \tfnk\times\mapd n1 \to \hom(\C_R^d,\C)
\end{equation}
induces a linear bundle map from the trivial bundle with fiber
$\mapd n1$ over $\fnk$ to $V^*$:
\[s:\fnk\to\hom(\mapd n1,V^*)
\]
such that for $\teps\in\fnk$, we have
$\ker(s(\teps))\tensor\C^k=\sol_{\teps}\subset\mapd nk$.
 \end{lemma}
\noindent The upshot of this identification is the following exact
sequence of vector bundles over $\fnk$:
\begin{diagram}[LaTeXeqno,labelstyle=\textstyle]
  \label{exact}
  0 & \rTo &  \solfnk & \rTo^{\mathrm{ev}}&  \mapd nk &
  \rTo^s&  V^*\tensor\C^k & \rTo & 0,
\end{diagram}
where the fiber of $\solfnk$ over $\teps$ is the subspace
$\sol_\teps$.

After these preparations we return to our main task: the replacement
of the Grassmannian in diagram \eqref{diagups} by a smaller
variety. Observe that \eqref{exact} induces a morphism
\begin{equation}
  \label{fnkgrass}
   \phic:\fnk \to \grnk .
\end{equation}
Lemmas \ref{thecorrespondence} and \ref{fdnquot} imply that
$\im(\tilde{\phi})=\im(\phigr)\subset\im(\alpha)$, and hence, were
$\alpha$ injective, we could argue that the map $\phip$
(cf. \eqref{homs}) induces an injective morphism from $\qmr$ to
$\fnk$.  This seems reasonable since $\tilde{\phi}$ clearly factors
through the map $\phic$. There is a subtlety here, however: the map
\eqref{fnkgrass} is not injective, thus we need to exercise some extra
care. Indeed, for example, let $d=3$, and take the points
\[\epsilon_1=(v_1,v_2,v_1^2)\text{ and }\epsilon_2=(v_1,v_1^2,v_2)\]
in $\widetilde{\mathcal{F}}_3(n)$. Then $
\sol_{\epsilon_1}=\sol_{\epsilon_2}$, hence
$\alpha(\teps_1)=\alpha(\teps_2)$, but $\teps_1\neq\teps_2$.

The following statement resolves our problem.
\begin{lemma}\label{inj}
We have
  \begin{itemize}
  \item    $\phip(\mapreg)\subset\tfnkreg$, and
  \item the map $\alpha$ (defined in \eqref{fnkgrass}) restricted to
    $\fnkreg$ is an injective algebraic map.  \end{itemize}
\end{lemma}

\begin{proof}
We have $v_1\neq 0$ for $(v_1,\ldots, v_d) \in \mapreg$, and hence the
term  $v_1^d$ in
\[\phip((v_1,\ldots, v_d))(e_d)=v_1^d+(d-1)v_1^{d-2}v_2+\ldots \]
does not vanish; this proves the first statement.
To show the second, recall from Remark \ref{remfibr} that $\fnk$ may
be thought of as a subvariety of the flag variety \eqref{flag}.
Now, given $\epsilon\in\tfnk$, we have
$\alpha(\teps)=\im(\epsilon)^\perp\otimes\C^k$, which clearly
determines the vector space $U=\im(\epsilon)$. This in turn defines a
sequence of vector spaces
\begin{equation}\label{inverse}
((U\cap \C^n) \subset (U\cap (C^n\oplus\sym^2\C^n)) \subset (U\cap
(\oplus_{l=1}^3\sym^l\C^n)) \subset \ldots \subset U).
\end{equation}
According to \eqref{fdnreg}, this is a flag when $\epsilon\in\tfnk$,
which means that we can recover $\teps$ form $\alpha(\teps)$ if
$\teps\in\fnk$.
\end{proof}
\begin{rem}
If $\epsilon\notin\tfnkreg$, then  \eqref{inverse} with
$U=\im(\epsilon)$ will not define a flag, as some of the subspaces in
the sequence will coincide.
\end{rem}
\begin{rem}
  Note that $\epsilon_1$ and $\epsilon_2$ in the example above are not
  in the image $\phip(\mapreg)$.
\end{rem}

Using Lemma \ref{inj}, we can define the map
\begin{equation}\label{defphit}
\phif= \alpha^{-1}\circ\phigr,
\end{equation}
where the domain of definition of $\alpha^{-1}$ is understood to be
$\mathrm{im}(\alpha|_{\fnkreg})$. This allows us to reformulate our
model as follows.

\begin{cor}
  \label{jemb}
  The map $\phip$ in \eqref{homs} induces an algebraic morphism
\[ \phif: \qmr \rightarrow \fnk,
\]
Moreover, $\phif{}^*(\solfnk) = \alpha^*(S)$, hence by
\eqref{upsdescr} and \eqref{upsdescr2}, we have
\begin{equation}
  \label{fmodel}
   \Upsilon_d = \mathrm{ev}_{\widetilde{\mathcal{F}}}
\left(\tau_{\widetilde{\mathcal{F}}}^{-1}[\im(\phif)]\right)
\end{equation}
and
\begin{equation}
  \label{fmodel1}
  \overline{\Theta}_d=\overline{\Upsilon}_d =
  \mathrm{ev}_{\widetilde{\mathcal{F}}}
\left(\tau_{\widetilde{\mathcal{F}}}^{-1}[\overline{\im(\phif)}]\right);
\end{equation}
finally, the map $\mathrm{ev}_{\widetilde{\mathcal{F}}}$ in
$\eqref{fmodel1}$ establishes a degree-1 map from
$\tau_{\widetilde{\mathcal{F}}}^{-1}[\im(\phif)]$ to
$\overline{\Upsilon}$.
\end{cor}

Combining diagram \eqref{diagups} and sequence \eqref{exact}, we
arrive at the following picture:
\begin{diagram}[LaTeXeqno,labelstyle=\textstyle]
\label{diagfirst}
S  & &\solfnk & \rInto^{\mathrm{ev}_{\widetilde{\mathcal{F}}}} &
\mapd nk & \rOnto^s & V^*\tensor\C^k  \\
&\rdTo && \rdTo^{\tau_{\widetilde{\mathcal{F}}}} &  & \ldTo & \\
&& \grnk  & \lTo_{\;\;\;\;\;\;\alpha}  & \fnk &  \lTo & \tfnk \\
& & & \luTo_{\phigr} & \uTo_{\phif} & &\uTo^{\phip} \\
& & & & \qmr & \lTo^{\rho} & \mapreg
\end{diagram}

\begin{rem}\label{zariskiopen}
  The closure $\overline{\im(\phif)}$ gives us a new compactification
  of the space $\qmr=\jddiff$.
\end{rem}

\subsection{Fibration over the flag variety}
\label{subsec:fibr}

In the previous paragraph we took advantage of the special
``filtered'' form of the system \eqref{modeleq}, and replaced the
Grassmannian from \eqref{diagups} with the space of linear systems
$\fnk$. In this second part of the section, we further refine this
construction.

We start with a closer look at the ``natural'' identification
\eqref{ident}. In fact, the two objects are rather different: $\mapd
1n$ is a module over $\diff n\times\diff 1$ while $\hom(\C^d,\C^n)$
is a module over $\gl n\times \gl d$; in addition, note that we have
the following somewhat odd inclusions:
\begin{equation}
  \label{cruc}
\diff 1\subset\gl d,\quad \gl n\subset\diff n.
\end{equation}
By a straightforward computation, the first of the two inclusions
may be made more precise as follows.
\begin{lemma}
  \label{hgroup}
  Under the identification \eqref{ident}, a substitution
\[\alpha_1t+\alpha_2t^2+\ldots +\alpha_dt^d \in \diff1\]
corresponds to the upper-triangular matrix
\[
\left(
\begin{array}{ccccc}
\alpha_1 & \alpha_2   & \alpha_3          & \ldots & \alpha_d\\
0        & \alpha_1^2 & 2\alpha_1\alpha_2 & \ldots & 2\alpha_1\alpha_{d-1}+\ldots \\
0        & 0          & \alpha_1^3        & \ldots & 3\alpha_1^2\alpha_ {d-2}+ \ldots \\
0        & 0          & 0                 & \ldots & \cdot \\
\cdot    & \cdot   & \cdot    & \ldots & \alpha_1^d
\end{array}
 \right);
\]
the coefficient in the $i$th row and $j$th column is
\[\sum_{\left\{\tau\in\Pi[j];\,|\tau|=i\right\}}\comb(\tau)\,\alpha_{\tau},\]
where the notation $\alpha_{\tau}=\prod_{i\in \tau}\alpha_i$ was
used. This correspondence establishes an isomorphism of $\diff1$
with a $d$-dimensional subgroup $H_d$ of the Borel subgroup
$B_d\subset\gl d$.
\end{lemma}

\begin{rem}
  In accordance with the convention introduced after
  Lemma \ref{thecorrespondence}, we will use the notation $H_L$ when
  working with the copy of the group $H_d$ in the ``left'' Borel
  subgroup $B_L$.
\end{rem}

Now we return to the identification \eqref{ident} of $\mapd1n$ with
$\hom(\C_L^d,\C^n)$, and consider consider the subspace of injective
linear maps:
\begin{equation}
\label{hnd} \hnd=\{\Gammac\in\hom(\C_L^d,\C^n);\;\ker(\Gammac)=0\}
\end{equation}
The following statements are standard:
\begin{lemma}
\label{flagquot}
\begin{itemize}
\item Under the identification \eqref{ident}, the space $\hnd$ is a
  dense, open subset of $\mapreg$.
\item The action of $B_L$ on $\hnd$ is free, and the quotient $\hnd/B_L$ is the
  compact, smooth variety of full flags of $d$-dimensional subspaces of
  $\C^n$:
\[  \fld=\left\{
0=F_0\subset F_1\subset\dots\subset F_d\subset\C^n,\; \dim
F_l=l\right\}.
\]
\item The residual action of $\gl n $ on $\fld$ is transitive.
\end{itemize}
  \end{lemma}
Since fibrations over $\fld$ will play a major role in what follows,
we introduce some notation related to the quotient described in
Lemma \ref{flagquot}.
\begin{defi}\label{refsequence}
 \begin{itemize}
  \item  Denote by $\eref$  the {\em reference} sequence
  \[
\eref=(e_1\ddd e_d)\in\homreg,
\]
where $e_i$ is the $i$th basis vector of $\C^n$, and and we use the
identification \eqref{ident}. Let $\reffl$ denote the corresponding
flag in $\fld$.
\item For a space $X$ endowed with a left $B_L$-action, denote by $\ind X$ the
induced space $\ind X=\hnd\times_{B_L}X$.
  \end{itemize}
\end{defi}

Note that, in particular, we have $\hnd = \ind{B_L}$, and, according
to Lemma \ref{hgroup},
\begin{equation}
  \label{indhl}
  \hnd/H_L = \ind{\eref B_L/H_L}.
\end{equation}

This equality means that we have managed to fiber a Zariski-open part
of $\qmr$ over $\fld$.  This suggests investigating the systems of
equations \eqref{modeleq} in a single fiber of this fibration; we will
take a closer look at the fiber $\eref B_L$ lying over the point
$\reffl\in\fld$.

To inspect these systems, we will write them down in the standard
basis of $\symdot$; using the notation introduced in
\S~\ref{sec:model}, this consists of the elements
\[   e_\tau=e_{i_1}\cdot \ldots \cdot e_{i_m},\; \text{ where }
\tau=[i_1,\dots,i_m],\,m=|\tau|\leq d,\text{ and } \max(\tau)\leq n.
\]
We will denote the corresponding components of $\Psi\in\mapd nk$ by
\[ \Psi_\tau=\Psi^m(e_{i_1},\dots,e_{i_m}).
\]
We start with the {\em reference system}
$\epsilon_\rf=\phip(\eref)$:
\begin{equation}
  \label{reference_eq}
  \epsilon_{\rf}=\left\{
\sum_{\um\tau=l} \comb(\tau)\, \Psi_\tau=0,\quad l=1,2,\dots,
d\right\}.
\end{equation}
With the convention of using the $m$th capital letter of the
alphabet for $\Psi^m$, the first four equations of $\epsilon_{\rf}$
look as follows:
\begin{eqnarray} \label{eqncomp}
& A_1=0 \\ \nonumber & A_2+B_{11}=0 \\ \nonumber &
A_3+2B_{12}+C_{111}=0 \\ \nonumber & A_4 +
2B_{13}+B_{22}+3C_{112}+D_{1111}=0
\end{eqnarray}

Now consider a general element of $\eref B_L$, a test curve over the
reference flag:
\[
\eref\cdot\left(\begin{array}{cccc}
\beta_{11} & \beta_{12} & \beta_{13} & \cdot \\
0        & \beta_{22} & \beta_{23} & \cdot \\
0        & 0       & \beta_{33} & \cdot \\
\cdot    & \cdot   & \cdot    & \cdot \\
\end{array} \right)
= (\beta_{11}e_1, \beta_{22}e_2+ \beta_{12} e_1,
\beta_{33}e_3+\beta_{23}e_2+\beta_{13}e_1 \ldots).
\]
 The first 3 equations of the corresponding system \eqref{modeleq} are
\begin{eqnarray} \label{eqnoverreffl}
  & \beta_{11}A_1=0 \\ \nonumber &
  \beta_{22}A_2+\beta_{12}A_1+(\beta_{11})^2B_{11}=0 \\
  \nonumber &
  \beta_{33}A_3+\beta_{23}A_2+\beta_{13}A_1+
2\beta_{11}\beta_{22}B_{12}+2\beta_{11}\beta_{12}B_{11}+(\beta_{11})^3C_{111}=0;
\end{eqnarray}
these are thus of the form
\begin{eqnarray} \label{eqngen}
& u^1_1A_1=0 \\ \nonumber & u^2_2A_2+u^2_1A_1+u^2_{11}B_{11}=0 \\
\nonumber &
u^3_3A_3+u^3_2A_2+u^3_1A_1+2u^3_{12}B_{12}+u^3_{11}B_{11}+u^3_{111}C_{111}=0,
\end{eqnarray}
with some complex coefficients of the form $u^m_\tau$, where $m$ is
the ordinal number of the equation, while $\tau$ marks the component
of $\Psi$.  We observe that in the $l$th equations of these systems,
only the components $\Psi_\tau$ satisfying $\um\tau\leq l$ appear.
This is in contrast with the equations of a general system
\eqref{modeleq}, which may be written in the components indexed by the
set \eqref{finiteq}.

\begin{lemma}
\label{eform}
  The system of equations \eqref{modeleq} corresponding to a test
  curve  $\Gammac\in\eref B_L$ is of the form
\begin{equation}
  \label{eqns}
 \sum_{\um\tau\leq l} \comb(\tau)\,u_{\tau}^{l}\,\Psi_\tau=0,\quad l=1,2\ddd d,
\end{equation}
where $u_{\tau}^{l}$, $\um\tau\leq l\leq d$, are some complex
coefficients.
 \end{lemma}
 \begin{rem}
 We will think of the complex numbers   $u_{\tau}^{l}$, $\um\tau\leq
 l\leq d$ as coordinates on $\homnewd$.
 \end{rem}
 We can formalize this simple point as follows: introduce a new
 filtered vector space $\newd$:
\begin{equation}\label{defy}
 \newd=\bigoplus_{\um\tau\leq d}\C e_\tau\supset\bigoplus_{\um\tau\leq
  d-1}\C e_\tau\supset\dots\supset\C e_{2}\oplus\C e_1^2\oplus\C e_1\supset
\C e_1;
\end{equation}
the notation is motivated by the fact that $\newd$ is a truncation
of $\symdot$.  Now recall the notation $\hom^\triangle(\cdot,\cdot)$ for
filtration preserving linear maps, and introduce the following analog of $\tfnk$:
\begin{equation}
  \label{defe}
  \eqns = \{\epsilon\in\homnewd;\;
\ker(\epsilon)=0\}.
\end{equation}
With this notation Lemma \ref{eform} says that $\phip(\eref
B_L)\subset\eqns$. This statement may be globalized as follows.
Observe that the space $\homnewd$ is a left-right
representation of the group $B_L\times B_R$, and consider the
commutative diagram
\begin{diagram}[LaTeXeqno,labelstyle=\textstyle]
  \label{diaghom}
&  &  \hofi \\
& \ruTo^{\phip} & \uTo_\kappa \\
\hom(\C^d_L,\C^n) & \rTo  & \hom(\C^d_L,\C^n)
\times_{B_L}\homnewd\\
\end{diagram}
where
\begin{itemize}
\item $\phip$ defined in
  \eqref{mapson},
\item
the horizontal arrow is the correspondence
$\Gammac\mapsto(\Gammac,\epsilon_{\rf})$,
\item $\kappa$ is obtained by composing the linear map
  $\C_R^d\to\newd$ with the substitution $\C_L^d\to\C^n$.
\end{itemize}
A key point here is that we represent the set of systems
\eqref{eqnoverreffl} as an orbit of the $B_L$-action on
$\homnewd$.

\begin{prop}\label{eqnstilquot}
 \begin{enumerate}
 \item The open subset $\eqns\subset\homnewd$ is
   invariant under the left-right action of $B_L\times B_R$.
 \item The quotient $\eqnstil=\eqns/B_R$ is a smooth, compact variety
   endowed with a left action of $B_L$.
 \item The map $\kappa$ in diagram \eqref{diaghom} is
   $B_R$-equivariant, and induces a map
   $\kappat:\ind{\eqnstil}\rightarrow\fnk$.
 \item The horizontal map in diagram \ref{diaghom} induces an
   algebraic embedding
 \[ \phie: \homreg/H_L \rightarrow \ind{\eqnstil},
\]
such that the restriction of the map $\phit$ to $\homreg/H_L\subset
\qmr$ factorizes as $\kappat\circ\phie$ (cf. diagram
\eqref{diagfirst}).
 \end{enumerate}
\end{prop}
\begin{proof}
  The first and the third statements are obvious, while the second may
  be proved the same way as Lemma \ref{fdnquot}.

  For proving the last statement, observe that $\eqnstil$ is naturally
  a subvariety of $\fnk$, and $\phip(\eref B_L)\subset\eqns$ implies
  that $\phif(\eref B_L/H_L)\subset\eqnstil\subset\fnk$.  Moreover,
  denoting by $\phi$ the restriction of $\phif$ to $\eref B_L/H_L$, it
  is clear that this injective map is an embedding, since it is an
  orbit of a point under a Lie group action.

  Now inducing over $\fld$, we obtain the embedding
  \[
  \phie:\ind{B_L/H_L} \hookrightarrow \ind{\eqnstil}.
  \]
  The second half of the last statement follows from the
  construction of $\phie$.
\end{proof}

\begin{cor}\label{stabilizer}
 Let $\buref\in\eqnstil$ be the reference point
  $\pre(\epsilon_\rf)$, where $\pre:\eqns \to \eqnstil$ is the projection.
  The stabilizer of the $B_L$-action on $\eqnstil$ of
  the point $\buref$ is the subgroup $H_L\subset B_L$.
\end{cor}

Combining the results of Proposition \ref{eqnstilquot} with diagram
\eqref{diagfirst}, we arrive at the following picture:
\begin{diagram}[LaTeXeqno,labelstyle=\textstyle]
\label{diagmain}
\eref B_L/H_L & \rInto^{\phi} &  \eqnstil & & & & \solfnk   \\
\dInto &  & \dInto & & & \ldTo& \dTo\\
\hnd/H_L & \rInto^{\phie} & \feq &  \rTo^{\kappat} & \fnk& & \mapd nk   \\
& \rdTo(1,2)_{\pfl}  \ldTo(1,2)^{\taufl} & & & & \luTo &\dOnto^s  \\
  & \fld &  & & & & V^*\tensor\C^k
\end{diagram}

Now we are ready to formulate our model in its final form.
  \begin{itemize}
  \item
Consider the fibered product $V=\hnd \times_{B_L} \eqns \times_{B_R}
\C^d_R$, resulting in the double fibration
\[
\fld \longleftarrow \ind\eqnstil \overset{\tau_V}\longleftarrow V
\]
where $\eqns$ is defined in \eqref{defe}, and $\eqnstil=\eqns/B_R$.
\item Let $\sol_{\eqnstil}=\tilde{\kappa}^*(\sol_{\fnk})$; then
  comparing the construction of the bundle $V$  given above
  with Lemma \ref{vident}, we see that  we can pull back the
 sequence from \eqref{exact} to an exact sequence over $\eqnstil$:
  \begin{diagram}\label{diageqnstil}
    0 & \rTo & \sol_{\eqnstil} & \rTo^{\ev_{\eqnstil}}& \mapd nk &
 \rTo^{s} &   V^*\tensor\C^k & \rTo &0\\
  \end{diagram}
\end{itemize}

We have the following analog of \eqref{fmodel1}.

  \begin{prop}
  \label{themodelprop}\mbox{} Let $\buref\in\eqnstil$ be the point
  corresponding to the system~\eqref{reference_eq} (cf. Corollary
  \ref{stabilizer}). Then the orbit $B_L\buref$ is an irreducible
  $B_L$-invariant subvariety in $\eqnstil$ of dimension $d\choose2$,
  and  $\ev_{\eqnstil}$ establishes a degree-1 map
  \[
\tau_{\eqnstil}^{-1}\left(\ind{\overline{B_L\buref}}\right)\longrightarrow
\overline{\Upsilon}_d=
  \overline{\Theta}_d.\]

 \end{prop}

\begin{proof}
The first half of the statement follows from Corollary
\ref{stabilizer} once we note that the image of the map $\phi$ is
exactly $B_L\buref$. For the second half consider the following
facts:
\begin{itemize}
\item The evaluation map $\ev_{\eqnstil}$ is proper.
\item According to Proposition \ref{eqnstilquot} (4), we have
  $\phit=\kappat\circ \phie$ on the Zariski open part $\homreg/H_L$ in
  $\jddiff$.
\item The closure of $\Theta_d$ coincides with that of $\Upsilon_d$.
\end{itemize}
Now the statement follows from our previous ``model'' construction,
\eqref{fmodel}.
\end{proof}

\section{Application of the localization formulas}
\label{sec:apploc}

Recall that our aim is the computation of the equivariant Poincar\'e
dual $\epd{\overline{\Theta}_d}$, where the subvariety
$\Theta_d\subset\mapd nk$ represents the $A_d$-singularity (cf.
\S~\ref{sec:basicsing}). The symmetry group of the problem is the
product of matrix groups $\gl n\times\gl k$; the respective
subgroups of diagonal matrices are $T_n$ with weights
$(\lambda_1,\dots,\lambda_n)$ and $T_k$ with weights
$(\theta_1,\dots,\theta_k)$, hence $\epd{\overline{\Theta}_d}$ is a bisymmetric
polynomial in these two sets of variables.

In this section,  we apply the localization techniques of
\S\ref{sec:loctech} to the computation of $\epd{\overline{\Theta}_d}$ using
the model described in  \S\ref{subsec:fibr}. As our model is a
double fibration, the application of the localization formula is a
2-step process.

Before we proceed, we set the following {\bf convention}: when
describing the action of $B_L$ on the $B_R$-quotient $\eqnstil$, we
will revert to the notation $B_d$, since here there is only one copy
of the Borel group is acting.

\subsection{Localization in $\fld$}
\label{subsec:locind}

The model of Proposition \ref{themodelprop} is an equivariant
fibration over the smooth homogeneous space $\fld$, hence, in this
case, we can use Proposition \ref{locsingfibr}  (cf.
\S~\ref{sec:nonlinfib}), which applies when the fibers of $S$ are not necessarily
linear and smooth. The result of our calculation is Proposition
\ref{flagcor} below.

The data needed for formula \eqref{loclin}) is
\begin{itemize}
\item the fixed point set of the $T_n$-action on $\fld$,
\item the weights of this action on the tangents spaces $\TT_p\fld$ at
  these fixed points,
\item the equivariant Poincar\'e duals of the fibers at these fixed
  points. 
\end{itemize}
The following general statement will be helpful in organizing our
fixed point data. Its proof is straightforward and will be omitted.
\begin{lemma}
\label{snaction} Assume that the torus action in Proposition
\ref{locsmooth} is obtained by a restriction of a $\gl n$-action to
its subgroup of diagonal matrices $T_n$. Then the Weyl group of
permutation matrices $\sg n$ acts on $M^{T_n}$, and we have
\[   \epd{S_{\!\sigma\cdot p},W} = \sigma\cdot \epd{S_{\!p},W}
\text{ and } \euler^{T_n}(\TT_{\!\sigma\cdot p}M) =
\sigma\cdot\euler^{T_n}(\TT_{\!p}M).
\]
for all $\sigma\in\sg n$ and $p\in M_{T_n}$.
\end{lemma}

Our situation is fortunate in the sense that the action of $\sg n$
on the fixed point set is transitive.  Indeed, the fixed point set
$\fld^{T_n}$ is the set of partial flags obtained from sequences of
$d$ elements of the basis $(e_1\ddd e_n)$ of $\C^n$; in particular,
$|\fld^{T_n}|=n(n-1)\dots(n-d+1)$.

Recall  the notation $\reffl$ for the reference flag associated to
the sequence $(e_1\ddd e_d)$. The stabilizer subgroup of $\reffl$ in
$\sg n$ is the subgroup $\sg{n-d}$ permuting the numbers starting
with $d+1$, and  the map $\sigma\mapsto \sigma\cdot\reffl$ induces a
bijection between $\fld^{T_n}$ and the quotient $\sg n/\sg{n-d}$.

According to Lemma \ref{snaction}, it is sufficient to
compute the equivariant Poincar\'e dual of the fiber and the weights
of the tangent space at the reference flag $\reffl$.
The weights of $\TT_{\!\reffl}\fld$ are well-known:
\[  \{\lambda_i-\lambda_m; \;1\leq m\leq
d,\,m<i\leq n\};
\]
the weights at the other fixed points are obtained by applying the
corresponding permutation this set.

The numerators of the summands of \eqref{loclin} in our case are
much harder to compute, although, thanks to Lemma \ref{snaction}, it
is suffices to compute the numerator for the fixed point $\reffl$.
The situation over $\reffl$ is reflected in the following diagram:
  \begin{diagram}[LaTeXeqno,labelstyle=\textstyle]\label{diageqns}
    \sol_{\eqnstil} & \rTo^{\ev_{{\eqnstil}}}& \mapd nk &
 \rTo^{s} &   V^*\tensor\C^k \\
 & \rdTo_{\tau_{{\eqnstil}}} & & \ldTo & \\
\O=\overline{B_d\buref}& \rInto & {\eqnstil}  & &
  \end{diagram}

The fiber of our model \eqref{fmodel} over the fixed point $\reffl$
is the set $\tau^{-1}_{\eqnstil}(\O)$, where we introduced the
notation $\O$ for the closure of the $B_d$-orbit of $\buref$. Using
this notation, we can write the numerator of the term corresponding
to $\reffl$ in the sum \eqref{loclin} as follows:
\begin{equation}
  \label{epdfiber}
 \mathrm{eP}\left[\ev_{{\eqnstil}}\left(\tau^{-1}_{\eqnstil}(\O)\right),\mapd
   nk\right].
\end{equation}
Recall that this is a polynomial in two sets of variables:
$\lambdab=(\lambda_1\ddd\lambda_n)$ and
$\thetab=(\theta_1\ddd\theta_k)$. Since $\O$ is invariant under
$B_d$ only, this polynomial is not necessarily symmetric in the
$\lambda$s. The following statement is straightforward.
\begin{lemma}
  \label{onlyd}
  The equivariant Poincar\'e dual \eqref{epdfiber} does not depend on
  the last $n-d$ basic $\lambda$-weights:
  $\lambda_{d+1},\dots,\lambda_n$.
\end{lemma}
\begin{proof}
Indeed, recall that
$\ev_{{\eqnstil}}\tau^{-1}_{\eqnstil}(B_d\buref)$ consists of all
possible solutions of the systems of equations of the form
$B_L\epsilon_{\rf}$.  We wrote down these systems explicitly in
(\ref{eqnoverreffl}), and saw in \S~\ref{subsec:fibr} that all these
systems are in $\eqns$.  The systems of equations in $\eqns$,
however, impose conditions only on those components of $\Psi$ which
do not have indices higher than $d$, and this implies the statement
of the Lemma.
\end{proof}

As a consequence of Lemma \ref{onlyd}, the equivariant Poincar\'e
dual \eqref{epdfiber} may be considered as being taken with respect
to the group $T_d\times T_k$, which has weights $\bz=(z_1,\dots,
z_d)$ and $\thetab=(\theta_1,\dots,\theta_k)$.

Putting together Lemmas \ref{snaction} and \ref{onlyd} and the
description of the fixed point set $\fld^{T_d}$ given above, we
arrive at the following form of  \eqref{loclin} applied to our
situation:
\begin{prop}
  \label{flagcor}
We have
\begin{equation} \label{flagloc}
\epd{\Theta_d}=\sum_{\sigma\in\sg n/\sg{n-d}}
\frac{\qfl(\lambda_{\sigma\cdot1}\ddd\lambda_{\sigma\cdot
d},\thetab)} {\prod_{1\leq m\leq
d}\prod_{i=m+1}^n(\lambda_{\sigma\cdot
    i}-\lambda_{\sigma\cdot m})},
\end{equation}
where
\begin{equation}
  \label{defqfl}
  \qfl(\bz,\thetab) =
  \mathrm{eP}\left[\ev_{{\eqnstil}}\left(\tau^{-1}_{\eqnstil}(\O)\right),\mapd
    nk\right]_{T_d\times T_k}.
\end{equation}
\end{prop}

\subsection{Residue formula for the cohomology pairings of $\fld$}
\label{sec:pairings}

Usually, formulas such as \eqref{flagloc} are difficult to use: they
have the form of a finite sum of rational functions, and only after
adding up the terms of this sum and performing some cancellations do
we obtain a polynomial. These computations often obscure the
underlying structures, and they are rather unwieldy as the number of
terms of the sum grows very quickly with $n$ and $d$.

In this paragraph, we derive an efficient residue formula for the
right hand side of \eqref{flagloc}. While the geometric meaning of
this formula is not entirely clear, our summation procedure yields
an effective, ``truly'' localized formula; by this we mean that for
its evaluation one only needs to know the behavior of a certain
function at a single point, rather than at a large, albeit finite
number of points.

To describe this formula, we will need the notion of an {\em iterated
  residue} (cf. e.g. \cite{szenes}) at infinity.  Let
$\omega_1,\dots,\omega_N$ be affine linear forms on $\C^d$; denoting
the coordinates by $z_1,\ldots, z_d$, this means that we can write
$\omega_i=a_i^0+a_i^1z_1+\ldots + a_i^dz_d$. We will use the shorthand
$h(\bz)$ for a function $h(z_1\ddd z_d)$, and $\dbz$ for the
holomorphic $d$-form $dz_1\wedge\dots\wedge dz_d$. Now, let $h(\bz)$
be entire function, and define the {\em iterated residue at infinity}
as follows:
\begin{equation}
  \label{defresinf}
 \res_{z_1=\infty}\ldots \res_{z_d=\infty}\frac{h(\bz)\,\dbz}{\prod_{i=1}^N\omega_i}
  \overset{\mathrm{def}}=\left(\frac1{2\pi i}\right)^d
\int_{|z_1|=R_1}\ldots
\int_{|z_d|=R_d}\frac{h(\bz)\,\dbz}{\prod_{i=1}^N\omega_i},
 \end{equation}
 where $1\ll R_1 \ll \ldots \ll R_d$. The torus $\{|z_m|=R_m;\;m=1\ddd
 d\}$ is oriented in such a way that $\res_{z_1=\infty}\ldots
 \res_{z_d=\infty}\dbz/(z_1\cdots z_d)=(-1)^d$.

 We will also use the following simplified notation:
 \[ \sires \overset{\mathrm{def}}=\ires.
 \]

In practice, the iterated residue \ref{defresinf} may be computed
using the following {\bf algorithm}: for each $i$, use the expansion
 \begin{equation}
   \label{omegaexp}
 \frac1{\omega_i}=\sum_{j=0}^\infty(-1)^j\frac{(a^{0}_i+a^1_iz_1+\ldots
   +a_{i}^{q(i)-1}z_{q(i)-1})^j}{(a_i^{q(i)}z_{q(i)})^{j+1}},
   \end{equation}
   where $q(i)$ is the largest value of $m$ for which $a_i^m\neq0$,
   then multiply the product of these expressions with $(-1)^dh(z_1\ddd
   z_d)$, and then take the coefficient of $z_1^{-1} \ldots z_d^{-1}$
   in the resulting Laurent series.

 We have the following {\em iterated residue theorem}.
\begin{prop} \label{propflag}
  \label{flagresidue} For a polynomial $Q(\bz)$ on $\C^d$, we have
\begin{equation}
  \label{flagres}
\sum_{\sigma\in\sg n/\sg{n-d}}
\frac{Q(\lambda_{\sigma\cdot1}\ddd\lambda_{\sigma\cdot d})}
{\prod_{1\leq m\leq d}\prod_{i=m+1}^n(\lambda_{\sigma\cdot
    i}-\lambda_{\sigma\cdot m})}=\sires
\frac{\prod_{1\leq m<l\leq d}(z_m-z_l)\,Q(\bz)\dbz}
{\prod_{l=1}^d\prod_{i=1}^n(\lambda_i-z_l)}
\end{equation}
\end{prop}
\begin{proof}
  We compute the iterated residue \eqref{flagres} using the Residue
  Theorem on the projective line $\C\cup\{\infty\}$.  The first
  residue, which is taken with respect to $z_d$, is a contour
  integral, whose value is minus the sum of the $z_d$-residues of the
  form in \eqref{flagres}. These poles are at $z_d=\lambda_j$,
  $j=1\ddd n$, and after canceling the signs that arise, we obtain the
  following expression for the right hand side of  \eqref{flagres}:
\[
\sum_{j=1}^n \frac{\prod_{1\leq m<l\leq
    d-1}(z_m-z_l)\,\prod_{l=1}^{d-1}(z_l-\lambda_j)\,Q(z_1\ddd
  z_{d-1},\lambda_j)\;dz_1\dots dz_{d-1}}
{\prod_{l=1}^{d-1}\prod_{i=1}^n(\lambda_i-z_l)\prod_{i\neq
    j}^n(\lambda_i-\lambda_j)}.
\]
After cancellation and exchanging the sum and the residue operation,
at the next step, we have
\[
(-1)^{d-1}\sum_{j=1}^n\res_{z_{d-1}=\infty} \frac{\prod_{1\leq
m<l\leq
    d-1}(z_m-z_l)\,Q(z_1\ddd z_{d-1},\lambda_j)\;dz_1\dots
  dz_{d-1}} {\prod_{i\neq j}^n
  \left((\lambda_i-\lambda_j)\prod_{l=1}^{d-1}(\lambda_i-z_l)\right)}.
\]
Now we again apply the Residue Theorem, with the only difference that
now the pole $z_{d-1}=\lambda_j$ has been eliminated. As a result,
after converting the second residue to a sum, we obtain
\[
(-1)^{2d-3}\sum_{j=1}^n\sum_{s=1,\,s\neq j}^n \frac{\prod_{1\leq
m<l\leq
    d-2}(z_l-z_m)\,Q(z_1\ddd
  z_{d-2},\lambda_s,\lambda_j)\;dz_1\dots dz_{d-2}}
{(\lambda_s-\lambda_j)\prod_{i\neq j,s}^n
  \left((\lambda_i-\lambda_j)(\lambda_i-\lambda_s)\prod_{l=1}^{d-1}(\lambda_i-z_l)\right)}.
\]
Iterating this process, we arrive at a sum very similar to
(\ref{flagloc}). The difference between the two sums will be the
sign: $(-1)^{d(d-1)/2}$, and that the $d(d-1)/2$ factors of the form
$(\lambda_{\sigma(i)}-\lambda_{\sigma(m)})$ with $1\le m<i\le d$ in
the denominator will have opposite signs. These two differences
cancel each other, and this completes the proof.
\end{proof}
\begin{rem}
  Changing the order of the variables in iterated residues, usually,
  changes the result. In this case, however, because all the poles are
  normal crossing, formula \eqref{flagres} remains true no matter in
  what order we take the iterated residues.
\end{rem}

\subsection{Localization in the fiber}
\label{subsec:locfib} Combining Proposition \ref{flagcor} with
Proposition \ref{propflag}, we arrive at the formula
\begin{equation}
  \label{epqres}
 \epd{\Thetab_d,\mapd nk}=\sires
\frac{\prod_{1\leq m<l\leq d}(z_m-z_l)\,\qfl(\bz,\thetab)\dbz}
{\prod_{l=1}^d\prod_{i=1}^n(\lambda_i-z_l)}.
\end{equation}
The ``only'' unknown here is the polynomial $\qfl(\bz,\thetab)$
defined in \eqref{defqfl}, and, therefore, we now turn to its
computation.

Let us briefly review the construction of $\qfl(\bz,\thetab)$ (cf.
diagram \eqref{diageqns} and Proposition \ref{flagcor}). This
polynomial is an equivariant Poincar\'e dual taken with respect to
the group $T_d\times T_k$, which has weights $(z_1,\dots, z_d)$ and
$(\theta_1,\dots,\theta_k)$.  Consider the $B_L\times B_R$-module
$\homnewd$, and endow it with coordinates $u^l_\tau\in\homnewd^*$,
indexed by pairs $(\tau,l)\in\Pi\times\Z_{>0}$ satisfying
$\um\tau\leq l\leq d$. We will consider the dual space spanned by
these coordinates as carrying a {\em right} action of $T_d\times
T_k$; accordingly,
\begin{equation}
  \label{uweights}
   \text{the weight of }
u^l_\tau\;=\;(z_{i_1}+z_{i_2}+\dots+z_{i_m},\theta_l),\; \text{where
}\tau=[i_1,i_2\ddd i_m].
\end{equation}

For each nondegenerate system $\epsilon\in\eqns\subset\homnewd$ we
denote the image $\pre(\epsilon)$ in the quotient
$\pre:\eqns\to\eqnstil=\eqns/B_R$ by $\teps$; in particular, we have
a reference point $\buref\in\eqnstil$ corresponding to the system
$\epsref$ given by
\begin{equation}
  \label{defpref}
  \cu(\epsilon_\rf)=
\begin{cases}
  1,\text{ if }\um\pi=l \\
  0,\text{ otherwise.}
\end{cases}
\end{equation}
The stabilizer subgroup of $\buref\in\eqnstil$ under the
$B_d$-action is  a $d$-dimensional subgroup $H_d\subset B_d$, hence
the orbit $B_d\buref\subset\eqnstil$ is a subvariety of dimension
$d(d-1)/2$; we denoted the closure of this subvariety by $\O$.

Next, consider the vector bundle
\[ V=\eqns\times_{B_R}\C_R^d\longrightarrow\eqnstil=\eqns/B_R\]
associated to the standard representation of $B_R$, and the
$T_d\times T_k$-equivariant linear bundle map from a trivial bundle
\[   s: \eqnstil\times\mapd nk \longrightarrow V^*\tensor\C^k
\]
defined by the natural composition \eqref{canpairing0}. Then,
according to Proposition \ref{themodelprop}, the polynomial
$\qfl(\bz,\thetab)$ is the {\em equivariant Poincar\'e dual in $\mapd
  nk$ of the union of the vector spaces $\ker(s)$ lying over
  $\O\subset\eqnstil$} (cf. \eqref{defqfl}).

While the variety $\O$ is highly singular, the set of $T_d$-fixed
points of $\O$ is finite -- as we will see shortly --  and hence we
can apply here the localization principle based on Rossmann's
integration formula: Proposition \ref{locepdsing}. The result is:
\begin{equation}
  \label{qfl1}
\qfl(\bz,\thetab) =\sum_{p\in \O^{T_d}} \frac{\euler^{T_d\times
T_k}(V_p^*\tensor
  \C^k)\,\emu_p[\O,\eqnstil]}
{\euler^{T_d\times T_k}(\TT_p\eqnstil)}.
  \end{equation}

Our task thus has reduced to the identification and computation of
the objects in this formula. These are:
\begin{itemize}
\item The set $\O^{T_d}$ of $T_d$-fixed points in $\O\subset\eqnstil$,
\item the weights of the $T_d$-action on the fibers $V_p$ for
  $p\in\O^{T_d}$,
\item the weights of the $T_d$-action on the tangent spaces
  $\TT_p\eqnstil$ for $p\in\O^{T_d}$,
\item the equivariant multiplicities of $\O$ in $\eqnstil$ at each
  fixed point $p\in\O^{T_d}$.
\end{itemize}

  The most immediate problem we face is that we do not have an
  effective description of the set $\O^{T_d}$ of $T_d$-fixed points in
  $\O$. There is a formal way
  around this: we replace the fixed point set $\O^{T_d}$ with the
  larger set $\eqnstil^{T_d}$, and define the equivariant multiplicity
  $\emu_p[\O,\eqnstil]$ to be zero in the case when $p\in
  \eqnstil^{T_d}\setminus \O^{T_d}$.

The set of fixed points $\eqnstil^{T_d}$ is fairly easy to
determine: these fixed points are given by those nondegenerate
systems $\epsilon\in\eqns\subset\homnewd$ for which the tensors
$\epsilon(e_m)\in\newd$, $m=1,\dots,d$ are of pure $T_d$-weight.
These, in turn, may be enumerated as follows.
\begin{defi}\label{admissible}
  We will call a sequence of partitions
  $\bipi=(\pi_1\ddd\pi_d)\in\Pi^{\times d}$ {\em admissible } if
\begin{enumerate}
\item $\um {\pi_l}\leq l$ for $l=1\ddd d$ and
\item $\pi_l\neq\pi_m$ for $1\leq l\neq m\leq d$. \end{enumerate}
We will denote the set of admissible sequences of length $d$ by
$\Bipi$; we also introduce the numerical characteristic:
\[  \mathrm{defect}(\bipi)=\sum_{l=1}^d(l-\um{\pi_l}).
\]
\end{defi}

As an example, we list the admissible sequences in the case $d=3$:
 \begin{multline*}
\boldsymbol{\Pi}_3=\{([1],[2],[3]),\, ([1],[2],[1,2]),\,
([1],[2],[1,1]),\,([1],[2],[1,1,1]) \\
([1],[1,1],[3]),\, ([1],[1,1],[1,1,1]), \, ([1],[1,1],[2]),\,
([1],[1,1],[1,2])\};
\end{multline*}

For $\bipi=(\pi_1\ddd\pi_d)\in\Bipi$ introduce the system
$\epsilon_\bipi$ given by
\begin{equation}\label{ubipi}
 u_{\tau}^{l}(\epsilon_\bipi)
=
\begin{cases}
  1 \text{ if } \tau=\pi_l,\\
0\text{ otherwise.}
\end{cases}
\end{equation}
As usual, the point corresponding to $\epsilon_\bipi$ in $\eqnstil$
will be denoted by $\teps_\bipi=\pre(\epsilon_\bipi)$.

The following statement follows from the definitions.
\begin{lemma}
  \label{eqnstilfixed}
  \begin{itemize}
  \item
  The correspondence $\bipi\mapsto\teps_\bipi$ establishes a
  bijection between the set $\Bipi$ of admissible
  sequences of partitions and the fixed point set $\eqnstil^{T_d}$.
\item For $\tau\in\Pi$, and an integer $i$, denote by $\mult(i,\tau)$
  the number of times $i$ occurs in $\tau$, and let
  $z_\tau=\sum_{i\in\tau}\mathrm{mult}(i,\tau)\,z_i$. Then, given an
  admissible sequence $\bipi\in\Bipi$, the weights of the $T_d$-action
  on the fiber of $V$ at the fixed point $\teps_\bipi$ are
 \[  z_{\pi_1}\ddd z_{\pi_d}.\]
\end{itemize}
\end{lemma}

\begin{cor}
  \label{dkweights}
  The weights of the $T_d\times T_k$ action on fiber
  $V^*_{\teps_{\bipi}}\tensor\C^k$ are
\[ \{\theta_j-z_{\pi_m};\;m=1\ddd d,\,j=1\ddd k\}.\]
\end{cor}

Next we turn to the 3rd item on our list: the weights of the
$T_d$-action on tangent space of $\eqnstil$ at the fixed points
$\teps_\bipi$;  we will use the simplified notation
$\TT_\bipi\eqnstil$ for this tangent space. To compute the answer,
it will be convenient to linearize the action near $\teps_\bipi$.

\begin{defi}
  \label{defnpi} For each $\bipi=(\pi_1,\dots,\pi_d)\in\Bipi$
  introduce the affine-linear subspace $\NN_\bipi\subset \hofi$ given
  by
\[\NN_{\bipi}=\left\{\epsilon \in \hofi;\;u_{\pi_l}^m(\epsilon)=
\begin{cases}
  1\,\text{ if }m=l\\
  0\,\text{ if }m> l
\end{cases}
  \text{ for } 1\le l \le d \right\};\]
Also, for  $\bipi \in \Bipi$ introduce the map
\[ \abi:\hofi \to \mathrm{Mat}^{d\times d}
\]
which associates to each system $\epsilon$ its $d\times d$ minor
corresponding to the sequence of partitions
$\bipi=(\pi_1,\dots,\pi_d)$.
\end{defi}

A few comments are in order.  First, we can rewrite the above
definition of $\NN_\bipi$ as follows:
\begin{equation}\label{nbipiother}
\NN_\bipi=\left\{ \epsilon \in \hofi;\;\abi(\epsilon)\in U_{-}
\right\}
\end{equation}
where $U_-$ is the subgroup of lower-triangular $d\times d$ matrices
with $1$s on the diagonal; this way it is apparent that
$\NN_\bipi\subset\eqns$.

Also, observe that $\epsilon_{\bipi}\in \NN_\bipi$, and considering
this special point to be the origin, we may think of $\NN_\bipi$ as
a {\em linear} space. Then $\NN_\bipi$ is endowed with a natural set
of coordinates:
\begin{equation}\label{npicoordinates}
   \uh_{\tau|\bipi}^l=u_\tau^l|\NN_\bipi, \,\um\tau\leq l\leq
d,\,\tau\neq\pi_1\ddd\pi_l.
\end{equation}
\begin{prop}
  \label{cover}
  Let $\bipi\in\Bipi$ be an admissible sequence of partitions. Then
  \begin{enumerate}
  \item the restriction of the projection $\pre:\eqns\to\eqnstil$ to
    $\NN_\bipi$ is an embedding and the collection
    $\left\{\pre(\NN_\bipi);\; \bipi \in \Bipi\right\}$ forms an open
    cover of $\eqnstil$.
  \item for {\em any} $\bipi\in \Bipi$, the image
    $\pre(\NN_\bipi)\subset\eqnstil$ is $T_d$-invariant, and the
    induced $T_d$-action on $\NN_\bipi$ is linear and diagonal with
    respect to the coordinates \eqref{npicoordinates}.
    Considering $T_d$ as acting on the {\em right} on these
    coordinates,
\begin{equation}
  \label{defwts}
  \text{ the weight of }\;\uh_{\tau|\bipi}^l\; =\;z_\tau-z_{\pi_l}.
\end{equation}
  \item If $\mathrm{defect}(\bipi)=0$, then
    $\pre(\NN_\bipi)\subset\eqnstil$ is $B_d$- invariant.
  \end{enumerate}
\end{prop}
\begin{rem}\label{justoneborel}
  We will denote by $T_\bipi$ and $B_\bipi$ the actions of $T_d$ and
  $B_d$ induced on $N_\bipi$ by the embedding $\pre$.
\end{rem}
\begin{proof}
  We first show that $\cup\left\{\pre(\NN_\bipi);\; \bipi \in
    \Bipi\right\}=\eqnstil$. This means that for an arbitrary element
  $\epsilon \in \eqns$, we have to find an admissible partition
  $\bipi\in \Bipi$ and an upper-triangular matrix
  $b_R=b_R(\epsilon,\bipi) \in B_R$ such that $\epsilon \cdot b_R \in
  \NN_{\bipi}$.  This can be done by elementary column operations:
  consider $\epsilon$ as a $\dim (\newd)\times d$ matrix whose columns
  are linearly independent, and whose rows are indexed by
  partitions. The only nonzero entry in the first column corresponds
  to the trivial partition $[1]$, hence we can multiply the first
  column by a constant to rescale this entry to 1, and then annihilate
  all other entries in the same row by adding multiples of the first
  column to the others.  Next, since $\epsilon$ is nonsingular, we can
  pick a nonzero entry in the second column of the resulting matrix --
  this entry will correspond to a partition $\pi_2$ -- and, again,
  using column operations, we annihilate all entries in this row
  starting form column 3 and so on. Continuing this process, we obtain
  an admissible $\bipi=(\pi_1,\ldots ,\pi_d)$, and the described
  sequence of column operations produces an upper-triangular $b_R\in
  B_R$ such that $\epsilon\cdot b_R\in\NN_\bipi$.

  The process described above finds an appropriate $\bipi \in \Bipi$
  for each $\epsilon$, and brings $\abi(\epsilon)$ to lower-triangular
  form. Moreover, if $\pre(\epsilon_1)=\pre(\epsilon_2)$ for
  $\epsilon_1,\epsilon_2 \in \NN_\bipi$, then $\epsilon_1 \cdot
  b_R=\epsilon_2$ for some $b_R \in B_R$, and therefore
  $\abi(\epsilon_1)\cdot b_R=\abi(\epsilon_2)$. Since
  $\abi(\epsilon_1),\abi(\epsilon_2)$ are lower-triangular with $1$s
  on the diagonal and $B_R$ is upper-triangular, this can only happen
  when $b_R$ is the unit matrix, so $\epsilon_1=\epsilon_2$. This
  proves that $\pre$ is injective on $\NN_\bipi$, hence the
  restriction $\pre|\NN_\bipi$ is an embedding.

To approach statements (2) and (3), we  write down the action
of $B_d$ on $\eqnstil$ in the  chart $\NN_\bipi$.  Recall that the
multiplication map $U_-\times B_d \to \gl d$ is injective. This
allows us to define the $B_d$-component $a^B$ for an element $a\in
U_-B_d$; in particular, for any such $a$, we have $a\cdot
(a^B)^{-1}\in U_-$. Then, for $b\in B_d$ and $\epsilon\in\NN_\bipi$
we can define the partial  action:
\begin{equation}\label{action}
(b,\epsilon)\mapsto b_\bipi\epsilon=b_L\cdot \epsilon \cdot
(\abi(b_L\cdot\epsilon)^B)^{-1},
\end{equation}
which is valid if $\abi(b_L\cdot\epsilon)\in U_-B_d$.

Now consider the case when $b=t\in T_d$ is a diagonal matrix. In
this case, $\abi(b_L\cdot\epsilon)$ remains lower-triangular, with
the numbers $(t^{\pi_1},\dots, t^{\pi_d})$ on the diagonal, where
$t^\tau$ is the character of $T_d$ corresponding to the weight
$z_\tau$. This means that $\abi(b_L\cdot\epsilon)\in U_-B_d$, and
the Borel factor $\abi(b_L\cdot\epsilon)^B$ is the diagonal matrix
with these same entries:
 \begin{equation}
    \label{diag_action}
    \abi(b_L\cdot\epsilon)^B = \mathrm{diag}[t^{\pi_1},\ldots,
    t^{\pi_d}].
  \end{equation}
Note that this matrix is independent of $\epsilon$. Now statement
(2) follows easily.

Finally, to prove (3), observe that if $\mathrm{defect}(\bipi)=0$,
then the filtration-preserving property implies that
$\abi(\epsilon)$ is upper-triangular for any $\epsilon\in\hofi$.
Hence for $\epsilon\in\NN_\bipi$ the matrix
$\alpha_\bipi(\epsilon)$ is the identity matrix, and thus, using the
condition $\mathrm{defect}(\bipi)=0$ once again, we can conclude
that $\abi(b_L\cdot\epsilon)$ is upper-triangular with the numbers
$(t^{\pi_1},\dots, t^{\pi_d})$ on the diagonal, where $t$ is the
diagonal part of $b$. This means that
$\abi(b_L\cdot\epsilon)^B=\abi(b_L\cdot\epsilon)\in B_d$, which
implies statement (3).
\end{proof}

\begin{rem}
\label{quad} Clearly, $\abi(b_L\cdot\epsilon)$ depends linearly on
  $\epsilon$.  In the case $\mathrm{defect}(\bipi)=0$, we have
  $\abi(b_L\cdot\epsilon)^B=\abi(b_L\cdot\epsilon)$, and hence the
  action \eqref{action} of $B_\bipi$ on $\NN_\bipi$ is quadratic, not
  linear as the $T_\bipi$-action. When $\mathrm{defect}(\bipi)>0$, the
  action of $B_\bipi$ is not defined on the whole of $\NN_\bipi$.
\end{rem}

Proposition \ref{cover} provides us with a linearization of the
$T_d$-action on $\eqnstil$ near every fixed point.  This allows us
to compute equivariant multiplicities in \eqref{qfl1} using
\eqref{emult}. Indeed, if we introduce the notation
\begin{equation}
    \label{opdef}
    \O_\bipi\overset{\mathrm{def}}=(\pre|\NN_\bipi)^{-1}(\O)
  \end{equation}
for the part of $\O$ in the local chart $\NN_\bipi$, then we can
write
 \begin{equation}
      \label{compemult}
      \emu_{\teps_\bipi}[\O,\eqnstil]=\epd{\O_\bipi,\NN_\bipi}.
    \end{equation}

Next, we take a closer look at the set $\O_{\bipi}$.
\begin{lemma}
  \label{qpi}
For every  $\bipi\in\Bipi$, we have
\begin{equation}
  \label{opnn}
  \O_\bipi=\overline{B_L\epsref B_R}\cap\NN_\bipi.
\end{equation}
Moreover, $\epsref\in\NN_\bipi$ if and only if $\defect(\bipi)=0$,
and in this case $\O_\bipi=\overline{B_\bipi\epsref}$, where
$B_\bipi$ stands for the action \eqref{action}.
\end{lemma}
\begin{proof}
 By definition, $\O_\bipi=\overline{B_L\epsref B_R\cap N_\bipi}$, and
 hence \eqref{opnn} follows from the fact that $B_d$ acts properly on
 the right on $U_-B_d\subset\gl d$. The second statement then
 immediately follows from the comparison of \eqref{defpref} and
 Definition \ref{defnpi}.
\end{proof}

Let us take stock of our results so far.
Substituting the weights from Corollary \ref{dkweights} and
\eqref{defwts} into (\ref{qfl1}), and taking into consideration
\eqref{compemult}, we obtain:
\begin{equation}
  \label{rossflag}
  \qfl(\lambdab,\thetab)=\sum_{\bipi\in\Bipi}
  \frac{ \prod_{m=1}^d\prod_{j=1}^k(\theta_j-z_{\pi_m})\,Q_\bipi(z_1\ddd
    z_d)}
  {\displaystyle \prod_{l=1}^d\prod_{\um\tau\leq l}^{\tau\neq\pi_1\ddd\pi_l}(z_\tau-z_{\pi_{l}})},
\end{equation}
where
\begin{equation}
  \label{qbipidef}
Q_\bipi=\begin{cases}
\epd{(\O_\bipi,\NN_\bipi} \text{ if }\tepsi\in\O,\\
0\text{ if }\tepsi\notin\O.
\end{cases}
\end{equation}
Combining this formula with \eqref{flagres}, and arrive at our first
formula for $\epd{\tcl}$:
  \begin{equation} \label{fixedone}
\epd{\tcl}= \sires \frac{\prod_{m<l}(z_m-z_l)\,\dbz}
{\prod_{l=1}^d\prod_{i=1}^n(\lambda_i-z_l)} \sum_{\bipi\in\Bipi}
\frac{\prod_{m=1}^d\prod_{j=1}^k(\theta_j-z_{\pi_m})
\,Q_\bipi(\bz)}{\prod_{l=1}^d\prod\{(z_\tau-z_{\pi_{l}});\,\um\tau\leq
l,\,\tau\neq\pi_1\ddd\pi_l\}}
  \end{equation}

  Now observe that the sum here is finite, hence we are free to
  exchange the summation with the residue operation. Rearranging the
  formula accordingly, we arrive at the following statement.

    \begin{prop}\label{proptwo}
      For each admissible series $\bipi=(\pi_1\ddd\pi_d)$ of $d$
      partitions, introduce the polynomial $Q_\bipi(\bz)$ defined by
      \eqref{qbipidef}, then
  \begin{equation} \label{fixedtwo}
\epd{\Theta_d}= \sum_{\bipi\in\Bipi} \sires \frac{\displaystyle
Q_\bipi(\bz)\,\prod_{m<l}(z_m-z_l) }{\displaystyle
\prod_{l=1}^d\prod_{\um\tau\leq l}^{\tau\neq\pi_1\ddd\pi_l}}
\frac{\displaystyle  \prod_{m=1}^d\prod_{j=1}^k(\theta_j-z_{\pi_m})}
{\displaystyle  \prod_{l=1}^d\prod_{i=1}^n(\lambda_i-z_l)} \,\dbz.
  \end{equation}
    \end{prop}

    This formula has the pleasant feature that the three
    parameters of our problem, $n,k$ and $d$, enter in it in a
    separate manner. The first fraction here only depends on $d$, the
    denominator of the second only depends on $n$, and the numerator
    of this latter fraction controls the $k$-dependence, with some
    interference from the sequence $\bipi$.

    While this formula is a step forward, it is rather difficult to
    use in practice, since the number of terms and factors in it grows
    with $d$ as the the number of elements in $\Bipi$. Also, the known
    properties of Thom polynomials listed in Proposition \ref{collect}
    are not manifest in \eqref{fixedtwo}.

    In the next section, we will see that this formula goes through
    two dramatic simplifications, which will make it easy to evaluate
    it for small values of $d$.

Before proceeding, we present a schematic diagram of the main
objects of our constructions. We hope this will help the reader to
navigate among the various spaces we have introduced.

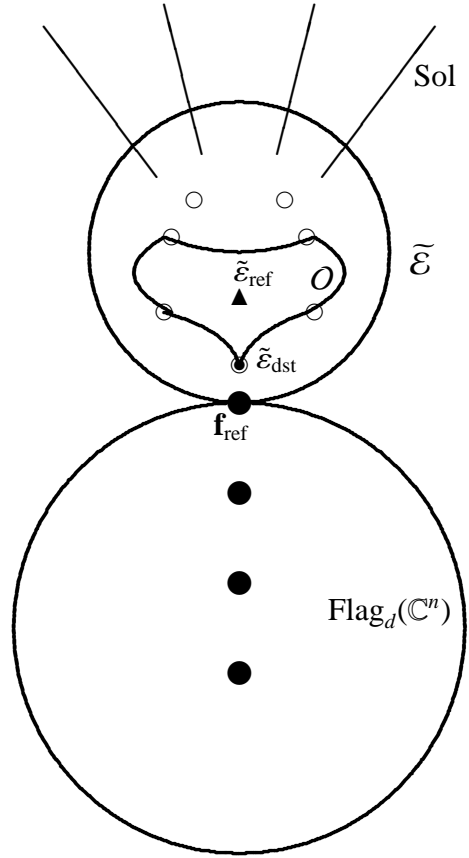
\begin{center}

\setlength{\unitlength}{1mm}
\begin{picture}(100,120)(-5,-5)
\linethickness{1pt}

\put(40,30){\bigcircle{60}} \put(40,80){\bigcircle{40}}

% % face features
\curve(31,72,26,77,30,82) \curve(30,82,40,80,50,82)
\curve(50,82,54,77,49,72) \curve(49,72,43,69,40,65)
\curve(40,65,37,69,30,72)

% empty circles in the head
% eyes
\put(34,87){\circle{2}} \put(46,87){\circle{2}}

\put(31,82){\circle{2}} \put(49,82){\circle{2}}

\put(50,72){\circle{2}} \put(30,72){\circle{2}}

% c point
\put(40,65){\circle*{1.5}} \put(40,65){\circle{2}}

% ref
% buttons
\put(40,60){\circle*{3}} \put(40,48){\circle*{3}}
\put(40,36){\circle*{3}} \put(40,24){\circle*{3}}

% hair
\thicklines \put(29,90){\line(-3,4){15}} \put(35,93){\line(-1,4){5}}
\put(45,93){\line(1,4){5}} \put(51,90){\line(3,4){15}}

% letters
\put(66,102){\makebox(0,0)[b]{$\sol$}}
\put(42,75.5){\makebox(0,0)[b]{$\teps_\rf$}}
\put(40,73){\makebox(0,0)[b]{$\blacktriangle$}}
\put(51,75){\makebox(0,0)[b]{$\O$}}
\put(45,64){\makebox(0,0)[b]{$\teps_\ds$}}
\put(39,55){\makebox(0,0)[b]{$\reffl$}}
\put(64,77){\makebox(0,0)[b]{$\eqnstil$}}
\put(60,30){\makebox(0,0)[b]{$\fld$}}

\end{picture}
\end{center}
Explanations:
\begin{itemize}
\item The lower circle is the flag variety $\fld$; the fat dots inside
  represent the $T_n$-fixed flags in $\fld$.
\item The upper circle is  $\eqnstil$, the fiber of the bundle
  $\ind\eqnstil$ over the reference flag $\reffl$. The small circles
  inside represent the $T_d$-fixed points in $\eqnstil$. One of these
  fixed points, $\teps_\ds\in\eqnstil$ will play an important role in what
  follows.
\item The region bounded by the curvy-linear pentagon represents the
  $B_d$-orbit of the reference point $\teps_\rf$, which is marked by a
  triangle. The closure of the orbit is $\O$; this is a singular
  subvariety of $\eqnstil$, which contains some of the fixed points of
  $\eqnstil$, but not all of them.
\item The straight lines on top are the linear solution spaces of the
  corresponding systems of equations in $\eqnstil$. The union of these
  solution spaces lying over those points of the fiber bundle
  $\ind\eqnstil$ which correspond to $\O$ form the closure of our
  singularity locus $\Theta_d$.
\end{itemize}

\section{Vanishing residues and the  main result}
\label{sec:vanishing}

The terms on the right hand side of formula \eqref{fixedtwo} are
enumerated by admissible sequences.  There is a simplest one among
these:
\begin{equation}\label{distinguished}
 \bipi_\ds= ([1],[2],\dots,[d]),
\end{equation}
which we will call {\em distinguished }. To avoid double indices,
below, we will use the simplified notation $Q_\ds$ instead of
$Q_{\bipi_\ds}$, and similarly $\teps_\ds,\NN_\ds,\O_\ds$, etc.

The following remarkable vanishing result holds.
\begin{prop}\label{vanishthm}
  Assume that $d\ll n\leq k$. Then all terms of the sum in
  \eqref{fixedtwo} vanish except for the term corresponding to the
  sequence of partitions $\bipi_\ds=([1],[2],\dots,[d])$. Hence,
  formula \eqref{fixedtwo} reduces to
\begin{equation}
  \label{fixedthree}
  \epd{\Thetab_d}=
\sires
\frac{Q_{\ds}(z_1\ddd z_n)\,\prod_{m<l}(z_m-z_l)\;\dbz}%dz_1\dots  dz_d}
{\prod_{l=1}^d\,\prod\{(z_\tau-z_l);\, \um\tau\leq l,\,|\tau|>1\}}
\frac{\prod_{l=1}^d\prod_{j=1}^k(\theta_j-z_l)}
{\prod_{l=1}^d\prod_{i=1}^n(\lambda_i-z_l)},
\end{equation}
where $Q_{\ds}=\epd{\O_\ds,\NN_\ds}$.
\end{prop}

Before turning to the proof, we make a few remarks.  First, note
that this simplification is dramatic: the number of terms in
\eqref{fixedtwo} grows exponentially with $d$, and of this sum now a
single term survives. This is fortunate, because computing all the
polynomials $Q_\bipi$, $\bipi\in\Bipi$ seems to be an insurmountable
task; at the moment, we do not even have an algorithm to determine
when $Q_\bipi=0$, i.e. when $\teps_{\bipi}\in\O$.

Our second observation is that after replacing in \eqref{fixedthree}
$z_l$ by $-z_l$, $l=1\ddd d$, we can rewrite \eqref{fixedthree} as
\begin{equation}
  \label{fixedfour}
  \epd{\Thetab_d}=
\sires \frac{(-1)^d\prod_{m<l}(z_m-z_l)\,Q_\ds(z_1\ddd z_n)}
{\prod_{l=1}^d\prod\left\{(z_\tau-z_l);\; \um\tau\leq
l,\,|\tau|>1\right\}}\; \prod_{l=1}^d
\ct\left(\frac1{z_l}\right)\,z_l^{k-n}\;dz_l,
\end{equation}
where $\ct(z)$ is the generating series of the relative Chern classes
introduced in \eqref{deftc}. Indeed, the denominator and the numerator
of the fraction in \eqref{fixedfour} are homogeneous polynomials of
the same degree, hence this substitution will leave the fraction
unchanged.  We thus obtain an explicit formula for the Thom polynomial
of the $A_d$-singularity in terms of the relative Chern classes.  This
is important, because the fact that \eqref{fixedfour} conforms to the
result of Thom-Damon, Proposition \ref{collect} (3), suggests that we
have the ``right'' formula.

Most of the present section will be taken up by the proof of
Proposition \ref{vanishthm}. In \S~\ref{sec:vanres}, we derive a
criterion  for the vanishing of iterated residues of the form
\eqref{defresinf}. Applying this criterion to the right hand side of
\eqref{fixedtwo} reduces Proposition \ref{vanishthm} to a statement
about the factors of the polynomials $Q_\bipi$, $\bipi\in\Bipi$:
Proposition \ref{critical}. According to Lemma \ref{divisible}, such
divisibility properties follow from the existence of relations of a
certain form in the ideal of the subvariety
$\O_\bipi\subset\NN_\bipi$.  We find a family of such relations in
\S~\ref{sec:homog} (see \eqref{relz}), and then convert the
condition in Lemma \ref{divisible} into a combinatorial condition on
$\bipi$ (cf. Lemma \ref{occurs}).  At the end of \S~\ref{sec:homog},
we show that if a sequence $\bipi$ does not satisfy this
combinatorial condition, then it is either $\bipi_\ds$ or
$\teps_\bipi\notin\O$, thus completing the proof of Proposition
\ref{vanishthm}.

Introduce the subset $\Bipio\subset\Bipi$ defined by
\begin{equation}
  \label{defpio}
\Bipio=\left\{\bipi\in\Bipi;\;\teps_{\bipi}\in\O\right\}.
\end{equation}
As we mentioned earlier, at the moment, we do not have an explicit
description of this set.  In the course of this proof, however, we
obtain a rather efficient, albeit incomplete criterion for a
sequence $\bipi\in\Bipi$ {\em not} to belong to $\Bipio$; we explain
this criterion in \S~\ref{sec:fixedpoints}.  Finally, in
\S\ref{sec:mainres}, we further simplify (\ref{fixedfour}), and
formulate our main result, Theorem \ref{dathm}.

Before embarking on this rather tortuous route, we give a few
examples below in \S\ref{sec:loc23}, which demonstrate the
localization formulas and the vanishing property explicitly. Note
that we devote the last chapter of the paper to the detailed study
of \eqref{fixedfour} for small values of $d$, and hence  the proofs
in \S\ref{sec:loc23} will be omitted.

\subsection{The localization formulas for $d=2,3$}
\label{sec:loc23}

The situation for $d=2$ and $3$ is simplified by the fact, that in
these cases the closure of the Borel-orbit
$\O=\overline{B_d\teps_\rf}\subset\eqnstil$ is smooth.  We will thus
use the Berline-Vergne localization formula \eqref{bvform} instead of
Rossmann's formula, and instead of \eqref{rossflag} we can work with
an explicit expression, not containing equivariant multiplicities
which need to be computed.  This allows us to write down the fixed
point formula for $\epd{\Thetab_d}$ obtained by substituting a
simplified version of\eqref{rossflag} into \eqref{epqres}, and then
compare it to the residue formula \eqref{fixedthree}. In these cases
we can describe the set $\Bipio$ easily as well. The formulas below
are justified in \S\ref{sec:howto}.

For $d=2$, we have $\O=\eqnstil\cong\P^1$. There are two fixed
points in $\eqnstil$:
\[
\Bipio={\boldsymbol{\Pi}}_2=\{([1],[2]),\  \ ([1],[1,1])\}.
\]
Then our fixed point formula reads as follows:
\begin{multline*}
\epd{\Theta_2} =\sum_{s=1}^n\sum_{t\neq s}^n \frac{1}{\prod_{i\neq
s}^n(\lambda_i- \lambda_s)\prod_{i\neq s,t}^n
(\lambda_i-\lambda_t)}\\ \times \left(
\frac{\prod_{j=1}^k(\theta_j-\lambda_s)\prod_{j=1}^k(\theta_j-\lambda_t)}{2\lambda_s-\lambda_t}
+\frac{\prod_{j=1}^k(\theta_j-\lambda_s)\prod_{j=1}^k(\theta_j-2\lambda_s)}{\lambda_t-2\lambda_s}
\right).
\end{multline*}
This is equal to the residue \eqref{fixedtwo}:

\begin{multline*}
\res_{z_1=\infty}\res_{z_2=\infty}\frac{z_1-z_2}{\prod_{i=1}^n(\lambda_i-z_1)\prod_{i=1}^n(\lambda_i-z_2)}\\
\times \left(
\frac{\prod_{j=1}^k(\theta_j-z_1)\prod_{j=1}^k(\theta_j-z_2)}{2z_1-z_2}
+\frac{\prod_{j=1}^k(\theta_j-z_1)\prod_{j=1}^k(\theta_j-2z_1)}{z_2-2z_1}
\right).
\end{multline*}
Proposition \ref{vanishthm} states that the residue of the second
term vanishes; this is easy to check by hand.

For $d=3$, the orbit closure $\O$ is a smooth 3-dimensional
hypersurface in $\eqnstil$. There are $6$ fixed points in $\O$,
namely
\begin{multline*}
\Bipio=\{([1],[2],[3]),\, ([1],[2],[1,2]),\, ([1],[2],[1,1]), \\
([1],[1,1],[3]),\, ([1],[1,1],[1,1,1]), \, ([1],[1,1],[2])\};
\end{multline*}
the remaining $2$ fixed points in $\eqnstil$ do not belong to $\O$
(see Proposition \ref{fixedpointsO}):
\[
([1],[2],[1,1,1]),\, ([1],[1,1],[1,2])\notin\Bipio.
\]
Hence the corresponding fixed point formula has $6$ terms:
\begin{multline*}
\epd{\Theta_3}=\sum_{s=1}^n\sum_{t \neq s}^n\sum_{u \neq s,t}^n
\frac{\prod_{j=1}^k(\theta_j-\lambda_s)}{\prod_{i\neq
s}^n(\lambda_i- \lambda_s)\prod_{i\neq s,t}^n
(\lambda_i-\lambda_t)\prod_{i\neq s,t,u}^n
(\lambda_i-\lambda_u)}\cdot \left[
\frac{\prod_{j=1}^k(\theta_j-\lambda_t)}{2\lambda_s-\lambda_t}\cdot \right.\\
\left(
\frac{\prod_{j=1}^k(\theta_j-\lambda_u)}{(2\lambda_s-\lambda_u)(\lambda_s+\lambda_t-\lambda_u)}
+\frac{\prod_{j=1}^k(\theta_j-\lambda_s-\lambda_t)}{(\lambda_u-\lambda_s-\lambda_t)(2\lambda_s-\lambda_s-\lambda_t)}
+\frac{\prod_{j=1}^k(\theta_j-2\lambda_s)}{(\lambda_u-2\lambda_s)(\lambda_s+\lambda_t-2\lambda_s)}
\right) +\\
\left.\frac{\prod_{j=1}^k(\theta_j-2\lambda_s)}{\lambda_t-2\lambda_s}\cdot
\left(
\frac{\prod_{j=1}^k(\theta_j-\lambda_u)}{(\lambda_t-\lambda_u)(3\lambda_s-\lambda_u)}
+\frac{\prod_{j=1}^k(\theta_j-3\lambda_s)}{(\lambda_u-3\lambda_s)(\lambda_t-3\lambda_s)}
+\frac{\prod_{j=1}^k(\theta_j-\lambda_t)}{(\lambda_u-\lambda_t)(3\lambda_s-\lambda_t)}
\right)\right].
\end{multline*}
The corresponding residue formula \eqref{fixedtwo} also has $6$
terms:
\begin{multline*}
\epd{\Theta_3}=\res_{z_1=\infty}\res_{z_2=\infty}\res_{z_3=\infty}
\frac{(z_1-z_2)(z_1-z_3)(z_2-z_3)\prod_{j=1}^k(\theta_j-z_1)}
{\prod_{i=1}^n(\lambda_i- z_1)\prod_{i=1}^n
(\lambda_i-z_2)\prod_{i=1}^n (\lambda_i-z_3)}\times \left[
\frac{\prod_{j=1}^k(\theta_j-z_2)}{2z_1-z_2}\cdot \right. \\
\left( \frac{\prod_{j=1}^k(\theta_j-z_3)}{(2z_1-z_3)(z_1+z_2-z_3)}
+\frac{\prod_{j=1}^k(\theta_j-z_1-z_2)}{(z_3-z_1-z_2)(2z_1-z_1-z_2)}
+\frac{\prod_{j=1}^k(\theta_j-2z_1)}{(z_3-2z_1)(z_1+z_2-2z_1)}
\right) +\\ \left.
\frac{\prod_{j=1}^k(\theta_j-2z_1)}{z_2-2z_1}\cdot \left(
\frac{\prod_{j=1}^k(\theta_j-z_3)}{(z_2-z_3)(3z_1-z_3)}
+\frac{\prod_{j=1}^k(\theta_j-3z_1)}{(z_3-3z_1)(z_2-3z_1)}
+\frac{\prod_{j=1}^k(\theta_j-z_2)}{(z_3-z_2)(3z_1-z_2)}
\right)\right].
\end{multline*}
Here, again, the last $5$ terms vanish, and only the one corresponding
to the distinguished fixed point $([1],[2],[3])$ remains, leaving us with
\eqref{fixedthree}.

For $d>3$, the variety $\O_d\subset\eqnstil_d$ is singular. This
means that the analogs of these formulas involve calculation of
equivariant multiplicities, which is a rather difficult problem. We
present some of these computations in \S~\ref{sec:howto}.

\subsection{The vanishing of residues}
\label{sec:vanres}

In this paragraph, we describe the conditions under which iterated
residues  of the type appearing in the sum in \eqref{fixedtwo}
vanish.

We start with the 1-dimensional case, where the residue at infinity
is defined by \eqref{defresinf} with $d=1$. By bounding the integral
representation along a contour $|z|=R$ with $R$ large, one can
easily prove
\begin{lemma}\label{1lemma}
  Let $p(z),q(z)$ be polynomials of one variable. Then
\[\res_{z=\infty}\frac{p(z)\,dz}{q(z)}=0\quad\text{if }\deg(p(z))+1<\deg(q).
\]
\end{lemma}

Consider now the multidimensional situation. Let $p(\bz),q(\bz)$ be
polynomials in the $d$ variables $z_1\ddd z_d$, and assume that
$q(\bz)$ is the product of linear factors $q=\prod_{i=1}^N L_i$, as
in \eqref{fixedthree}. We continue to use the notation $\dbz=dz_1\dots dz_d$.
We would like to formulate conditions under which the iterated
residue
\begin{equation}
  \label{ires}
\ires\frac{p(\bz)\,\dbz}{q(\bz)}
\end{equation}
vanishes. Introduce the following notation:
\begin{itemize}\label{notations}
\item For a set of indices $S\subset\{1\ddd d\}$, denote by $\deg(p(\bz);S)$
  the degree of the one-variable polynomial $p_S(t)$ obtained from $p$
  via the substitution $z_m\to
  \begin{cases}
t\text{ if }m\in S,\\ 1\text{ if }m\notin S.
  \end{cases}$
\item For a nonzero linear function $L=a_0+a_1z_1+\ldots +a_dz_d$,
  denote by $\coeff(L,z_l)$ the coefficient $a_l$;
\item finally, for $1\leq m\leq d$, set
  \[\lead(q(\bz);m)=\#\{i;\;\max\{l;\;\coeff(L_i,z_l)\neq0\}=m\},\]
which is the number of those factors $L_i$ in which the coefficient
of $z_m$ does not vanish, but the coefficients of $z_{m+1},\ldots,
z_d$ are $0$.
\end{itemize}
Thus we group the $N$ linear factors of $q(\bz)$ according to the
nonvanishing coefficient with the largest index; in particular, for
$1\leq m\leq d$ we have
\[   \deg(q(\bz);m)\geq\lead(q(\bz);m),\, \text{ and } \sum_{m=1}^d\lead(q(\bz);m)=N.
\]
Now applying Lemma \ref{1lemma} to the first residue in
\eqref{ires}, we see that
\[ \res_{z_d=\infty}\frac{p(z_1,\ddd,z_{d-1},z_d)\dbz}{q(z_1,\ddd,z_{d-1},z_d)}=0
\]
whenever $\deg(p(\bz);d)+1<\deg(q(\bz),d)$; in this case, of course,
the entire iterated residue \eqref{ires} vanishes.

Now we suppose the residue with respect to $z_d$ does not vanish,
and we look for conditions of vanishing of the next residue:
\begin{equation}
  \label{2res}
\res_{z_{d-1}=\infty}\res_{z_d=\infty}\frac{p(z_1,\ddd,z_{d-2},z_{d-1},z_d)\dbz}
{q(z_1,\ddd,z_{d-2},z_{d-1},z_d)}.
\end{equation}
 Now the condition $\deg(p(\bz);d-1)+1<\deg(q(\bz),d-1)$ will
{\em insufficient};  for example,
\begin{equation}
\res_{z_{d-1}=\infty}\res_{z_d=\infty}\frac{dz_{d-1}dz_d}{z_{d-1}(z_{d-1}+z_d)}=
\res_{z_{d-1}=\infty}\res_{z_d=\infty}\frac{dz_{d-1}dz_d}{z_{d-1}z_d}
\left(1-\frac{z_{d-1}}{z_d}+\ldots\right)=1.
\end{equation}
After performing the expansions \eqref{omegaexp} to $1/q(\bz)$, we
obtain a Laurent series with terms $z_1^{-i_1}\ldots z_d^{-i_d}$ such that
$i_{d-1}+i_d\geq\mathrm{deg}(q(z);d-1,d)$, hence the condition
\begin{equation}\label{toprove}
\deg(p(\bz);d-1,d)+2<\deg(q(\bz);d-1,d)
\end{equation}
will suffice for the vanishing of \eqref{2res}.

There is another way to ensure the vanishing of \eqref{2res}: suppose
that for $i=1\ddd N$, every time we have $\coeff(L_i,z_{d-1})\neq0$,
we also have $\coeff(L_i,z_{d})=0$, which is equivalent to the
condition $\deg(q(\bz),d-1)=\lead(q(\bz);d-1)$.  Now the Laurent
series expansion of $1/q(\bz)$ will have terms $z_1^{-i_1}\ldots
z_d^{-i_d}$ satisfying $i_{d-1}\geq
\deg(q(\bz),d-1)=\lead(q(\bz);d-1)$, hence, in this case the vanishing
of \eqref{2res} is guaranteed by $
\deg(p(\bz),d-1)+1<\deg(q(\bz),d-1)$.  This argument easily
generalizes to the following statement.
\begin{prop}
  \label{vanishprop}
Let $p(\bz)$ and $q(\bz)$ be polynomials in the variables $z_1\ddd
z_d$, and assume that $q(\bz)$ is a product of linear factors:
$q(\bz)=\prod_{i=1}^NL_i$; set $\dbz=dz_1\dots dz_d$. Then
\[ \ires\frac{p(\bz)\dbz}{q(\bz)} = 0
\]
if for some $l\leq d$, either of the following two options hold:
\begin{itemize}
\item $\deg(p(\bz);d,d-1,\dots,l)+d-l+1<\deg(q(\bz);d,d-1,\dots,l),$
\\ or
\item  $\deg(p(\bz);l)+1<\deg(q(\bz);l)=\lead(q(\bz);l)$.
\end{itemize}
\end{prop}
Note that for the second option, the equality
$\deg(q(\bz);l)=\lead(q(\bz);l)$ means that
  \begin{equation}
    \label{op2cond}
\text{ for each }i=1\ddd N\text{ and }m>l,\,
\coeff(L_i,z_{l})\neq0\text{ implies }\coeff(L_i,z_{m})=0.
  \end{equation}

Recall that our goal is to show that all the terms of the sum in
\eqref{fixedtwo} vanish except for the one corresponding to
$\bipi_\dist=([1]\ddd [d])$. Let us apply our new-found tool,
Proposition \ref{vanishprop}, to the terms of this sum, and see what
happens.

Fix a sequence $\bipi=(\pi_1,\dots,\pi_d)\in\Bipi$, and consider the
iterated residue corresponding to it on the right hand side of
\eqref{fixedtwo}. The expression under the residue is the product of
two fractions:
\[\frac{p(\bz)}{q(\bz)}=\frac{p_1(\bz)}{q_1(\bz)}\cdot\frac{p_2(\bz)}{q_2(\bz)},\]
where
\begin{equation}
  \label{pq}
\frac{p_1(\bz)}{q_1(\bz)}=\frac{\displaystyle
Q_\bipi(\bz)\,\prod_{m<l}(z_m-z_l) }{\displaystyle
\prod_{l=1}^d\prod_{\um\tau\leq
l}^{\tau\neq\pi_1\ddd\pi_l}(z_\tau-z_{\pi_{l}})}\text{ and
  }
\frac{p_2(\bz)}{q_2(\bz)} = \frac{\displaystyle
\prod_{m=1}^d\prod_{j=1}^k(\theta_j-z_{\pi_m})} {\displaystyle
\prod_{l=1}^d\prod_{i=1}^n(\lambda_i-z_l)}.
\end{equation}

Note that $p(\bz)$ is a polynomial, while $q(\bz)$ is a product of
linear forms, and that $p_1(\bz)$ and $q_1(\bz)$ are independent of
$n$ and $k$, and depend on $d$ only.

As a warm-up, we show that if the last element of the sequence is
not the trivial partition, i.e. if $\pi_d\neq[d]$, then already the
first residue in the corresponding term on the right hand side of
\eqref{fixedtwo} -- the one with respect to $z_d$ -- vanishes.
Indeed, if $\pi_d\neq[d]$, then $\deg(q_2(\bz);d)\geq n$, while
$z_d$ does not appear in $p_2(\bz)$. Then, assuming that $d\ll n$,
we have $\deg(p(\bz);d)\ll\deg(q(\bz);d)$, and this, in turn,
implies the vanishing of the residue with respect to $z_d$ (see
Proposition \ref{vanishprop}).

We can thus assume that $\pi_d=[d]$, and proceed to the study of the
next residue, the one taken with respect to $z_{d-1}$. Again, assume
that $\pi_{d-1}\neq[d-1]$. As in the case of $z_d$ above, $d\ll n$
implies $\deg(p(\bz);d-1)\ll\deg(q(\bz);d-1)$. However, now we
cannot use the first option in Proposition \ref{vanishprop}, because
$\deg(p_2(\bz);d-1,d )=k\geq n$. In order to apply the second
option, we have to exclude all linear factors from $q_1(\bz)$ which
have nonzero coefficients in front of both $z_{d-1}$ and $z_d$. The
fact that $\pi_d=[d]$, and the restrictions $\um{\pi_l}\leq l$,
$l=1\ddd d$, tell us that there are two troublesome factors:
$(z_d-z_{d-1})$ and $(z_d-z_{d-1}-z_1)$ which come from the two
partitions: $\tau=[d-1]$ and $\tau=[d-1,1]$ in the $l=d$ part of
$q_1(\bz)$. The first of the two fortunately cancels with a factor
in the Vandermonde determinant in the numerator; as for the second
factor: our only hope is to find it as a factor in the polynomial
$Q_\bipi$.

Continuing this argument by induction, we can reduce Proposition
\ref{vanishthm} to the following statement about the equivariant
multiplicities $Q_\bipi$, $\bipi\in\Bipi$.
\begin{prop}\label{critical}
  Let $l\geq1$, and let $\bipi$ be an admissible sequence of
  partitions of the form \eqref{bipiform}, where $\pi_l\neq[l]$. Then
  for $m>l$, and every partition $\tau$ such that
  $l\in\tau$, $\um\tau\leq m$, and $|\tau|>1$, we have
  \begin{equation}
    \label{zdividesq}
 (z_\tau-z_m)|Q_\bipi.
  \end{equation}
\end{prop}
This statement will be proved in the next paragraph:
\S\ref{sec:homog}. For now, we will assume that it is true, and give
a quick proof of the result with which we started this section.
\\{\em
Proof of Proposition \ref{vanishthm}}: Let $\bipi\neq\bipi_\ds$ be
an admissible sequence of partitions. This means that there is $l>1$
such that $\pi_l\neq [l]$, but $\pi_m=[m]$ for $m>l$:
\begin{equation}
  \label{bipiform}
\bipi = (\pi_1,\dots\pi_l,[l+1],[l+2]\ddd[d]).
\end{equation}
Note that $l$ does not appear anywhere in $\bipi$, and thus we can
conclude $\deg(p(\bz);l)\ll\deg(q(\bz);l)$ from $d\ll n$, as usual.
This allows us to apply the second option of Proposition
\ref{vanishprop} to the residue taken with respect to $z_l$ as long
as we can cancel from $q_2(\bz)$ all factors which do not satisfy
condition \eqref{op2cond}.

These factors are of the form $z_\tau-z_m$, where $m>l$ and
$l\in\tau$. If $|\tau|=1$, i.e. if $\tau=[l]$, then we can find this
factor in the Vandermonde determinant in the numerator.  We can use
Proposition \ref{critical} to cancel the rest of the factors, as
long as we make sure that such factors occur in $q_1(\bz)$ with
multiplicity 1. This is straightforward in our case, since the
variable $z_m$ with $m\geq l$ may appear only in the $m$th factor of
$q_1(\bz)$. \qed

\subsection{The homogeneous ring of $\eqnstil$ and factorization of
  $Q_\bipi$}
\label{sec:homog}

Now we turn to the proof of Proposition \ref{critical}. Let
$\bipi\in\Bipi$ be an admissible sequence of partitions. Recall (cf.
\eqref{qbipidef}) that $Q_\bipi$ is the $T_d$-equivariant Poincar\'e
dual of the part $\O_\bipi=\pre^{-1}(\O)\cap\NN_\bipi$ of the orbit
closure $\O$ in the linear chart $\NN_\bipi$ (cf.
\eqref{compemult}); this latter linear space is endowed with
coordinates $\uh_{\tau|\bipi}^l$ defined in \eqref{npicoordinates}.

Our plan is to use Lemma \ref{divisible}, which, when applied to our
situation, says that the divisibility relation \eqref{zdividesq}
follows if we find a relation in the ideal of the subvariety
$\O_\bipi\subset\NN_\bipi$ expressing the appropriate variable
$\uh_{\tau|\bipi}^m$ as a polynomial of the rest of the variables.

We will lift the calculation from $\eqnstil$ to the vector space
$\homnewd$. Denote by $\polye$ the ring of polynomial functions on
$\homnewd$, i.e. the space of polynomials in the variables
$u^l_\tau$, $1\leq l\leq d$, $\um\tau\leq l$. As one can see from
Definition \ref{defnpi}, and \eqref{npicoordinates}, the relations on the two
spaces are connected as follows:

\begin{lemma}
  \label{prescr}
  Let $Z\in \polye$ be a polynomial on $\hofi$, and let
  $M\subset\hofi$ be a closed subvariety, such that $Z|M$ vanishes.
  Then the restricted polynomial $\hat{Z}=Z|\NN_\bipi$, written in terms of the
  coordinates $\uh_{\cdot|\bipi}$, may be obtained from $Z$ as follows:
\begin{itemize}
\item setting $u_{\pi_l}^{l}$ to $1$, for $l=1\ddd d$,
\item setting $u_{\pi_l}^{m}$ to $0$, for $1\leq l\leq m\leq d$,
\item replacing the remaining variables $u_{\tau}^{l}$ by $\uh_{\tau|\bipi}^l$.
\end{itemize}
In addition, $\hat{Z}$ vanishes on $M\cap\NN_\bipi$.
\end{lemma}

Eventually, using  this lemma with $M=\overline{B_L\epsref B_R}$
and  $M\cap\NN_\bipi=\O_\bipi$, we will be able to produce the
necessary relations in the defining ideal of $\O_\bipi\subset\NN_\bipi$.
As most of the action will take space in $\polye$, our next task is to
set up some convenient notation for this ring.

The ring $\polye$ carries a right action of the group $B_L$, and a
left action of the group $B_R$. In particular, it has two
multigradings induced from the $T_L$ and $T_R$ actions: the
$L$-multigrading is the vector of multiplicities
$(\mult(i,\pi),\,i=1\ddd d)$, while the $R$-multigrading is the $l$th
basis vector in $\Z^d$. A combination of these gradings will be
particularly important for us (cf. Definition \ref{admissible}):
\begin{equation}
  \label{defdefect}
\defect(u_{\pi}^{l})=l-\um\pi;
\end{equation}
this induces a $\Z^{\geq0}$-grading on $\polye$.

Recall that the projection $B_d\to T_d$ is a group homomorphism,
whose kernel is the subgroup of unipotent matrices. We denote the
corresponding nilpotent Lie algebras of strictly upper-triangular
matrices by $\n_R$ and $\n_L$ for $B_R$ and $B_L$, respectively..

The two Lie algebras, $\n_L$ and $\n_R$ are generated by the simple
root vectors
\[\Delta_L=\{E^L_{l,l+1};\;l=1\ddd d-1\},\;\text{and}\;
\Delta_R=\{E^R_{l,l+1};\;l=1\ddd d-1\},
\]
respectively, where $E_{l,l+1}$ is the matrix whose only nonvanishing
entry is a 1 in the $l$th row and $l+1$st column. Let us write down
the action of these root vectors on $\polye$ in the coordinates
$u_{\tau}^{l}$, $|\tau|\leq l\leq d$. We first define certain
operations on partitions:
\begin{itemize}
\item given a positive integer $m$ and a partition $\tau\in\Pi$, denote
  by $\tau\cup m$ the partition with $m$ added to $\tau$,
  e.g. $[2,3,4]\cup 3=[2,3,3,4]$
\item if $m\in\tau$, then denote by $\tau-m$ the partition $\tau$ with
  one of the $m$s deleted, e.g. $[2,4,4,5,5,5,6]-5=[2,4,4,5,5,6]$;
\item more generally, we will write $[2,4,5,5]\cup[3,4]=[2,3,4,4,5,5]$,
  and $[2,4,5,5]-[4,5]=[2,5]$.
\end{itemize}

Returning to the Lie algebra actions, we have
\begin{equation}
  \label{euact}
  \begin{cases}
 \n_Ru_{\tau}^{l}=u_{\tau}^{l}\n_L=0, \text{ if } \um\tau=l,\\
  E^R_{m,m+1}u_{\tau}^{l}=\delta_{l,m+1}u_{\tau}^{l-1},\;
u_{\tau}^{l}E^L_{m,m+1}=\mult(m,\tau)\,u_{\tau-m\,\cup\, m+1}^{l},
\text{ if } \um\tau<l.
  \end{cases}
\end{equation}
where $\delta_{a,b}$ is the Kronecker delta.  Observe that both $\n_R$ and
$\n_L$ act compatibly with the $T_R\times T_L$-multigrading, and they
both decrease the defect \eqref{defdefect}.

The following subspace will play a key role in our calculations:
\begin{equation}
  \label{defio}
  \io=\left\{ Z\in\polye;\; \n_R Z=0 \text{ and } [Z\n_L^N](\epsilon_{\rf})=0 \text{ for }
N=0,1,2,\dots\right\},
\end{equation}
where $\n_L^N$ is the subset $\{X_1\cdot\dots\cdot X_N;\;
X_i\in\n_L,\,i=1\ddd N\}$ of the universal enveloping algebra of
$\n_L$.

\begin{prop}
  \label{io}
  If $Z\in\io$, then $Z(\epsilon)=0$ for every $\epsilon\in B_L\epsref
  B_R$.
\end{prop}
\begin{proof}
First, observe that the actions of $\n_R$ and $\n_L$ described in
\eqref{euact} are compatible with the multigrading induced by the
$T_R\times T_L$-action, and hence, if $Z$ is in $\io$, then so are all
of its $T_R\times T_L$-homogeneous components. This means that without
loss of generality we may assume that $Z$ is a homogeneous element of
$\io$.

For such $Z$, clearly, $Z(\epsilon)=0\Leftrightarrow t_R
Zt_L(\epsilon)=0$ for any $t_L\in T_L$, $t_R\in T_R$.  Combining this
with the condition $\n_RZ=0$ we can conclude that the zero set of $Z$
is $B_R$-invariant, hence it is sufficient to show $Z(\epsilon)=0$ for
$B_L\epsref$.  Now, since $\ker(B_L\to T_L)=\exp(\n_L)$, the
definition of $\io$ also implies $Z(b\epsilon_\rf)=0$ for all $b\in
B_L$., and this completes the proof.
\end{proof}

\begin{rem}
 \label{geomean}
Before we proceed, we make a comment on the geometric meaning of
$\io$. The space $\{ Z\in\polye;\; \n_R Z=0 \}$ is the homogeneous
coordinate ring of $\eqnstil$, corresponding to the line bundles
induced by the characters of $T_R$. Then Proposition \ref{io} may be
interpreted as saying that $\io$ is contained in the ideal of
functions vanishing on $\O$. In fact, is not difficult to show that
$\io$ is exactly this ideal.
\end{rem}

We will be looking for polynomials $Z\in\io$ in a particular subspace
of $\polye$. To describe this space, introduce for each
$\bipi\in\Bipi$ the monomial
\begin{equation}
  \label{up}
\bu^\bipi = \prod_{l=1}^d u_{\pi_l}^{l};\text{ these satisfy }
\bu_\bipi(\epsilon_{\bipi'})=
\begin{cases}
  1,\,\text{if }\bipi=\bipi'\\
0,\,\text{otherwise.}
\end{cases}
\end{equation}
Now consider the linear span of these monomials:
\begin{equation}
  \label{deflambda}
\Lambda=\left\{\sum_{\bipi\in\Bipi}\alpha_\bipi\bu^\bipi\in\polye;\;
\alpha_\bipi\in\C\right\}.  
\end{equation}

In order to write down our formulas for certain elements of
$\Lambda\cap\io$, we need to introduce two operations on $\Bipi$.  For
a sequence of partitions $\bipi=(\pi_1\ddd\pi_d)$ and a permutation
$\sigma\in\sg d$ define the the permuted sequence
\[
\bipi\cdot\sigma=(\pi_{\sigma(1)},\dots,\pi_{\sigma(d)});
\]
this defines a natural right action of $\sg d$ on $\Pi^{\times
d}$. Note that permuting an admissible sequence $\bipi\in\Bipi$ does
not necessarily result in an admissible sequence.

The second operation modifies just one entry of $\bipi$: for
$\bipi\in\Bipi$ and $\tau\in\Pi$, define
\[ \bipi\cup_m\tau=(\pi_1\ddd,\pi_{m-1},\pi_m\cup\tau,\pi_{m+1}\ddd\pi_d).
\]

Now we are ready to write down our relations.
\begin{prop}
  \label{fundament}
  Let $\bipi\in\Bipi$ be an admissible sequence of partitions and let
  $\tau\in\Pi$ be any partition. Then following polynomial is an
  element of $\io$:
  \begin{equation}
    \label{relz}
 \relz\bipi\tau = \sum \sign(\sigma)\,\bu^{\bipi\cdot\sigma\cup_m\tau},
  1\leq m\leq d,\, \sigma\in\sg d,\, \bipi\cdot\sigma\cup_m\tau\in\Bipi,
      \end{equation}
\end{prop}

\begin{rem} \label{termsaa} The sum in \eqref{relz}  may be
  empty. This happens when there are no pairs $(\sigma,m)$ satisfying
  the conditions in \eqref{relz}. Note, however, that no two terms of
  this sum may cancel each other.
\end{rem}
\begin{proof}
  We begin by noting that $\relz\bipi\tau$ is of pure $T_R\times T_L$
  weight. Indeed, the torus $T_R$ acts on the whole space $\Lambda$
  with the same weight $(1,1,\dots,1)$, while the $l$th component of
  the $T_L$-weight of a term of $\relz\bipi\tau$ is equal to
  $\mult(l,\tau)+\sum_{m=1}^d\mult(l,\pi_m)$.

  Next, we show that
  \begin{equation}
    \label{erkills}
E_{l,l+1}^R\relz\bipi\tau=0,\; l=1\ddd d-1,
  \end{equation}
  which implies that $\n_R\relz\bipi\tau=0$. Let us fix $l$; the
  terms of $\relz\bipi\tau$ in \eqref{relz} are indexed by pairs
  $(\sigma,m)$, and we can ignore those pairs for which
  $\um{\pi_{l+1}}+\delta_{m,l+1}\um\tau\ge l+1$, since in this case
  $E_{l,l+1}^R\bu^{\bipi\cdot\sigma\cup_m\tau}=0$. Then the vanishing
  \eqref{erkills} clearly follows if, on the set of the remaining
  pairs contributing to \eqref{relz}, we find an involution
  $(\sigma,m)\mapsto(\sigma',m')$ such that
\[
E_{l,l+1}^R\bu^{\bipi\cdot\sigma\cup_m\tau}=E_{l,l+1}^R\bu^{\bipi\cdot\sigma'\cup_{m'}\tau}
\;\text{and}\; \sign(\sigma')=-\sign(\sigma).
\]
Indeed, it is easy to check that this holds for the involution
  \[(\sigma',m')=(\sigma\cdot\trp l{l+1}, \trp l{l+1}(m)),\]
where $\trp l{l+1}\in\sg d$ is the transposition of $l$ and
$l+1$. This proves \eqref{erkills}.

Our second task is to show that $\relz\bipi\tau$ is in the linear
space
\[ \lrel=
\left\{Z\in\polye;\;\left[Z\n_L^N\right](\epsilon_{\rf})=0\text{ for
  }N=0,1,\dots\right\}.
\]
Using the Leibniz rule, it is easy to see see that
$\lrel\subset\polye$ is an ideal.

First we show that for partitions $\rho,\tau\in\Pi$ and
$m\geq\um\rho+\um\tau$ the polynomial
\begin{equation}\label{defzrel}
  Z_{\rho\tau}^m = u_{\rho\cup\tau}^{m}-\sum u_{\rho}^{t}u_{\tau}^{r},\;
t+r=m,\, t\geq\um\rho,\,r\geq\um\tau
\end{equation}
is in $\lrel$. Indeed, a quick computation produces the equality
\[ Z_{\rho\tau}^m E^L_{l,l+1}=\mult(l,\rho)Z_{\rho'\tau}^m+
\mult(l,\tau)Z_{\rho\tau'}^m,\text{ where }\rho'=\rho- l\cup
[l+1],\, \tau'=\tau- l\cup [l+1].
\]
This equality implies that it is sufficient for us to prove
$Z_{\rho\tau}^m(\epsilon_{\rf})=0$ for the case $m=\um\rho+\um\tau$.
In this case we have
\begin{equation}
Z_{\rho\tau}^m =
u_{\rho\cup\tau}^{m}-u_\rho^{\um\rho}u_\tau^{\um\tau},
\end{equation}
and this polynomial clearly vanishes on $\epsilon_{\rf}$, because
all three coordinates appearing in this relation are equal to 1
according to \eqref{defpref}.

Now we return to the proof of $\relz\bipi\tau\in\lrel$. Using the
fact that $Z_{\rho\tau}^m$ is in the ideal $\lrel$, modulo the
$\lrel$, we can replace all the factors of the form
$u_{\pi_{\sigma(m)}\cup\tau}^m$ in all the terms of $\relz\bipi\tau$
by the appropriate sum of quadratic terms in \eqref{defzrel}. Our
claim is that the resulting polynomial is identically zero, which
implies that $\relz\bipi\tau\in\lrel$.

Indeed, let us perform this substitution; the terms of the resulting
sum are parametrized  by a triple $(\sigma,m,r)$, which is obtained
by applying \eqref{defzrel} to the term of $\relz\bipi\tau$ indexed
by $(\sigma,m)$ and taking the term corresponding to $r$ in
\eqref{defzrel}. The correspondence is thus
\begin{equation}
\label{oneterm} (\sigma,m,r)\longrightarrow
u_{\pi_{\sigma(1)}}^{1}\ldots
  u_{\pi_{\sigma(m-1)}}^{m-1}u_{\pi_{\sigma(m)}}^{m-r}
u_{\tau}^{r}u_{\pi_{\sigma(m+1)}}^{m+1}\ldots
u_{\pi_{\sigma(d)}}^{d}.
\end{equation}
Just as above, we can see that the involution $(\sigma,m,r)\mapsto
(\sigma\cdot\trp{m}{m-r},m,r)$ provides us with a complete pairing
of the terms of the sum described above; each pair consists of
identical monomials with opposite signs. This implies that indeed,
the result is zero, hence $\relz\bipi\tau$ vanishes modulo $\lrel$,
i.e. $\relz\bipi\tau\in\lrel$.
\end{proof}

Armed with these relations, we are ready to {\em prove Proposition}
\ref{critical}. Recall that according to the strategy described at
the beginning of this paragraph, given $\bipi\in\Bipi$, $m$ and
$\tau$ as in Proposition \ref{critical}, we need to find a relation
of the form $\relz\cdot\cdot$, which, when restricted to $\NN_\bipi$,
expresses the variable $\uh^m_{\tau|\bipi}$ in terms of the rest of
the variables.

Thus the first thing is to study the conditions under which
$\uh^m_{\tau|\bipi}$ appears as the restriction of a monomial of the form
$\bu^{\bipi'}$. The following statement immediately follows form the
prescription Lemma \ref{prescr}.
\begin{lemma}
\label{urestr}
Given $\bipi=(\pi_1\ddd\pi_d)\in\Bipi$, a positive integer $m\leq d$,
and a partition $\tau\in\Pi\setminus\{\pi_1\ddd\pi_d\}$ satisfying
$\um\tau\leq m$, we have $\bu^{\bipi'}|\NN_\bipi=\uh^m_{\tau|\bipi}$
for some $\bipi'\in\Bipi$ if and only if
\[   \bipi'= (\pi_1\ddd\pi_{m-1},\tau,\pi_{m+1}\ddd\pi_d).
\]
\end{lemma}

Now let us take a closer look at the conditions of
Proposition \ref{critical}. We are  given $1\leq l<m\leq d$ and
$\tau\in\Pi$ satisfying
\[  \um\tau\leq m,\, l\in\tau\text{ and }|\tau|>1,
\]
and a sequence $\bipi$ of the form \eqref{bipiform} with
$\pi_l\neq[l]$. In view of Lemma \ref{urestr}, the variable
$\uh^m_{\tau|\bipi}$ will appear as the restriction to $\NN_\bipi$ of
the term $\bu^{\birho\cup_m\tau\setminus[l]}$ of a relation
$\relz\birho{\tau\setminus[l]}$ as long as
\[ \birho=
(\pi_1\ddd,\pi_l,[l+1],[l+2]\ddd,[m-1],[l],[m+1],\ddd,[d-1],[d])
\]
is admissible, which is obvious. We leave it to the reader to check is
that the rest of the terms of $\relz\birho{\tau\setminus[l]}$ cannot
contain $\uh^m_{\tau|\bipi}$ as a factor. This completes the proof of
Proposition \ref{critical} and thus also the proof of Proposition
\ref{vanishthm}. \qed

This proof suggests a simple criterion for finding out for which
$\bipi\in\Bipi$ the monomial $\bu^\bipi$ appears in one of the
relations \eqref{relz}.

\begin{defi} \label{defcomplete} We will call an admissible sequence
  of partitions $\bipi=(\pi_1\ddd\pi_d)$ {\em complete} if for every
  $l\in\dod$ and every nontrivial subpartition $\tau\subset\pi_l$,
  there is $m\in\dod$ such that $\pi_m=\tau$.
\end{defi}
Taking into account Remark \ref{termsaa}, we have the following
criterion.

\begin{lemma} \label{occurs}
  A monomial $\bu^{\bipi}$ appears in a relation $\relz\birho\tau$
  for some $\birho\in\Bipi$ and $\tau\in\Pi$ if and only if $\bipi$ is
  {\em not} complete.
\end{lemma}

\subsection{The fixed points of the $T_L$-action on $\O$}
\label{sec:fixedpoints}

As a small detour, based on the results of the previous paragraph,
we obtain a rather powerful criterion for $\bipi\in\Bipi$ {\em not}
to belong to $\Bipio$, i.e. we will construct a large number of
$T_L$-fixed points which do not lie in $\O$. We note, however, that
describing the set $\Bipio$ remains an interesting open problem. Our
starting point is \eqref{up}.
\begin{lemma}
If the monomial $\bu^{\bipi}$ appears with nonzero coefficient in a
polynomial from $\Lambda\cap\io$, then the fixed point
$\teps_\bipi\notin\O$, i.e. $\bipi\notin\Bipio$.
\end{lemma}
\begin{proof}
  Indeed, let $Z$ be such a polynomial. According to Proposition
  \ref{io}, a polynomial in $\io$ vanishes at all points of
  $\O$. On the other hand, it is clear from \eqref{up} that all but
  exactly one of the terms of $Z$ vanishes at $\epsilon_{\bipi}$, and
  hence $Z(\epsilon_{\bipi})\neq0$.
  \end{proof}

  Combining this statement with Lemma \ref{occurs} we have the
  following.
\begin{prop}\label{fixedpointsO}
If $\bipi\in\Bipio$. i.e. if $\teps_\bipi\in\O$, then the sequence
$\bipi$ is complete.
\end{prop}

This Proposition provides us a rather strict necessary , although,
as an example below shows, not sufficient condition for $\bipi$ to
be in $\Bipio$.
\begin{exa}
  \begin{enumerate}
  \item  The  sequence
    \[ ([1],[2],\dots,[d-1],[l,m]),\;\text{ where }l+m\leq d.\] is
    complete, and, in fact, it corresponds to a fixed point.
\item For $d=3,4$, the reverse of Proposition \ref{fixedpointsO}
  holds: if $\bipi$ is complete then the fixed point $\teps_{\bipi}$ lies in
  the orbit closure $\O_d$, see section \S\ref{sec:howto}.
\item The completeness of $\bipi$ is a necessary but not sufficient
  condition for $\bipi$ to be in $\Bipio$. An example is the
  following zero-defect sequence of partitions: let $d=60$,
  $\tau=[1,12,12,15,20]$ and set
\[
\pi_l =
\begin{cases}
  \rho,\text{ if } \rho\subset\tau\text{ and }\um\rho=l,\\
[l],\text{ otherwise.}
\end{cases}
\]
By definition, this is a complete sequence of partitions, but it is
not in $\O$, which is left as an exercise.
  \end{enumerate}
\end{exa}

\subsection{The distinguished fixed point and the main result}
\label{sec:mainres}

Now we turn our attention to our much simplified formula
\eqref{fixedthree} for the Thom polynomial of the $A_d$-singularity.

The proof of the vanishing of the contributions to \eqref{fixedtwo},
naturally, fails at the fixed point $\teps_\ds$. Indeed, for the
for the factors \eqref{pq} in the case of  the
distinguished sequence $\bipi_\ds$, we have
$\deg(p_2(\bz);l)>\deg(q_2(\bz);l)$ for $l=1\ddd d$, and hence we
cannot apply Proposition \ref{vanishprop}.

The factorization arguments of \S\ref{sec:homog} may be partially
saved, however.
Indeed, for the case of the distinguished partition $\bipi_\ds$,
each $T_L$-weight $z_\tau-z_l$ of $\NN_\ds$ appears with
multiplicity one (cf. end of \S\ref{sec:vanres}). Hence, again, we
can apply Lemmas \ref{divisible}, \ref{urestr} and \ref{occurs} to
conclude that for $|\tau|>1$,
\[  (z_\tau-z_l)\,|\,Q_\ds\;\text{ if }\;
 ([1],[2]\ddd[l-1],\tau,[l+1]\ddd[d-1],[d])\text{ is not complete}.
\]
Clearly, such a sequence is complete if and only if $|\tau|=2$, and
this means that in the fraction on the right hand side of
\eqref{fixedfour}, we can cancel all factors between the numerator
and the denominator corresponding to partitions $\tau$ with
$|\tau|>2$. This reduces the denominator to the product of the
factors with $|\tau|=2$:
\[
\prod (z_m+z_r-z_l), \;1\leq m\leq r,\,m+r\leq l\leq d,
\]
while $Q_\ds$ is replaced by a polynomial $\qd d$, whose degree is
much smaller than that of $Q_\ds$. Note that in this case no factors
of the Vandermonde in the numerator are canceled; the fraction in
\eqref{fixedfour} thus simplifies to
\[  \frac{(-1)^d\prod_{m<l}(z_m-z_l)\,\QQ_d(z_1\ddd z_d)}
{\prod_{l=1}^d\prod_{m=1}^{l-1}\prod_{r=1}^{\min(m,l-m)}(z_m+z_r-z_l)}
\]
The polynomial $\qd d$, just as $Q_\ds$, only depends on $d$; we mark
its $d$-dependence explicitly.

All that remains to do before we can formulate our final result,
is to  describe the geometric meaning of this cancellation, and that of
the polynomial $\qd d$ itself.

First, note that $\bipi_\ds$ is of the defect-0 type, hence, according
to Proposition \ref{cover} (3) and Lemma \ref{qpi}, we have an action
of the upper-triangular group  $B_\ds$ on $\NN_\ds$ given by
\eqref{action}; moreover, $\epsref\in\NN_\ds$
and $\O_\ds=\overline{B_\ds\cdot\epsref}$. Remarkably, this action is
also linear (cf. Remark \ref{quad}), because the $B_L\times
B_R$-action on $\hofi$ preserves the length of the partitions, and
$\bipi_\ds$ contains all the partitions of length 1.

Next, define the linear subspace $\Nh_d\subset\NN_\ds$:
\begin{equation}\label{nhd}
 \Nh_d=\{\epsilon\in\NN_\ds;\,\uh_{\tau|\ds}^{m}(\epsilon)=0\text{ for
}|\tau|>2\}\subset\hom(\C^d,\sym^2\C^d),
\end{equation}
and let $\phat:\NN_\ds\to\Nh_d$ be the natural projection.  Then
 (cf. Remark \ref{geominterpr}) we can conclude that
 \begin{equation}
   \label{ohdef}
   \qd d = \epd{\OO_d,\Nh_d},\;
\text{ where }\,  \OO_d=\phat(\OO_\ds).
    \end{equation}

    In addition, it is easy to see that $\phat$ commutes with the
    $B_\ds$-action, in particular, $\Nh_d$ in $\NN_\ds$ is
    $B_\ds$-invariant. The linear representation of $B_\ds$ on $\Nh_d$
    is easily identified with an action of degree-3 tensors (see the
    Theorem below). In any case, we have
\[   \OO_d = \overline{B_d\ehat},\;\text{where }\ehat=\phat(\epsref).
\]

Stripping our formulas of extraneous notation, we can formulate our
main result in a self-contained manner as follows:
\begin{thm}
  \label{dathm}
  Let $T_d\subset B_d\subset\gl d$ be the subgroups of invertible
  diagonal and upper-triangular matrices, respectively; denote the
  diagonal weights of $T_d$ by $z_1\ddd z_d$.  Consider the $\gl
  d$-module of 3-tensors $\hom(\C^d,\sym^2\C^d)$; identifying the
  weight-$(z_m+z_r-z_l)$ symbols $q^{mr}_l$ and $q^{rm}_l$, we can
  write a basis for this space as follows:
\[ \hom(\C^d,\sym^2\C^d)=\bigoplus \C q^{mr}_l,\;  1\leq m,r,l \leq d.
\]
Consider the reference element
\[ \erf=\sum_{m=1}^d\sum_{r=1}^{d-m}q_{mr}^{m+r},
\]
in  the $B_d$-invariant subspace
\begin{equation}
  \label{nhmodule}
    \Nh_d = \bigoplus_{1\leq m+r\leq l\leq d} \C q^{mr}_l\subset
\hom(\C^d,\sym^2\C^d).
\end{equation}
Set the notation $\OO_d$ for the orbit closure
$\overline{B_d\erf}\subset\Nh_d$, and consider its $T_d$-equivariant
Poincar\'e dual
\[    \QQ_d(z_1\ddd z_d) = \epd{\OO_d,\Nh_d}_{T_d},
\]
which is a homogeneous polynomial of degree
$\dim(\Nh_d)-\dim(\OO_d).$

Then for arbitrary integers $n\leq k$, the Thom polynomial for the
$A_d$-singularity with $n$-dimensional source space and
$k$-dimensional target space is given by the following iterated
residue formula:
\begin{equation}
  \label{thethom}
   \epd{\Theta_d}=
\sires \frac{(-1)^d\prod_{m<l}(z_m-z_l)\,\QQ_d(z_1\ddd z_d)}
{\prod_{l=1}^d\prod_{m=1}^{l-1}\prod_{r=1}^{\min(m,l-m)}(z_m+z_r-z_l)}
\prod_{l=1}^d \ct\left(\frac1{z_l}\right)\,z_l^{k-n}\;dz_l,
\end{equation}
where $\ct(\cdot)$ is the generating function of the relative Chern
classes given in \eqref{deftc}.
\end{thm}
Let us briefly review our the proof of this theorem. We began by
interpreting the Thom polynomial as an equivariant Poincar\'e dual of
a variety $\Thetab_d$ in the space of map-jets (cf. \eqref{deftp} and
Proposition \ref{jets}).  Next, we constructed a birational model for
$\Theta_d$ in Proposition \ref{themodelprop}, and then we applied a
localization formula \eqref{locsing} to this model, which resulted in
expression \eqref{fixedtwo} for the Thom polynomial. Finally, by
studying certain explicit relations and under the assumption that
$d\ll n$, we uncovered a cancellation phenomenon, which lead to the
simplified formula \eqref{thethom}.

Note that the formulation of Theorem \ref{dathm} is more general than
to what we seem to be entitled: Proposition \ref{vanishthm} includes
the assumption $d\ll n$, while here we claim that our statement holds
for any $d$ and $n\leq k$. To finish the proof, we simply need to
point our that according to Proposition \ref{collect}, an expression
of a Thom polynomial in the relative Chern classes holds for large
$n$, then the same expression works for any $n$.\qed

Let us make a few final comments.  It is not difficult to see that
formula \eqref{thethom} manifestly satisfies all properties listed in
Proposition \ref{collect}. In particular, it only depends on the
codimension $k-n$, and reducing the codimension by 1 leads to shifting
the indices of the relative Chern classes down by 1.  Another benefit
of the result is that it shows that the {\em Thom series} introduced in
\cite{fr}, which, in principle has infinitely many parameters, is
governed by a finite object: $\qd d$. A detailed study of the
polynomial $\qd d$ will be given in a later publication
\cite{berczi}. In the final section of our paper, we turn to examples,
and explicit calculations.

\section{How to calculate $\QQ_d$? Explicit formulas for Thom
  polynomials}
\label{sec:howto}

Theorem \ref{dathm} reduces the computation of the Thom polynomials
of the algebra $A_d$ to that of the polynomial $\qd d$, which is the
\ePd of a $B_d$-orbit in a certain $B_d$-invariant subspace of
3-tensors in $d$ dimensions. Note that the parameters $n$ and $k$ do
not enter this picture; in particular, $\qd d$ only depends on $d$.

Clearly, in principle, the computation of $\qd d$ is a finite
problem in commutative algebra, which, for each value of $d$, can be
handled by a computer algebra package such as Macaulay. However, the
number of variables and the degree of $\qd d$ grow rather quickly:
they are of order $d^3$. More importantly, computer algebra programs
have difficulties dealing with parametrized subvarieties already in
very small examples.

At this point, we do not have an efficient method of computation for
$\qd d$ in general. The purpose of this section is to show how to
compute $\qd d$ for small degrees: $d=2,3,4,5,6$. At the end, we also
present an application of our result to the conjectured positivity of
the coefficients of the Thom polynomials in Section \S
\ref{sec:positive}.

\subsection{The degree of $\qd d$}\label{sec:codimension}
The degree of the polynomial $\qd d$ is the codimension of the orbit
$B_d\epsilon_\rf$, or that of its closure $\OO_d$, in $\Nh_d$.

Recall that $\Nh_d$ has a basis indexed by the set of indices
$\{m+r\leq l\leq d\}$. 
An elementary computation shows that $\dim \Nh_d$ is given by a
cubic quasi-polynomial in $d$ with leading term $d^3/24$.

On the other hand, we have
\[ \dim(B_d\erf) = \dim(B_d) - \dim(H_d) = {d+1 \choose 2}-d={d \choose 2}.
\]

Next, denote by $\Nh_d^0$ the {\em minimal} or defect-zero part of
$\Nh_d$ spanned by the vectors $\{q_{mr}^{l};\;m+r=l\leq d\}$, and let
$\prz:\Nh_d\to \Nh_d^0$ be the natural projection; note that $\erf\in
\Nh_d^0$. Recall that $B_d=T_dU_d$, where $U_d\subset B_d$ is the
subgroup of unipotent matrices. It is easy to check that $U_d$ acts
trivially on $\Nh_d^0$, and its action commutes with the projection
$\prz$. Now introduce the toric orbit
$T_d\erf\subset \Nh_d^0$ and its closure $\KK\subset \Nh_d^0$. The
following is a simple consequence of the preceding arguments.
\begin{lemma} \label{bovert}
The projection $\prz$ restricted to the orbit $B_d\erf$ establishes
a fibration over the toric orbit $T_d\erf$. This map extends to a
map between the closures  $\OO\to\KK$, where
$\KK=\overline{T_d\erf}$.
\end{lemma}
\begin{rem}
  We note that there are standard algorithms to compute the \ePd of a
  toric orbit -- we presented some of these in the example of the
  toric orbit in \S\ref{sec:basic} -- but no such algorithm is known
  for Borel orbits. The fibration in Lemma \ref{bovert} suggests that,
  in our situation, one might be able to reduce this latter problem to
  the former. We will pursue this idea in a later publication.
\end{rem}

Lemma \ref{bovert} implies, in particular, that the codimension of
$B_d\erf$ is the sum of the codimensions of $\KK$ in $\Nh_d^0$ and
the codimension in the fiberwise directions. We collect the
appropriate numeric values in the following table: \newline \newline
\begin{tabular}{|c|c|c|c|c|c|c|}
\hline $d$ & $\dim\OO={d \choose 2}$ & $\dim \Nh_d$ & $\deg
\QQ_d=\codim (\OO)$ & $\dim(\KK)={d-1}$ & $\dim \Nh_d^0$ &
$\codim(\KK)$
\\
\hline 1 & 0 & 0 & 0 & 0 & 0 & 0\\
\hline 2 & 1 & 1 & 0 & 1 & 1 & 0\\
\hline 3 & 3 & 3 & 0 & 2 & 2 & 0\\
\hline 4 & 6 & 7 & 1 & 3 & 4 & 1\\
\hline 5 & 10 & 13 & 3 & 4 & 6 & 2\\
\hline 6 & 15 & 22 & 7 & 5 & 9 & 4
\\
\hline
\end{tabular}
\newline
\newline
The first 3 columns list the codimension of the closure of the Borel
orbit $\OO$ in $\Nh_d$, while the last three - the codimension of
the closure of the toric orbit $\KK$ in $\Nh_d^0$.

Now we are ready for the computations.

\subsection{The cases d=1,2,3}
In these cases $\deg\QQ_d=0$ and thus $\QQ_d=1$; geometrically, this
means that $\O_d=\eqnstil_d$, and thus $\OO_d=\Nh_d$. The case of
$d=1$ was described in \S\ref{sec:interlude}.

For $d=2$ we obtain
\begin{equation}
  \epd{\Thetab_2}=\res_{z_1=\infty}\res_{z_{2}=\infty}\frac{z_1-z_2}{2z_1-z_2}
  \ct\left(\frac1{z_1}\right)\ct\left(\frac1{z_2}\right)\,z_1^{k-n}\,z_2^{k-n}\;dz_1dz_2.
\end{equation}
Expanding the iterated residue, one immediately recovers Ronga's
formula \cite{sigma11}:
\begin{equation}\label{ronga}
\epd{\Thetab_2}=c_{k-n+1}^2+\sum_{i=1}^{k-n+1}2^{i-1}c_{k-n+1-i}c_{k-n+1+i}.
\end{equation}

For $d=3$, the formula is
\begin{multline}
\epd{\Thetab_3}=(-1)\res_{z_1=\infty}\res_{z_{2}=\infty}
\res_{z_{3}=\infty}\frac{(z_1-z_2)
(z_1-z_3)(z_2-z_3)}{(2z_1-z_2)(z_1+z_2-z_3)(2z_1-z_3)}\\
\ct\left(\frac1{z_1}\right)\ct\left(\frac1{z_2}\right)
\ct\left(\frac1{z_3}\right)
\,z_1^{k-n}\,z_2^{k-n}\,z_3^{k-n}\;dz_1dz_2dz_3.
\end{multline}
This is a more compact and conceptual formula for
$\epd{\Theta_3}$ than the one given in \cite{bfr}.

\subsection{The basic equations in
general}\label{sec:basicequations}

As our table in \S\ref{sec:codimension} shows, the polynomial $\qd
d$ is not trivial when $d>3$. As a step towards its computation, we
describe a set of equations satisfied by $\OO\subset \Nh_d$ and
$\KK\subset \Nh_d^0$. We will call these equations {\em basic}.

The equations will be written in terms of the coordinates
$\uh_{\tau|\ds}^{l}$ on $\NN_\ds$ introduced in
\eqref{npicoordinates}, where now we assume that $|\tau|=2$. Clearly,
these variables form a dual basis to the basis $\{q_{mr}^l\}$ of
$\Nh_d$. We will streamline our notation by writing $\uh_{mr}^l$
instead of $\uh_{[m,r]|\ds}^l$; naturally, we have
$\uh_{mr}^l=\uh_{rm}^l$, and $r+m\leq l$.

The construction is as follows.  If $i+j+m\leq l$, then the sequence
\[\bipi(i,j,m;l)=([1],[2]\ddd[l-1],[i,j,m],[l+1]\ddd[d-1],[d])\]
is admissible but not complete, hence $\bu^{\bipi(i,j,m;l)}$ will
appear as a term of some of the relations $\relz\birho\tau$
introduced in Proposition \ref{fundament}. In fact, it appears in
three different relations:
\[  \text{for } \tau=[i],\,\rho_l=[j,m],\;
 \text{for } \tau=[j],\,\rho_l=[i,m],\;
 \text{and for } \tau=[m],\,\rho_l=[i,j];
\]
in all cases $\rho_r=[r]$ for $r\neq l$.  Next, we reduce the
relation $\relz\birho\tau$ according to the prescription of Lemma
\ref{prescr}. After the reduction, only the terms corresponding to
the identity permutation and those corresponding to the
transpositions of the form $(s,l)$ survive; for example, in the case
$\tau=[m]$, we obtain the ``localized'' relation
\begin{equation}\label{uijkl}
\hat{u}_{ijm}^{l}=\sum_{s=j+m}^{l-i}\hat{u}_{jm}^{s}\hat{u}_{is}^{l}.
\end{equation}
Note that the number of terms on the right hand side is
$l-(i+j+m)+1$, which is the defect of $\hat{u}_{ijm}^{l}$ plus $1$.

We obtain two other expressions for $\hat{u}_{ijm}^{l}$ when we
choose $\tau$ to be $[j]$ or $[k]$, and the resulting equalities
provide us with quadratic relations among our variables
$\uh_{mr}^{l}$, $m+r\leq l\leq d$.
\begin{prop}\label{basicrelations} Let  $(i,j,m;l)$ be a quadruple of
  nonnegative integers satisfying  $i+j+m\le l\le d$. Then the ideal
  of the variety $\OO\subset \Nh_d$ contains the relations
  \begin{equation}
    \label{localrels} R(i,j,m;l):\quad
\sum_{s=j+m}^{l-i}\hat{u}_{jm}^{s}\hat{u}_{is}^{l}=\sum_{s=i+m}^{l-j}
\hat{u}_{im}^{s}\hat{u}_{js}^{l}=
\sum_{s=i+j}^{l-m}\hat{u}_{ij}^{s}\hat{u}_{ms}^{l}.
  \end{equation}
\end{prop}

\begin{rem}
\begin{itemize}
\item In general, the quadruple $(i,j,m;l)$ gives us 2
  relations. If $i=j\neq m$, then the number of relations reduces to
  1, and if $i=j=m$, then \eqref{localrels} is vacuous.
\item The equalities $R(i,j,m;l)$ with $i+j+m=l$ are relations of the toric
  orbit closure $\KK\subset \Nh_d^0$. We will call these equations
  {\em toric}.
\end{itemize}
\end{rem}

\subsection{d=4,5,6}\label{negyothat}

The first nontrivial case is $d=4$: here $\deg \qd 4=1$, i.e.
$\OO_4=\overline{B_4\erf}$ is a hypersurface in $\Nh_4$. Checking the
table at the end of \S~\ref{sec:codimension}, we see that the
codimension of the toric piece $\KK_4$ in $\Nh_4^0$ is the same as the
codimension of $\OO_4$ in $\Nh_4$. This means that
$\qd4=\epd{\KK_4,\Nh_4^0}$.

It is not surprising then to find that the only basic equation is a
toric one, corresponding to the quadruple $(1,1,2,4)$:
\begin{equation}\label{d4}
 R(1,1,2;4):\quad \uh_{11}^{2}\uh_{22}^{4}=\uh_{12}^{3}\uh_{13}^{4}.
\end{equation}
We note that this toric hypersurface is essentially our
example from \S\ref{sec:basic}. The variety defined by
\eqref{d4} in $\Nh_4$ is irreducible, and has the same dimension as
$\OO_4$, therefore it coincides with $\OO_4$. We have already
determined the equivariant Poincar\'e dual in this case in a number
of ways: it is the sum of the weights of any of the monomials in the
equation. This brings us to the formula
\begin{equation}
\QQ_4(z_1,z_2,z_3,z_4)=(2z_1-z_2)+(2z_2-z_4)=2z_1+z_2-z_4.
\end{equation}
As a result we obtain
\begin{multline*}
\epd{\Theta_4}=\res_{z_1=\infty}\res_{z_2=\infty}\res_{z_3=\infty}
\res_{z_{4}=\infty} \prod_{l=1}^4
\ct\left(\frac1{z_l}\right)\,z_l^{k-n}\;dz_l
\\
\frac{(z_1-z_2)
(z_1-z_3)(z_1-z_4)(z_2-z_3)(z_2-z_4)(z_3-z_4)(2z_1+z_2-z_4)}
{(2z_1-z_2)(z_1+z_2-z_3)(2z_1-z_3)
(z_1+z_3-z_4)(2z_2-z_4)(z_1+z_2-z_4)(2z_1-z_4)}.
\end{multline*}

\textbf{d=5}: Again, we consult our table. We have $\dim \Nh_5=13$
and $\codim\, \OO_5=3$, while $\dim \Nh_5^0=6$ and $\codim\,
\KK_5=2$.

Let us list our variables.
\begin{eqnarray*}
\text{6 toric}: &
\uh_{14}^{5},\uh_{23}^{5},\uh_{13}^{4},\uh_{22}^{4},\uh_{12}^{3},\uh_{11}^{2}  \\
\text{4 defect-1}: &\uh_{13}^{5},\uh_{22}^{5},\uh_{12}^{4},\uh_{11}^{3},  \\
\text{2 defect-2}: &\uh_{12}^{5},\uh_{11}^{4}, \text{and }  \\
\text{1 defect-3}: &\uh_{11}^{5}.
\end{eqnarray*}

There are 3 toric equations, which necessarily involve the toric
variables only:
\begin{eqnarray}
R(1,1,2;4):\quad\uh_{12}^{3}\uh_{13}^{4}=\uh_{11}^{2}\uh_{22}^{4} \notag\\
R(1,1,3;5):\quad \uh_{14}^{5}\uh_{13}^{4}=\uh_{23}^{5}\uh_{11}^{2} \label{d5a}\\
R(1,2,2;5):\quad\uh_{14}^{5}\uh_{22}^{4}=\uh_{23}^{5}\uh_{12}^{3}
\notag
\end{eqnarray}
and one defect-1 equation:
\begin{equation}\label{d5d}
 R(1,1,2;5):\quad
\uh_{13}^{5}\uh_{12}^{3}+\uh_{14}^{5}\uh_{12}^{4}=\uh_{11}^{2}\uh_{22}^{5}+\uh_{23}^{5}\uh_{11}^{3}
\end{equation}

We observe that the toric equations \eqref{d5a} describe the
vanishing of the $3$ maximal minors of a $2\times 3$ matrix. This is
an irreducible toric variety, thus we can again argue that it
coincides with $\KK_5$. Fortunately, this variety is a special case
of the $A_1$-singularity, this time with $n=2$ and $k=3$.
Substituting the appropriate weights into (\ref{sigma1}), we obtain:
\begin{multline}
\epd{\KK_5,\Nh_d^0}=\\
=\frac{(z_1+z_2-z_3)(2z_1-z_2)(z_1+z_4-z_5)-
(2z_2-z_4)(z_1+z_3-z_4)(z_2+z_3-z_5)}{z_1+z_4-z_2-z_3}=\\
=2z_1^2+3z_1z_2-2z_1z_5+2z_2z_3-z_2z_4-z_2z_5-z_3z_4+z_4z_5.
\end{multline}

Let $M_5$ denote the variety determined by the basic equations.
Notice that for fixed $\uh_{11}^2,\uh_{12}^3,\uh_{14}^5,\uh_{23}^5$
\eqref{d5d} is linear in the remaining variables. This means that
outside the codimension-2 subvariety $\KK_5'$ in $\KK_5$ where these 4
variables vanish, the natural projection $M_5 \to\KK_5$ is the
projection of a vector bundle onto its base, which implies that $M_5$
is irreducible, and thus $M_5=\OO_5$; the fibers of this vector bundle
are hyperplanes in the 7-dimensional complement of $\Nh_5^0$ in
$\Nh_5$. It is also clear from \eqref{d5d} that the variety determined
by the relation $R(1,1,2,5)$ is transversal to $\prz^{-1}(\KK_5)$
outside the part lying over $\KK_5'$, and hence we can conclude that
$\epd{\OO_5,\Nh_5}$ is the product of $\epd{\KK_5,\Nh_5^0}$ and the
weight of the relation $R(1,1,2;5)$. The latter equals $2z_1+z_2-z_5$,
hence the final result is
\[
\QQ_5(z_1,z_2,z_3,z_4,z_5)=
(2z_1+z_2-z_5)(2z_1^2+3z_1z_2-2z_1z_5+2z_2z_3-z_2z_4-z_2z_5-z_3z_4+z_4z_5).
\]

\textbf{d=6}

Now $\qd6$ is a degree-$7$ polynomial in 6 variables, and one needs
the help of a computer algebra program to do the calculations. Here
we summarize our computations with Macaulay.

Let $M_6$ denote, again, the variety defined by the basic equations.
It turns out, that the codimension of $M_6$ in $\Nh_6$ is equal to
the codimension of $\OO_6$, however, $M_6$ contains two maximal
dimensional components, namely,
\[M_6^1=\langle
\uh_{11}^2,\uh_{12}^3,\uh_{11}^3,\uh_{14}^5,\uh_{14}^6,\uh_{15}^6,\uh_{24}^6
\rangle \] and
\[M_6^2=\langle \mathrm{basic\ equations}, R \rangle,\] where the
extra relation is
\begin{multline*}
R=\uh_{12}^4\uh_{12}^4\uh_{23}^5\uh_{33}^6+\uh_{22}^4\uh_{13}^4\uh_{12}^5\uh_{33}^6
+\uh_{13}^4\uh_{13}^4\uh_{22}^5\uh_{23}^6+\uh_{22}^4\uh_{13}^4\uh_{23}^5\uh_{13}^6
\\
-\uh_{22}^4\uh_{11}^4\uh_{23}^5\uh_{33}^6-\uh_{13}^4\uh_{12}^4\uh_{22}^5\uh_{33}^6
-\uh_{22}^4\uh_{13}^4\uh_{13}^5\uh_{23}^6-\uh_{13}^4\uh_{13}^4\uh_{23}^5\uh_{22}^6=0
\end{multline*}
The weight of $R$ is $2z_1+3z_2+3z_3-2z_4-z_5-z_6$. Since $\OO_6$ is
irreducible, we have $\OO_6=M_6^2$. The other component, $M_6^1$, is
a linear subspace, and we obtain $\QQ_6$ as
\[\QQ_6=\epd{M_6}-\epd{M_6^1}.\]
Having described the vanishing ideal of $\OO_6$ by explicit relations,
using Macaulay, one then obtains $\QQ_6$; this formula is too long to
present here.

\subsection{An application: the positivity of Thom
  polynomials} \label{sec:positive} It is conjectured in
\cite[Conjecture 5.5]{rimanyi} that all coefficients of the Thom
polynomials $\tp d^{n\to k}$ expressed in terms of the relative
Chern classes are nonnegative. Rim\'anyi also proves that this
property is special to the $A_d$-singularities. In this final
paragraph, we would like to show that our formalism is well-suited
to approach this problem.  We will also formulate a more general
positivity conjecture, which will imply this statement.

We start with a comment about the sign $(-1)^d$ in our main formula
\eqref{thethom}. Recall from \eqref{defresinf} in
\S\ref{sec:pairings} that, according to our convention, the iterated
residue at infinity may be obtained by expanding the denominators in
terms of $z_i/z_j$ with $i<j$ and then {\em
  multiplying the result  by $(-1)^d$}. This sign appears because of
the change of orientation of the residue cycle when passing to the
point at infinity. This means that if we compute \eqref{thethom} via
expanding the denominators, then the sign in the formula cancels.

Now we are ready to formulate our positivity conjecture.
\smallskip

\noindent{\bf Conjecture}: {\em Expanding
  the rational function
\[ \frac{\prod_{m<l}(z_m-z_l)\,\QQ_d(z_1\ddd z_d)}
{\prod_{l=1}^d\prod_{m=1}^{l-1}\prod_{r=1}^{\min(m,l-m)}(z_m+z_r-z_l)}
\]
in the domain $|z_1|\ll\dots\ll|z_d|$, one obtains a Laurent series
with nonnegative coefficients.}

This statement clearly implies the nonnegativity of the
coefficients of the Thom polynomial.

At the moment we do not know how to prove this conjecture in
general. However, we observe that the expansion of a fraction of the
form $(1-f)/(1-(f+g))$ with $f$ and $g$ small has positive
coefficients. Indeed, this follows from the identity
\[   \frac{1-f}{1-f-g}=1+\frac{g}{1-f-g}.
\]
Now, introducing the variables $a=z_1/z_2$ and $b=z_2/z_3$, we can
rewrite the above fraction in the $d=3$ case as follows:
\[\frac{(z_1-z_2)
(z_1-z_3)(z_2-z_3)}{(2z_1-z_2)(z_1+z_2-z_3)(2z_1-z_3)}
=\frac{1-a}{1-2a}\cdot\frac{1-ab}{1-2ab}\cdot\frac{1-b}{1-b-ab}.
\]
Applying the above identity to the right hand side of this formula
immediately implies our conjecture for $d=3$. As a token reward for
having followed our paper this far, we offer to the reader the rather
amusing exercise of proving the same statement for $d=4$.

\newpage
\section{List of notations}
\label{sec:listnot}
\begin{itemize}
\item $\mathcal{J}(n)$: algebra of power series in $n$ variables,
  without constant term [\S\ref{sec:setup}].\\
  $\rd n$ : $d$-jets of holomorphic functions on $\C^n$ near
  the origin [\S\ref{sec:setup}]. \\
  $\mapd nk$: map-jets, i.e. {\em $d$-jets} of maps
  $(\C^n,0)\to(\C^k,0)$ [\S\ref{sec:setup}].

\item $\lin$: linear part of a germ or jet [\S\ref{sec:setup}].

\item $\diff n$: the group  of $d$-jets of diffeomorphisms of
  $\C^n$ fixing the origin [\S\ref{sec:setup}].

\item $A_\Psi$: the nilpotent algebra of the map germ $\Psi$
  [\S\ref{sec:setup}].\\ 
  $A_d$: the nilpotent algebra $t\C[t]/t^{d+1}$
  [\S\ref{sec:adsing}].

\item $\Theta_A,\Theta_A^{n\to k}$: set of jets with nilpotent
  algebra $A$
  [\eqref{defconea}].\\
  $\Theta_d,\Theta_d^{n\to k}$  notation for $\Theta_{A_d}$
  [\S\ref{sec:adsing}].

\item $\mathcal{K},\mathcal{K}_d(n,k)$: the contact group [\eqref{defk}].

\item $\epd{\Sigma,W}_T$: $T$-equivariant Poincar\'e dual  of $\Sigma
  \subset W$   [\S\ref{sec:epdmult}-\S\ref{subsec:axiomatic}].
\\
 $\euler^T(W)$: the equivariant Euler class of the $T$-module $W$ 
[\eqref{eec}].\\
 $\emu_p[M,Z]$: equivariant multiplicity of $M$ in $Z$ at $p\in M$ 
 [\eqref{emult}].

\item $\ct(q)$: the generating function of the relative Chern
classes [\eqref{deftc}].

\item $\tp A^{n\to k}(\boldsymbol{\lambda},\thetab)$: the Thom polynomial
of a nilpotent algebra $A$, [Definition \ref{deftp}].\\
$\tp d^{n\to k}$: the Thom polynomial of $A=A_d$.\\
$\DD^{j}_d$: the Thom-Damon polynomial  [Proposition
\ref{collect}].

\item $|\pi|, \um\pi$ , $\max(\pi)$ , $\comb(\pi)$: the length, the sum, the maximal element and the number of different
permutations of the partition $\pi$ [Notation
\ref{partition}].\\
$\Pi[m]$: the set of all partitions of $m$, [\eqref{pim}].

\item $\mapreg$: set of  curve-jets with nonvanishing linear
part [\eqref{mapreg}].\\
  $\Gammac$:  test curve $\mapd 1n$ [\eqref{gammaform}].\\
$\blp=(\blp^1,\ldots, \blp^d)=(A,B,C,\ldots)$: map-jet in $\mapd
nk$ [\eqref{psidef}].

\item $\qmr$: the  quotient $\jddiff$ [diagram
\eqref{diagups} and Proposition \ref{propmodel}].
\item $\grnk$:  the Grassmannian of   {\em codimension}-$dk$ linear
  subspaces in $\mapd nk$ [diagram \eqref{diagups}].

\item
  $\sol_\epsilon, \sol_\teps\subset \mapd nk$: the linear subspace of
  solutions of $\epsilon$ [Definition \ref{epsdef}, also \eqref{solepsilon}].\\
$\solfnk, \sol_{\eqnstil}$:  vector bundles with fibers $\sol_\teps$ and base
  $\fnk$ and $\eqnstil$ [\eqref{exact}].

\item $\hom^{\triangle}(\cdot,\cdot)$: filtration preserving linear
  maps  between two filtered vector spaces  [\eqref{filtered}].

\item $\phip$: the  map $\hom(\C_L^d,\C^n)
\longrightarrow\hom(\C_R^d,\symdot)$ defined in \eqref{homs}.

\item $\tfnkreg\subset\tfnk\subset\hofi$: [\eqref{fdn}--\eqref{fdnreg}].\\
$\fnk,\fnkreg$: quotients by $B_R$-action [Lemma \ref{fdnquot}].

\item $\epsilon$: element of $\tfnk$ thought of as a nonsingular system of
  linear equations.\\
  $\teps$: image of $\epsilon$ under the projection is $\tfnk\to
  \fnk$ [Definition \ref{epsdef}].

\item $V$: bundle over $\fnk$ and $\eqnstil$ associated to the standard
representation of $B_R$ [Lemma \ref{vident}].

\item $\hnd$: the maximal-rank elements of $\hm dn$ [\eqref{hnd}].

\item $\fld$: variety of full flags of $d$-dimensional subspaces of
  $\C^n$ [Lemma \ref{flagquot}].

\item $\ind X$: the induced space $\ind
X=\hnd\times_{B_L}X$ [Definition \ref{refsequence}].

\item $\newd$: the filtered subspace of $\symdot$ introduced in
  \eqref{defy}. 

\item $\homnewd$: the space of filtration-preserving maps with respect to the filtrations \eqref{defy} and
\eqref{filtcd}.

\item $\eqns$: the nondegenerate part of
$\homnewd$ [\eqref{defe}],\\
$\eqnstil$: the quotient $\eqns/B_R$ [Proposition
\ref{eqnstilquot}]; $\pre :\eqns \to \eqnstil$: the projection.
\item $\phigr$, $\phif$, $\phie$ and $\phi$: injective morphisms
  [\eqref{diagups},\eqref{defphit}, Proposition \ref{eqnstilquot}].
\item $\eref$: the sequence $(e_1\ddd e_d)\in\hnd$,\\
  $\reffl$: the corresponding flag in $\fld$ [Definition
  \ref{refsequence}], \\
  $\epsilon_{\rf}$: the reference system $\phip(\eref)$ in $\eqns$
  [\eqref{reference_eq}],\\
  $\buref=\pre(\epsilon_\rf)$: the corresponding point in $\eqnstil$
  [Definition \ref{epsdef}].

\item $\Bipi$ : the set of admissible sequences of partitions [Definition
\ref{admissible}],\\
$\defect(\bipi)$: integer defined for $\bipi\in\Bipi$ [Definition \ref{admissible}],\\
 $\Bipio$: the set of admissible sequences
corresponding to fixed points in $\O$  [\eqref{defpio}].

\item $\epsilon_\bipi$: the system $\eqns$ corresponding to the
  admissible sequence $\bipi$ [\eqref{ubipi}],\\
  $\buref=\pre(\epsilon_\bipi)\in \eqnstil$: the corresponding
  $T$-fixed point in $\eqnstil$.
\item $\O=\overline{B_d\buref}\subset\eqnstil$, the closure of the
Borel orbit of $\buref$ [diagram \eqref{diageqns}].

\item $\NN_\bipi$: the affine-linear subspace of $\eqns$ associated to
  $\bipi\subset\Bipi$ [Definition \ref{defnpi}].\\
$\O_\bipi$: the piece of the
  orbit closure $\O$ in the chart $\NN_\bipi$
  [\eqref{opdef}].

\item $u_\pi^l$:  coordinates on $\homnewd$ [\eqref{eqns}],
\\  $\defect(u^l_\pi)$: integer defined for $\um\pi\leq l$
  [\eqref{defdefect}],\\
$\uh_{\tau|\bipi}^l$: coordinates on $\NN_\bipi$ 
[\eqref{npicoordinates}].

\item $\bipi_\ds$: the distinguished sequence of partitions
  [\eqref{distinguished}],\\ 
 $\teps_\ds$, $\NN_\ds$, $\O_\ds$, etc.: simplified notation, replacing
 $\bipi_\ds$ by ``$\ds$'' in the indices.
\item $\Nh_d\subset \NN_\ds$: a linear subspace [\eqref{nhd}],\\
$\phat:\NN_\ds\to\Nh_d$: linear projection,\\
 $\OO_d=\phat(\O_\ds)\subset\Nh_d$ [\eqref{ohdef}],\\
$\uh_{mr}^l$: coordinates on $\Nh_d$ obtained as  the restriction of
$\uh_{[m,r]|\ds}^l$ [\S\ref{sec:basicequations}].
\item $\qfl$: the equivariant Poincar\'e dual of the fiber of our mode
  over $\reffl$ [\eqref{defqfl}],\\
 $Q_\bipi$: the equivariant Poincar\'e dual of $\O_\bipi$ in
 $\NN_\bipi$ [\eqref{qbipidef}],\\
 $Q_\ds$: simplified notation for the equivariant Poincar\'e dual
 of $\O_\ds$ in $\NN_\ds$,\\
 $\qd d$: The equivariant Poincar\'e dual of $\OO_d$ in $\Nh_d$
 [\eqref{ohdef}].
\item $\io$: the ideal of the subvariety $\O \subset
 \eqnstil$ [Definition \ref{defio}].
\item $\deg(p(\bz);S)$, $\coeff(L,z_l)$, $\lead(q(\bz);m)$:
 [\S\ref{sec:vanres} after
Lemma \ref{1lemma}].
\item $\polye$: polynomial functions on
  $\homnewd$ [\S\ref{sec:homog} before Lemma \ref{prescr}]. \\
  $\bu^\bipi$: a monomial in $\polye$ depending on $\bipi\in\Bipi$
  [\eqref{up}].\\ 
  $\Lambda$: subspace of $\polye$ [\eqref{deflambda}].\\
  $\relz{\birho, \tau}$: the relation \eqref{relz} in $\io$.

\end{itemize}

\end{document}